\documentclass[a4paper,10pt]{article}
\usepackage[utf8]{inputenc}
\usepackage[T1]{fontenc}
\usepackage{textcomp}
\usepackage{fancyhdr}
\usepackage{fullpage}
\usepackage{lmodern}
\usepackage{amsmath,amssymb,bm}
\usepackage{float}
\usepackage{graphicx}
\usepackage{caption}
\usepackage{subfig}
\usepackage{titling}
\IfFileExists{yhmath.sty}{\usepackage{yhmath}}{}
\providecommand{\wideparen}[1]{\overset{\frown}{#1}}
\usepackage{mathrsfs}
\usepackage{cite}
\usepackage{footnote}
\usepackage{amsthm}
\usepackage{color}
\usepackage{slashed}
\usepackage{verbatim}
\usepackage[colorlinks=true,linkcolor=red,citecolor=blue]{hyperref}
\usepackage{amsfonts,euscript}
\usepackage{latexsym, multicol}
\usepackage{epstopdf}
\usepackage{psfrag}
\usepackage{xcolor}
\usepackage{comment}
\usepackage{authblk}
\usepackage{indentfirst}
\usepackage{lipsum}
\usepackage{tikz}
\DeclareMathOperator{\err}{Err}
\newcommand{\widecheck}[1]{\mathpalette\widecheckaux{#1}}
\newcommand{\widecheckaux}[2]{%
  \begingroup
  \mathsurround=0pt
  \setbox0=\hbox{$#1#2$}%
  \mathord{\vbox{\offinterlineskip
    \hbox{\resizebox{\wd0}{0.3ex}{$#1\vee$}}%
    \kern0.1ex
    \hbox{\copy0}%
  }}%
  \endgroup
}
\numberwithin{equation}{section}

\allowdisplaybreaks[3]

\def\vr{\varrho}
\DeclareMathOperator{\Vol}{Vol}
\DeclareMathOperator{\Area}{Area}
\def\NNN{\mathbb{N}}
\def\SSS{\mathbb{S}}
\def\RRR{\mathbb{R}}

\def\W{{\mathbf{W}}}
\def\Wb{{\underline{\W}}}

\def\xja{\langle x \rangle}

\def\Des{\slashed{\De}}
\def\slg{\slashed{g}}

\def\sdivs{{\slashed{\sdiv}\,}}
\def\curls{{\slashed{\curl}\,}}
\def\kah{\widehat{\ka}}
\def\thh{\widehat{\th}}
\def\th{\theta}
\def\eps{\epsilon}

\def\nabs{\slashed{\nab}}

\def\sic{\wideparen{\si}}

\def\mk{\mathcal{K}}

\renewcommand{\c}{\cdot}

\newcommand{\pa}{\partial}
\newcommand{\pr}{\pa}

\newcommand{\una}{\underline{\alpha}}
\newcommand{\unb}{\underline{\beta}}
\newcommand{\bg}{\mathbf{g}}
\newcommand{\mo}{\mathcal{O}}
\newcommand{\unchi}{\underline{\chi}}
\newcommand{\omb}{\underline{\omega}}
\newcommand{\etab}{\underline{\eta}}
\newcommand{\sld}{\slashed{d}}

\def\ujp{{\langle u \rangle}}

\def\Vphi{\varPhi}
\def\vphi{\varphi}

\def\NN{\mathcal{N}}

\def\slu{{\slashed{U}}}

\def\sfr{\mathfrak{s}}
\def\sk{\sfr}

\newcommand{\mr}{\mathcal{R}}

\newcommand{\M}{\mathcal{M}}
\newcommand{\D}{\mathbf{D}}

\newcommand{\dd}{{\mathfrak{d}}}

\newcommand{\Ric}{{\ric}}
\DeclareMathOperator{\sRic}{Ric}
\newcommand{\g}{\bg}
\newcommand{\R}{{\mathbf{R}}}
\newcommand{\K}{\mathbf{K}}

\def\as{\slashed{\a}}
\def\hot{\widehat{\otimes}}

\def\II{\mathcal{I}}
\def\la{\lambda}
\def\rhoc{\wideparen{\rho}}

\def\fc{\widecheck{f}}

\def\MM{\mathcal{M}}

\def\chib{\unchi}

\def\om{\omega}
\def\ze{\zeta}
\def\Om{\Omega}
\def\aa{\una}
\def\bb{\unb}
\def\Si{\Sigma}
\def\si{\sigma}

\def\dk{\dd}
\def\dkb{\slashed{\dk}}

\def\Ga{\Gamma}
\def\xib{\underline{\xi}}
\def\a{\alpha}
\def\b{\beta}
\def\hch{\widehat{\chi}}
\def\hchb{\widehat{\chib}}
\def\trch{\tr\chi}
\def\trchb{\tr\chib}
\def\trchc{\widecheck{\trch}}

\def\ka{\kappa}

\def\de{\delta}
\def\De{\Delta}
\def\nab{\nabla}
\def\ov{\overline}

\def\Gab{\Ga_b}
\def\Gag{\Ga_g}

\def\ep{\varepsilon}
\def\les{\lesssim}

\DeclareMathOperator{\curl}{curl}

\DeclareMathOperator{\grad}{grad}

\DeclareMathOperator{\tr}{tr}
\DeclareMathOperator{\sdiv}{div}
\def\bdiv{\mathbf{Div}}
\def\ric{\mathbf{Ric}}

\newtheorem{thm}{Theorem}[section]
\newtheorem{prop}[thm]{Proposition}
\newtheorem{lem}[thm]{Lemma}
\newtheorem{cor}[thm]{Corollary}
\newtheorem{rk}[thm]{Remark}
\newtheorem{df}[thm]{Definition}

\title{Global stability of Minkowski spacetime with borderline decay}
\author{Dawei Shen\footnote{Email address: ds4350@columbia.edu \par\indent\hspace{0.235cm} Department of Mathematics, Columbia University, New York, NY, 10027.}}
\begin{document}
\maketitle
\begin{abstract}
We prove the global stability of Minkowski spacetime for initial data with borderline decay. The maximal globally hyperbolic development is future causally geodesically complete and admits a global maximal-null foliation in its future region. The main difficulty in the case of borderline decay is the logarithmic divergence of certain nonlinear error terms when estimated using pointwise decay alone. We overcome this difficulty by establishing weighted flux and spacetime bulk estimates for both curvature and connection coefficients, which allow the nonlinear interactions to be controlled without any loss of decay.
\end{abstract}
\tableofcontents
\section{Introduction}
\subsection{Einstein vacuum equations and the Cauchy problem}
A Lorentzian $4$--manifold $(\MM,\g)$ is called a vacuum spacetime if it satisfies the Einstein vacuum equations:
\begin{equation}\label{EVE}
\Ric(\g)=0\quad \mbox{ in }\MM,
\end{equation}
where $\Ric$ denotes the Ricci tensor of the Lorentzian metric $\g$. The Einstein vacuum equations are invariant under diffeomorphisms; therefore, one considers equivalence classes of solutions. Expressed in general coordinates, \eqref{EVE} is a non-linear geometric coupled system of partial differential equations of order 2 for $\g$. In the wave coordinates, it can be shown that \eqref{EVE} is hyperbolic and hence admits an initial value formulation.

The corresponding initial data for the Einstein vacuum equations are given by specifying a triplet $(\Si,g,k)$, where $(\Si,g)$ is a Riemannian $3$--manifold with metric $g$, and $k$ is a symmetric $2$--tensor on $\Si$ satisfying the constraint equations:
\begin{align}\label{constraintk}
    R=|k|^2-(\tr k)^2,\qquad\qquad\nab^jk_{ij}=\nab_i\tr k,
\end{align}
where $R$ denotes the scalar curvature of $g$, $\nab$ denotes the Levi-Civita connection of $g$ and
\begin{align*}
    |k|^2:=g^{im}g^{jl}k_{ij}k_{lm},\qquad\qquad\tr k:=g^{ij}k_{ij}.
\end{align*}
In the future development $(\MM,\g)$ of such initial data $(\Si,g,k)$, $\Si\subset \MM$ is a spacelike hypersurface with induced metric $g$ and  second fundamental form $k$.

The seminal well-posedness results for the Cauchy problem obtained in \cite{cb,cbg} ensure that for any smooth initial data $(\Si,g,k)$ satisfying \eqref{constraintk}, there exists a unique smooth maximal globally hyperbolic development $(\MM,\g)$ solution of the Einstein equations \eqref{EVE} such that $\Si\subset \MM$ and $g$, $k$ are respectively the first and second fundamental forms of $\Si$ in $\MM$.

The prime example of a vacuum spacetime is Minkowski space $(\MM,\eta)$:
\begin{equation*}
    \MM=\mathbb{R}^4,\qquad \boldsymbol{\eta}=-(dt)^2+(dx^1)^2 +(dx^2)^2+(dx^3)^2,
\end{equation*}
for which the initial data are given by
\begin{equation*}
    \Si=\mathbb{R}^3,\qquad e=(dx^1)^2+(dx^2)^2+(dx^3)^2,\qquad k=0.
\end{equation*}
A basic question for the Einstein vacuum equations is whether a small perturbation of the Euclidean initial data $(e,0)$ has a global future development that remains close to Minkowski spacetime. We first recall the main results on this stability problem.
\subsection{Previous works on the stability of Minkowski spacetime}
\label{ssec6.1}
In 1986, Friedrich \cite{Friedrich86} established future null geodesic completeness for small perturbations of Minkowskian hyperboloidal data. The global nonlinear stability of Minkowski spacetime for asymptotically flat Cauchy data was proved in 1993 in the seminal work of Christodoulou--Klainerman \cite{Ch-Kl}.

In 2003, Klainerman--Nicol\`o \cite{KlNi99,Kl-Ni} developed a double null proof of stability in the exterior of an outgoing null cone. Klainerman--Nicol\`o \cite{knpeeling} also proved curvature peeling in the exterior under stronger initial assumptions. Lindblad-Rodnianski \cite{lr1,lr2} gave a proof of global stability in wave coordinates using the weak null structure of the Einstein equations. Bieri \cite{Bieri} extended \cite{Ch-Kl} to weaker decay and regularity.

Huneau \cite{huneau} proved stability under a translational spacelike Killing symmetry. Hintz--Vasy \cite{hv} gave a proof of global stability using Melrose's $b$--analysis. Keir \cite{Keir18} gave another proof in harmonic coordinates using $r^p$--weighted energy estimates of Dafermos--Rodnianski \cite{Da-Ro}. Graf \cite{graf} proved global stability for spacelike-characteristic initial data.

In \cite{Shen22}, the author gave a double-null proof of exterior stability using $r^p$--weighted estimates. Hintz \cite{Hintz} proved exterior stability of Minkowski spacetime in a modified wave coordinate gauge using edge-$b$ energy estimates. In \cite{Shen23}, the author proved global stability under minimal decay assumptions using a maximal-null foliation, and subsequently established exterior stability at the borderline decay rate in \cite{Shen24}. Luk--Oh--Warnick \cite{LukOhWarnick} proved global stability in a centre-normalized Newman--Unti gauge.

For further results on the stability of Minkowski spacetime for the Einstein equations coupled to matter fields, we refer to the introductions of \cite{Shen22,Shen23} and the survey \cite{ShenSurvey26}.
\subsection{Initial data and comparison of stability results}
\label{ssec6.2}
\begin{df}\label{def6.1}
A vacuum initial data set $(\Si_0,g,k)$ is called maximal if
$\tr_g k=0$. In this case, the constraint equations
\eqref{constraintk} reduce to
\begin{equation}
R(g)=|k|_g^2,\qquad \sdiv_g k=0,\qquad \tr_g k=0.
\end{equation}
\end{df}
\begin{df}[Weighted Sobolev spaces]\label{dffunctional}
Fix smooth global coordinates $x=(x^1,x^2,x^3)$ on $\Si_0$. We denote $\xja=(1+|x|^2)^{1/2}$.
For $m\in\mathbb N_0$ and $\ell\in\RRR$, let $H_b^{m,\ell}(\Si_0)$ be the completion of $C_c^\infty(\Si_0)$
under the norm
\begin{equation}\label{weighted-b-sobolev}
\|F\|_{H_b^{m,\ell}(\Si_0)}^2:=\sum_{|\a|\leq m}\int_{\RRR^3}\xja^{2\ell+2|\a|}|\pr^\a F|^2\,dx.
\end{equation}
Let $W_b^{m,\infty,\ell}(\Si_0)$ be the space of locally integrable fields whose weak derivatives through order $m$ have finite norm
\begin{equation}\label{weighted-b-sobolev-infinity}
\|F\|_{W_b^{m,\infty,\ell}(\Si_0)}:=\sum_{|\a|\leq m}\big\|\xja^{\ell+|\a|}\pr^\a F\big\|_{L^\infty(\Si_0)}.
\end{equation}
For tensor fields, both norms are taken componentwise in these
coordinates.
\end{df}
Let $e$ denote the Euclidean metric and $g_M$ the spatial
Schwarzschild metric of mass $M\geq 0$, given in the asymptotic
region by
\begin{equation*}
g_M=\left(1-\frac{2M}{\varrho}\right)^{-1}d\varrho^2+\varrho^2d\sigma_{\SSS^2},\qquad\quad \varrho:=|x|.
\end{equation*}
For global norms, $g_M$ is understood to have a fixed smooth
interior extension, with $g_0=e$.
\begin{df}\label{defsq}
Let $s\in\RRR$ and $q\in\NNN_0$. A complete vacuum initial data set
$(\Si_0,g,k)$ is called \emph{$(s,q)$--asymptotically flat of size at most
$\ep_0$} if there exist global coordinates $x$ extending an asymptotic
coordinate system and a constant $M\geq0$, with $M=0$ when $s<3$,
such that the following holds. In these coordinates, set
\[
h=\begin{cases}
g-e, & s<3,\\
g-g_M, & s\geq3.
\end{cases}
\]
We require
\begin{equation}\label{initial-data-weighted-smallness}
\|h\|_{H_b^{q+1,\frac{s-4}{2}}(\Si_0)}
+\|k\|_{H_b^{q,\frac{s-2}{2}}(\Si_0)}
+\|h\|_{W_b^{q+1,\infty,\frac{s-1}{2}}(\Si_0)}
+\|k\|_{W_b^{q,\infty,\frac{s+1}{2}}(\Si_0)}
+M\leq\ep_0.
\end{equation}
Under the same coordinate and background conventions, the initial data
are called \emph{pointwise $(s,q)$--asymptotically flat of size at most
$\ep_0$} if
\begin{equation}\label{pointwiseasym}
\|h\|_{W_b^{q+1,\infty,\frac{s-1}{2}}(\Si_0)}+\|k\|_{W_b^{q,\infty,\frac{s+1}{2}}(\Si_0)}+M\leq\ep_0.
\end{equation}
\end{df}
\begin{rk}
The smallness condition \eqref{initial-data-weighted-smallness}
includes both weighted $L^2$ and pointwise bounds. More precisely,
it implies
\begin{equation}\label{sFlux}
\sum_{|\a|\leq q+1}\int_{\RRR^3}\xja^{s-4+2|\a|}|\pr^\a h|^2\,dx+\sum_{|\a|\leq q}\int_{\RRR^3}\xja^{s-2+2|\a|}|\pr^\a k|^2\,dx\leq\ep_0^2,
\end{equation}
while \eqref{pointwiseasym} gives
\[
h=O_{q+1}\left(\frac{\ep_0}{\xja^{(s-1)/2}}\right),\qquad
k=O_q\left(\frac{\ep_0}{\xja^{(s+1)/2}}\right).
\]
Here $F=O_m(f)$ means that $\xja^{|\a|}|\pr^\a F|\les f$ for all $|\a|\leq m$.\footnote{In particular, for any $s'<s$, pointwise $(s,q)$--asymptotically flat initial data are also $(s',q)$--asymptotically flat, with size at most $C_{s,s',q}\,\ep_0$.}
\end{rk}
We first recall the classical global stability theorem of Christodoulou--Klainerman \cite{Ch-Kl}.
\begin{thm}[Christodoulou--Klainerman \cite{Ch-Kl}]\label{ckmain}
Let $(\Si_0,g,k)$ be a maximal $(4,3)$--asymptotically flat initial data set of sufficiently small size in the sense of Definition~\ref{defsq}. Then its maximal globally hyperbolic development is geodesically complete and admits a global maximal foliation and an optical function in the exterior region. Moreover, the curvature decays to zero at infinity.
\end{thm}
The decay requirement was subsequently weakened in \cite{Shen23}, yielding global stability for all $s>1$ sufficiently close to the borderline threshold.
\begin{thm}[Global stability with minimal decay \cite{Shen23}]\label{Shen23main}
Let $s\in(1,2]$, and let $(\Si_0,g,k)$ be a maximal $(s,2)$--asymptotically flat vacuum initial data set of sufficiently small size in the sense of Definition~\ref{defsq}. Then the future development of $(\Si_0,g,k)$ exists globally.
\end{thm}
At the borderline threshold $s=1$, exterior stability was established in \cite{Shen24} under pointwise asymptotic flatness assumptions.
\begin{thm}[Exterior stability with borderline decay \cite{Shen24}]\label{Shen24main}
Let $(\Si_0,g,k)$ be a vacuum initial data set that is pointwise $(1,3)$--asymptotically flat of size at most $\ep_0$ on $\Si_0\setminus K$, where $K\subset\Si_0$ is compact and $\ep_0$ is sufficiently small. Then the future development of $\Si_0\setminus K$ exists globally and admits a double null foliation with complete outgoing null cones.
\end{thm}
In the spherically symmetric setting, stability in the interior region can be established at both the critical decay threshold $s=1$ and the critical regularity $q=\tfrac12$. In recent joint work with Wan \cite{ShenWanSpherical}, we proved the following result for the Einstein--scalar field system.
\begin{thm}[Interior stability in spherical symmetry \cite{ShenWanSpherical}]\label{ShenWanSphericalmain}
Consider the spherically symmetric Einstein--scalar field system with sufficiently small characteristic initial data that are $(1,\tfrac12)$--asymptotically flat in the spherically symmetric sense. Then the corresponding solution exists globally in the interior region and converges asymptotically to Minkowski spacetime.
\end{thm}
Here the fractional regularity $q=\tfrac12$ is understood in the Sobolev sense for characteristic initial data. The proof in \cite{ShenWanSpherical} relies essentially on the spherical reduction of the field equations. Remark~\ref{failure}\textup{(\ref{failure:spherical})} explains the obstruction to extending this argument to the nonsymmetric setting. For further stability results in spherical symmetry, we refer to the introduction of \cite{ShenWanSpherical}.
\begin{rk}
The following table compares the decay and regularity assumptions in selected stability results for the Einstein equations, without symmetry assumptions, including both vacuum and matter-coupled systems. The indices $s$ and $q$ represent the decay and regularity requirements, respectively, although the precise initial data norms may differ from those in Definition~\ref{defsq}. Here $0<\de_*\ll1$.
\begin{center}
\begin{tabular}{lccccc}
\hline
Authors & Reference & $s$ & $q$ & System & Region \\
\hline
Klainerman--Nicol\`o & \cite{knpeeling} & $7+\de_*$ & $5$ & Vacuum & Exterior \\
Christodoulou--Klainerman & \cite{Ch-Kl} & $4$ & $3$ & Vacuum & Global \\
Klainerman--Nicol\`o & \cite{Kl-Ni} & $4$ & $3$ & Vacuum & Exterior \\
Graf & \cite{graf} & $4$ & $3$ & Vacuum & Interior \\
Hintz--Vasy & \cite{hv} & $3+\de_*$ & $26$ & Vacuum & Global \\
Lindblad--Rodnianski & \cite{lr2} & $3+\de_*$ & $6$ & Vacuum & Global \\
Hintz & \cite{Hintz} & $(3,5)$ & $2433$ & Vacuum & Exterior \\
Shen & \cite{Shen22} & $(3,\infty)$ & $2$ & Vacuum & Exterior \\
Ionescu--Pausader & \cite{ionescu} & $3-\de_*$ & $200$ & Klein--Gordon & Global \\
Wang & \cite{WangEV} & $3-\de_*$ & $200$ & Massive Vlasov & Global \\
LeFloch--Ma & \cite{leflochMa} & $2+\de_*$ & $20$ & Klein--Gordon & Global \\
Bieri & \cite{Bieri} & $2$ & $2$ & Vacuum & Global \\
Luk--Oh--Warnick & \cite{LukOhWarnick} & $(1,\infty)$ & $29$ & Vacuum & Global \\
Shen & \cite{Shen23} & $(1,2]$ & $2$ & Vacuum & Global \\
Shen & \cite{Shen24} & $1$ & $3$ & Vacuum & Exterior \\
\multicolumn{2}{r}{Present paper} & $1$ & $2$ & Vacuum & Global \\
\hline
\end{tabular}
\end{center}
In the absence of symmetry, the approach of the present paper requires angular regularity to control the nonlinear terms in the commuted Bianchi and null structure equations. This is achieved through angular commutations and Sobolev inequalities on the spheres. We expect $q=\tfrac32$ to be the limiting regularity threshold for this approach, although the endpoint remains open.

For comparison, in terms of derivative count, $q=1$ corresponds to $H^2\times H^1$ regularity for the initial data $(g,k)$. At this level, the bounded $L^2$ curvature theorem of Klainerman--Rodnianski--Szeftel \cite{KRS} provides a local existence theory that could serve as a starting point for a global stability argument. More generally, Klainerman's strong stability conjecture \cite{Klainerman01} asks whether global stability of Minkowski spacetime holds under a scale-invariant smallness condition on the initial data, at the regularity level of $\dot H^{3/2}\times\dot H^{1/2}$, corresponding to $q=\tfrac12$. Establishing global stability for $\tfrac12\leq q<1$ would go beyond the bounded $L^2$ curvature theory and likely require substantially new ideas.
\end{rk}
\subsection{Main results}\label{ssec6.3}
We first state our main stability theorem in qualitative form. The precise estimates are given in Theorem~\ref{mainth}.
\begin{thm}[Global stability with borderline decay]\label{maintheoremintro}
Let $(\Si_0,g,k)$ be a maximal $(1,2)$--asymptotically flat vacuum
initial data set in the sense of Definition~\ref{defsq}, with
\begin{equation}\label{ini}
\|g-e\|_{H_b^{3,-3/2}}+\|k\|_{H_b^{2,-1/2}}+\|g-e\|_{W_b^{3,\infty,0}}+\|k\|_{W_b^{2,\infty,1}}\leq\ep_0.
\end{equation}
If $\ep_0$ is sufficiently small, then its maximal globally hyperbolic development $(\M,\g)$ is future causally geodesically complete and admits a global maximal-null foliation $(C_u,\Si_t)$.
\end{thm}
The next theorem shows that $s=1$ is the sharp decay threshold for
stability under the smallness condition \eqref{initial-data-weighted-smallness}.
\begin{thm}[Sharpness of the decay threshold]\label{sharp-decay-threshold}
For every $s<1$, $q\in\NNN_0$, and $\ep_0>0$, there exists a smooth, complete, asymptotically flat maximal vacuum initial data set $(\mathbb R^3,g,k)$ such that
\begin{equation*}
\|g-e\|_{H_b^{q+1,\frac{s-4}{2}}}+\|k\|_{H_b^{q,\frac{s-2}{2}}}+\|g-e\|_{W_b^{q+1,\infty,\frac{s-1}{2}}}+\|k\|_{W_b^{q,\infty,\frac{s+1}{2}}}<\ep_0,
\end{equation*}
while its maximal globally hyperbolic development is future null geodesically incomplete.
\end{thm}
\begin{proof}
The constructions in \cite{BeigOMurchadha96,LiYu15,ShenWan25,ShenWan26ADM,ChenKlainermanTS26} provide smooth, asymptotically flat vacuum initial data whose future developments contain trapped surfaces. We choose such a complete maximal data set $(g,k)$ on $\mathbb R^3$, as furnished by \cite{BeigOMurchadha96} or the time-symmetric examples in \cite{ChenKlainermanTS26}. After a compactly supported change of coordinates, we may assume $g(0)=e$. For $\lambda\geq1$, set
\[
g_\lambda(x)=g(x/\lambda),\qquad k_\lambda(x)=\lambda^{-1}k(x/\lambda).
\]
For $0<s<1$, rescaling and the decay of the fixed data give
\[
\|g_\lambda-e\|_{H_b^{q+1,\frac{s-4}{2}}}+\|k_\lambda\|_{H_b^{q,\frac{s-2}{2}}}+\|g_\lambda-e\|_{W_b^{q+1,\infty,\frac{s-1}{2}}}+\|k_\lambda\|_{W_b^{q,\infty,\frac{s+1}{2}}}\les_{s,q}\lambda^{\frac{s-1}{2}}\to0.
\]
The case $s\leq0$ follows by monotonicity of the weights. The rescaled data remain complete and maximal, and their future developments still contain trapped surfaces. Penrose's singularity theorem \cite{Penrose65}, applied to the rescaled data, gives future null geodesic incompleteness.
\end{proof}
\begin{rk}[Constructions of initial data]
The construction of asymptotically flat initial data with prescribed decay has been extensively studied. Corvino--Schoen \cite{CorvinoSchoen06} constructed vacuum data agreeing with Kerr data outside a compact region, while Carlotto--Schoen \cite{CarlottoSchoen16} obtained initial data localized to cones. Mao--Tao \cite{MaoTao22} developed explicit solution operators for the constraint equations, obtaining localized vacuum data with optimal decay and new support properties. LeFloch--Nguyen \cite{LeFlochNguyen23} developed a seed-to-solution method allowing prescribed asymptotic behavior. Fang--Szeftel--Touati \cite{FangSzeftelTouati25} constructed nontrivial vacuum initial data with arbitrarily prescribed decay beyond the mass term. In particular, their time-symmetric examples provide initial data in the range $1<s\leq2$ covered by Theorem~\ref{Shen23main}.

At the borderline decay threshold $s=1$, in joint work with Wan \cite{ShenWan26Data}, the author constructed nontrivial $(1,q)$--asymptotically flat vacuum initial data for any $q\in\NNN$ by adapting the conical solution operators of Mao--Tao \cite{MaoTao22}. These examples satisfy the corresponding weighted decay and regularity conditions in Theorem~\ref{maintheoremintro}. Independently, Chen--Klainerman \cite{ChenKlainerman26} constructed a large class of smooth vacuum initial data on exterior ends by integrating the constraint equations outward. Their method also yields nontrivial maximal $(1,q)$--asymptotically flat vacuum initial data satisfying the assumptions of Theorem~\ref{maintheoremintro}. Their construction also yields pointwise $(1,q)$--asymptotically flat data for arbitrarily large $q\in\NNN$, including genuinely borderline examples that fail to satisfy the weighted integrability condition \eqref{sFlux} with $s=1$ required for $(1,q)$--asymptotic flatness. Theorem~\ref{Shen24main} establishes exterior stability for such data, while the stability of the corresponding interior regions remains open.
\end{rk}
\subsection{Outline of the proof and main ingredients}\label{intro:proof-outline}
At the borderline decay rate $s=1$, pointwise decay alone gives a
logarithmically divergent bound for certain nonlinear curvature
errors at $p=1$. We retain this weight by combining curvature flux
and bulk estimates with strong sphere flux bounds, in particular
for $\eta$. After establishing the initial estimates, we
successively improve the curvature norm $\mr$, the maximal
connection norm $\mk$, and the null connection norm $\mo$, all
defined in Section~\ref{ssec8.1}.
\subsubsection{Curvature estimates}\label{sec:introR}
We use the $r^p$--weighted energy method of Dafermos--Rodnianski \cite{Da-Ro}, adapted by Holzegel to the Bianchi equations near null infinity \cite{holzegel}. For Minkowski stability, Keir applied the method to the reduced Einstein equations in wave coordinates \cite{Keir18}, and Luk--Oh--Warnick used it for the Teukolsky equation in Newman--Unti gauge \cite{LukOhWarnick}. In \cite{Shen22}, the method was applied directly to the Bianchi equations in the exterior. The $r^p$--weighted estimates were extended to the interior up to the axis $r=0$ in \cite{Shen23}, and the exterior case of borderline decay was treated in \cite{Shen24}.

We work at the critical weight $p=1$. To see how the positive bulk
arises, consider the linearized Bianchi pair $(\a,\b)$ on Minkowski
spacetime,
\[
 \nabs_3\a-\frac1r\a=-2\sld_2^*\b,\qquad\nabs_4\b+\frac4r\b=\sld_2\a.
\]
Pairing the equations with $2r\a$ and $4r\b$, respectively, cancels
the angular terms after integration by parts on each sphere. The
identities $e_3r=-1$, $e_4r=1$, $\bdiv e_3=-2/r$, and
$\bdiv e_4=2/r$ then give the current
$r|\a|^2e_3+2r|\b|^2e_4$ and the bulk $|\a|^2+10|\b|^2$.
On a region $V$ bounded to the future by $C_u$ and $\Si_t$,
the corresponding estimate for the nonlinear equations has the
schematic form
\[
\begin{aligned}
\int_{C_u}r|\a|^2\,d\mu_{C_u}+\int_{\Si_t}r\bigl(|\a|^2+|\b|^2\bigr)\,d\mu_g+\int_{V}\bigl(|\a|^2+|\b|^2\bigr)\,d\mu_{\g}\les\int_{\Si_0}r\bigl(|\a|^2+|\b|^2\bigr)\,d\mu_g+\text{NL},
\end{aligned}
\]
where $\mbox{NL}$ denotes the nonlinear errors. We combine the estimates for the Bianchi pairs with coefficients chosen so that the positive bulk of a component in one estimate absorbs its negative bulk in an adjacent estimate. Section~\ref{exteriorest} establishes positivity of the combined bulk at $p=1$.

The methods of \cite{Shen23,Shen24} do not close the global estimates at $p=s=1$, as explained in Remark~\ref{failure}. For the representative interaction $r^pr\nabs\rhoc\c\eta\c r\nabs\bb$ in the exterior region $\ujp\les r$, we combine the outgoing curvature flux of $r\nabs\rhoc$, the strong weighted sphere flux of $\eta$, and the spacetime bulk of $r\nabs\bb$. Coarea and Cauchy--Schwarz on each outgoing cone give
\begin{align*}
&\int_V r|r\nabs\rhoc|\,|\eta|\,|r\nabs\bb|\,d\mu_{\g}=\int du\int dt\int_{S(t,u)}\phi ar|r\nabs\rhoc|\,|\eta|\,|r\nabs\bb|\,d\slg\\
&\les\int du\left(\int dt\int_S\phi r|r\nabs\rhoc|^2\,d\slg\right)^{1/2}\left(\sup_t r^{1/2}\|a^{1/2}\eta\|_{L^\infty(S)}\right)\left(\int dt\int_S\phi a|r\nabs\bb|^2\,d\slg\right)^{1/2}\\
&\les\left(\sup_u\int_{C_u\cap V}r|r\nabs\rhoc|^2\,d\mu_{C_u}\right)^{1/2}\left(\int du\,\sup_t\|a^{1/2}\eta\|_{W^{1,4}(S)}^2\right)^{1/2}\left(\int_V|r\nabs\bb|^2\,d\mu_{\g}\right)^{1/2}\les\ep^3.
\end{align*}
The second inequality uses sphere Sobolev on $a^{1/2}\eta$ and Cauchy--Schwarz in $u$. The three factors are the outgoing curvature flux of $r\nabs\rhoc$, the strong weighted sphere flux of $\eta$, and the spacetime bulk of $r\nabs\bb$. The factor $a$ in the coarea measure is accounted for by $a^{1/2}$ in the sphere norm of $\eta$ and by the spacetime measure in the $\bb$ bulk. Thus the estimate closes at $p=1$ without pointwise bounds for $a$ or $r\eta$. Theorem~\ref{Junkman} uses the same combination of curvature energies and strong connection fluxes to control the nonlinear errors in the commuted curvature equations.
\begin{rk}
Mixed norms in which the supremum is taken before the transverse $L^2$ norm appear in \cite[Remarks~7--8]{Luk18}. In that setting, suitable weights and short null intervals provide the $L^1$ integrability needed to control the nonlinear errors. The context here is different: the evolution takes place in an unbounded region, where the short-interval gain from $L^2$ to $L^1$ is no longer available uniformly. A new ingredient of the present approach is to combine strong flux estimates with spacetime bulk estimates, providing the additional integrability required to close the nonlinear error estimates at the borderline decay $s=1$.
\end{rk}
\begin{rk}\label{failure}
We explain why the methods of \cite{Shen23,Shen24,ShenWanSpherical} do not extend directly to the global problem at $p=s=1$. Throughout this remark, we restrict $V$ to the exterior region $\ujp\les r$.
\begin{enumerate}
\item \emph{Global stability for $s>1$} \cite{Shen23}.
The stronger decay assumption $s\in(1,2]$ gives
\[
\|\rhoc\|_{L^4(S)}\les\frac{\ep}{r^{(s+2)/2}},\qquad \|\eta\|_{L^\infty(S)}\les\frac{\ep}{r\ujp^{(s-1)/2}},\qquad\|\bb\|_{L^4(S)}\les\frac{\ep}{r^{3/2}\ujp^{(s-1)/2}}.
\]
For $p=s$, these bounds give
\begin{align*}
\int_V r^p|\rhoc|\,|\eta|\,|\bb|\,d\mu_{\g}&\les \int du\,dr\,r^{p+1}\|\rhoc\|_{L^4(S)}\|\eta\|_{L^\infty(S)}\|\bb\|_{L^4(S)}\\
&\les \ep^3\int du\,dr\,\frac{r^{(s-5)/2}}{\ujp^{s-1}}\les \ep^3\int\ujp^{-(s+1)/2}\,du.
\end{align*}
The last integral converges for $s>1$ and diverges logarithmically at $s=1$.
\item \emph{Exterior stability at $s=1$} \cite{Shen24}.
For $0<\de\ll1$, the relevant bounds are
\[
\begin{alignedat}{2}
 \|\rhoc\|_{L^4(S)}&\les\frac{\ep}{r^{3/2}},
 &\qquad \|\eta\|_{L^\infty(S)}&\les\frac{\ep}{r^{1-\de}\ujp^\de},\\
 \|\bb\|_{L^4(S)}&\les\frac{\ep}{r^{3/2-\de}\ujp^\de},
 &\qquad \|\Om\|_{L^\infty(S)}&\les\left(\frac{r}{\ujp}\right)^{C\ep},
\end{alignedat}
\]
with $\Om$ the double null lapse function in \cite{Shen24}. At $p=1$, these bounds give a radial integrand proportional to
$(r/\ujp)^{2C\ep}r^{-2+2\de}\ujp^{-2\de}$.
Integrating over $r\gtrsim\ujp$ leaves the logarithmically divergent
integral $\int\ujp^{-1}\,du$.
Lowering the weight to $p=1-2\de$ restores integrability, giving
\begin{align*}
 \int_V r^{1-2\de}|\rhoc|\,|\eta|\,|\bb|\,d\mu_{\g}\les \ep^3\int du\,dr\,\frac{(r/\ujp)^{2C\ep}}{r^2\ujp^{2\de}}\les \ep^3\int\ujp^{-1-2\de}\,du
 \les_\de\ep^3.
\end{align*}
Suitable incoming transport equations recover the lost decay for some components at the cost of an additional derivative. For $\eta$, no suitable $\nabs_3$ equation is available for this step. Its bound retains the factor $(r/\ujp)^\de$, which grows like $r^\de$ in the wave zone $\ujp\simeq1$ and prevents the argument from closing the global estimates at the critical weight.
\item\label{failure:spherical}
\emph{Interior stability in spherical symmetry} \cite{ShenWanSpherical}.
Let $L$ be the envelope constructed from the initial data as in \cite[Lemma~3.1]{ShenWanSpherical}. The construction guarantees
\[
 \int_1^\infty \frac{dr}{rL(r)^2}<\infty.
\]
The integral $\int_1^\infty dr/(rL(r))$ need not be finite. In spherical symmetry, the geometric transport coefficient is quadratic in the scalar field and thus carries the integrable weight $1/(rL(r)^2)$. In the nonsymmetric problem, the outgoing transport equation for $\eta$ takes the schematic form
\[
 \nabs_4(r\eta)
 =-r\b+\Gag+r\Gag\c\Gag+\Gag\c(r\eta),
\]
with $\Gag$ as in Definition~\ref{gammag}. The desired bound $|\Gag|\les\ep/(rL(r))$ would make $r\Gag\c\Gag$ integrable along the cone. In the mixed term $\Gag\c(r\eta)$, the coefficient of $r\eta$ carries only one factor of $L^{-1}$. A direct Gronwall estimate would therefore require the integrability of $1/(rL(r))$ along the cone. The curvature source $-r\b$ poses a separate difficulty. The outgoing flux controls $\int r\|\b\|_{L^2(S)}^2\,dr$. Integrating the normalized sphere $L^2$ norm of $r\b$ along the cone amounts to controlling $\int\|\b\|_{L^2(S)}\,dr$, and Cauchy--Schwarz leaves the divergent factor $(\int_1^\infty dr/r)^{1/2}$.
\end{enumerate}
\end{rk}
\subsubsection{Maximal connection estimates}
We control the maximal connection coefficients by three-dimensional elliptic estimates on the maximal hypersurfaces. The main quantity is the normal contraction $i_Nk:=N\c k$ of the second fundamental form. Writing $(i_Nk)_N=\de$ and $(i_Nk)_A=\eps_A$, we have
\begin{equation}\label{intro:maximal-div-curl}
\begin{split}
 \sdiv(i_Nk)&=-r^{-1}\de+\text{NL},\\
 (\curl(i_Nk))_N&=\sic,\\
 (\curl(i_Nk))_A&=r^{-1}{\in_A}^{B}\eps_B
 +\tfrac12{}^*(\b-\bb)_A+\text{NL}.
\end{split}
\end{equation}
The linear $r^{-1}$ terms require a coercive estimate for $r^{-1}i_Nk$. We combine the radial Hardy inequality, whose constant is $\tfrac14$, with an angular spectral gap. The angular form of the maximal momentum constraint controls the $\ell=0$ modes of $i_Nk$ by curvature and nonlinear errors, allowing the angular estimate to contribute an additional $2$.

For the slice flux, set $F=r^{1/2}i_Nk$. The radial Hardy
inequality controls $\tfrac14\int_{\Si_t}r^{-2}|F|^2\,d\mu_g$
by the energy of $\nab_NF$ up to geometric errors. The angular
spectral gap gives an additional coefficient of $2$ from the
tangential derivatives. Combining the two estimates yields
\[
 \frac94\int_{\Si_t}r^{-2}|F|^2\,d\mu_g
 \leq \int_{\Si_t}|\nab F|^2\,d\mu_g
 +C\int_{\Si_t}r|(\b,\bb)|^2\,d\mu_g+\text{NL}.
\]
Combining the div--curl and radial Hardy identities and applying
Young's inequality, we obtain
\[
\begin{aligned}
 &\int_{\Si_t}\left(
 \sum_A|\nab_A F|^2
 +r\left|\nab_N(i_Nk)+r^{-1}i_Nk\right|^2
 \right)d\mu_g\\
 &\qquad\leq
 \de_0\int_{\Si_t}r^{-2}|F|^2\,d\mu_g
 +C_{\de_0}\int_{\Si_t}r|(\b,\sic,\bb)|^2\,d\mu_g
 +\text{NL}.
\end{aligned}
\]
Here the radial Hardy contribution with coefficient $\tfrac14$ cancels the corresponding contribution from squaring the leading linear div--curl terms. The geometric corrections are included in $\text{NL}$. The angular spectral gap then allows us to absorb the remaining $r^{-2}|F|^2$ term for sufficiently small $\de_0$.
\begin{rk}
The argument of \cite[Proposition~6.6]{Shen23} uses the radial Hardy constant $\tfrac14$. At $s=1$, the squared linear terms in the slice estimate reach this same coefficient, so the absorption argument cannot be extended directly to the endpoint. The angular estimate in Lemma~\ref{normal:angular-gap} raises the coercivity constant to $\tfrac94$, which allows us to absorb these terms.
\end{rk}
For the spacetime bulk, we apply the unweighted identities directly to $i_Nk$. The linear terms in \eqref{intro:maximal-div-curl} have coefficients of magnitude $r^{-1}$, so the corresponding estimates take the schematic form
\begin{align*}
\tfrac94\int_V r^{-2}|i_Nk|^2\,d\mu_{\g}&\leq \int_V|\nab(i_Nk)|^2\,d\mu_{\g}+C\int_V|(\b,\bb)|^2\,d\mu_{\g}+\text{NL},\\
\int_V|\nab(i_Nk)|^2\,d\mu_{\g}&\leq (1+\de_0)\int_Vr^{-2}|i_Nk|^2\,d\mu_{\g}+C_{\de_0}\int_V|(\b,\sic,\bb)|^2\,d\mu_{\g}+\text{NL}.
\end{align*}
Since $\tfrac94>1$, a sufficiently small $\de_0$ also permits absorption in the bulk. The curvature fluxes and bulks are controlled by the estimates established in Section~\ref{sec:introR}.
The estimates for $\nab\log\phi$ follow by the same div--curl and
Hardy argument, using $\curl(\nab\log\phi)=0$ in the angular estimate.
\subsubsection{Null connection estimates}
We control the null connection coefficients by combining weighted transport along the outgoing null cones with Hodge estimates on the spheres. A typical transport equation is
\[
\nabs_4X+\la\trch\,X=F,\qquad \la\geq\tfrac12.
\]
For $p\in[1,2]$, we pair this equation with $p r^{p-3}|X|^{p-2}X$. Along each cone, the sphere area form evolves by $e_4(d\slg)=\trch\,d\slg$, see Lemma \ref{dint}. Using also $e_4\log a=-2\om$ and $e_4(t)=\phi^{-1}$, we obtain
\begin{align*}
 &\frac{d}{dt'}\int_{S(t',u)}ar^{p-3}|X|^p\,d\slg+\int_{S(t',u)}\phi ar^{p-3}
 \left((p\la-1)\trch-(p-3)\frac{e_4r}{r}+2\om\right)|X|^p\,d\slg\\
 &=p\int_{S(t',u)}\phi ar^{p-3}|X|^{p-2}\langle X,F\rangle\,d\slg.
\end{align*}
The coefficient in parentheses is bounded below by $\bigl(p(2\la-1)+1\bigr)/(2r)>0$, so the bulk controls $\phi ar^{p-4}|X|^p$ on each cone. We integrate from the initial cut or from the vertex, where the sphere term vanishes by regularity. Take the supremum over terminal cuts for the sphere term and integrate up to the last cut for the bulk term. Integrating in $u$ and using the coarea identities gives
\begin{equation}\label{intro:null-transport-estimate}
\begin{split}
 &\int du\,\sup_{t'}\int_{S(t',u)}ar^{p-3}|X|^p\,d\slg
 +\int_Vr^{p-4}|X|^p\,d\mu_{\g}\\
 &\quad\les\int_{\Si_0\cap V}r^{p-3}|X|^p\,d\mu_g
 +\int_Vr^{p-3}|X|^{p-1}|F|\,d\mu_{\g}.
\end{split}
\end{equation}
The supremum in $t'$ is taken separately on each cone before integration in $u$, as required in the nonlinear curvature estimate above. We use \eqref{intro:null-transport-estimate} to derive bulk estimates for $r\mu$, $r\bb$, $\eta$, and $\hchb$. Applying localized trace estimates along the outgoing null cones to $a^{1/2}(\eta,\hchb)$ and $a^{1/2}r(\mu,\bb)$, together with the corresponding angular derivatives, we obtain
\[
\int du\,\sup_{t'}\left(\|a^{1/2}(\eta,\hchb)\|_{W^{1,4}(S(t',u))}^2+\|a^{1/2}r(\mu,\bb)\|_{L^4(S(t',u))}^2\right)\les\ep_0^2,
\]
see \eqref{null:incoming-shear-strong-W14}. These bounds improve the bootstrap norms used to control the nonlinear curvature terms.
\subsection{Structure of the paper}\label{intro:proof-structure}
Section~\ref{sec2} introduces the maximal-null foliation and records the geometric equations used throughout the paper. Section~\ref{secmainthm} states the main theorem and sets up the bootstrap argument, and Section~\ref{firstboot} establishes the preliminary estimates. Section~\ref{secR} proves the curvature energy estimates. Section~\ref{seck} controls the maximal connection coefficients using spatial div--curl estimates. Section~\ref{secO} controls the null connection coefficients using outgoing transport estimates and Hodge estimates on the spheres. Finally, Section~\ref{sec:conclusions} proves future causal geodesic completeness and asymptotic roundness at null infinity.
\subsection*{Acknowledgements and AI statement}
\noindent\textit{Acknowledgements.}
The author would like to thank Xuantao Chen, Sergiu Klainerman, J\'er\'emie Szeftel, and Jingbo Wan for many helpful discussions.

\medskip
\noindent\textit{Disclosure of AI use.}
The author used ChatGPT to assist with lengthy formula calculations and technical estimates, as well as to improve the language and presentation of the manuscript. The overall mathematical framework, main ideas, and proof strategies were developed by the author. All AI-assisted calculations and arguments were reviewed and verified by the author, who takes full responsibility for the correctness and content of the manuscript.
\section{Preliminaries}\label{sec2}
\subsection{Geometric set-up}
\subsubsection{Maximal foliation}\label{maximalfoliation}
We foliate spacetime $(\MM,\g)$ by spacelike hypersurfaces $\Si_t$ as level hypersurfaces of a time function $t$. We denote by $T$ the unit, future oriented, normal to $\Si_t$. We then introduce a coordinate system $(x^1,x^2,x^3)$ on $\Si_0$, and extend it to $\M$ by $T(x^i)=0$ for $i=1,2,3$. Relative to this foliation, the spacetime metric takes the form:
\begin{equation}\label{metricSi}
    \g=-\phi^2 dt^2+\sum_{i,j=1}^3g_{ij}dx^idx^j.
\end{equation}
We define the \emph{lapse function} of the foliation by
\begin{equation*}
    \phi(t,x):=(-\g(\D t,\D t))^{-\frac{1}{2}},
\end{equation*}
with $\D t$ the spacetime gradient of $t$. Then, we have $T=-\phi\D t$.
We define the second fundamental form:
\begin{align}\label{dfk}
    k_{ij}:=-\g(\D_i T,\pr_j).
\end{align}
We denote by $\nab$ and $\sRic$ respectively the induced covariant derivative and the corresponding Ricci curvature tensor on the leaves $\Si_t$. Relative to an orthonormal frame $e_1,e_2,e_3$ tangent to the leaves of the foliation, we have the formulae:
\begin{align*}
    \D_i e_j&=\nab_ie_j-k_{ij} T,&
    \D_i T&=-k_{ij}e_j,\\
    \D_T e_i&=\nab_Te_i+(\phi^{-1}\nab_i\phi)T,&
    \D_T T&=(\phi^{-1}\nab_i\phi)e_i,
\end{align*}
where $\nab_Te_i$ denotes the projection of $\D_T e_i$ to the tangent space of $\Si_t$.
\begin{prop}\label{maximal}
We impose a maximal foliation on $\Si_t$ in the sense of Definition \ref{def6.1}. Then, we have the maximal constraint equations
    \begin{align}\label{eqtrk=0}
        \tr k=0,\qquad\qquad\sdiv k=0,\qquad\qquad R=|k|^2,
    \end{align}
and the lapse equation:
\begin{align}\label{lapseequation}
     \De\phi=\phi|k|^2.
\end{align}
\end{prop}
\begin{proof}
See (1.0.11)--(1.0.13) in \cite{Ch-Kl}.
\end{proof}
\subsubsection{Maximal-null foliation}\label{nullfoliation}
Let $O$ be the origin of the initial coordinates of
Definition~\ref{defsq}. A change of coordinates supported in a fixed
ball puts $g_{ij}|_O=e_{ij}$ and preserves
\eqref{initial-data-weighted-smallness} up to a universal constant.
In constructing the initial cuts, we write $\varrho=|x|$ for the coordinate
radius on $\Si_0$.

Let $\Vphi$ be the integral curve of $T$ through $O$, called the
\emph{timelike axis}. We normalize the retarded optical function
along $\Vphi$ by
\[
    T(u)=1 \quad \text{on }\Vphi,\qquad u(O)=0.
\]
Extend the axis $\Vphi$ slightly to the past of $\Si_0$. Label each future null cone by the value of $u$ at its vertex on the axis. Every point $x\in\Si_0$ sufficiently close to $O$, with $x\neq O$, lies on a unique cone with a vertex in the past of $\Si_0$. Define $w_{ret}(x):=-u|_{\Si_0}$ and set $w_{ret}(O)=0$. In the normalized coordinates, $w_{ret}=\vr+O(\vr^2)$ near $O$.

Choose a smooth cutoff $\chi_*$ equal to $1$ on $[0,1]$ and vanishing on $[2,\infty)$. Take $\vr_*>0$ small enough that $w_{ret}$
is defined in a neighborhood of $\{\vr\leq 2\vr_*\}$, and set
\[
w=\vr+\chi_*(\vr/\vr_*)(w_{ret}-\vr),
\]
with $w=\vr$ outside the support of the cutoff. We prescribe
$u|_{\Si_0}=-w$. Thus the initial cuts agree with the local null
cone cuts for $\vr\leq \vr_*$ and with the coordinate spheres for
$\vr\geq 2\vr_*$. For $\vr_*$ sufficiently small, $\pr_\vr w>0$ away from $O$. Thus the positive level sets
\[
    S_{(0)}(w_1):=\{p\in\Si_0:w(p)=w_1\},\qquad w_1>0,
\]
are radial graphs over $\mathbb S^2$. We set $S_{(0)}(0):=\{O\}$. The cuts agree with the local
retarded cuts for $r\leq \varrho_*$ and with the coordinate spheres
for $\vr\geq 2\vr_*$.

In $D^+(\Si_0\setminus O)$, we solve the \emph{eikonal equation}
\begin{equation}\label{uSi0}
    \g^{\mu\nu}\pr_\mu u\pr_\nu u=0,\qquad
    u|_{\Si_0\setminus O}=-w,
\end{equation}
choosing the future outgoing branch
\begin{equation*}
    T(u)|_{\Si_0\setminus O}=|\nab w|>0.
\end{equation*}
For $u<0$, the corresponding outgoing null hypersurface $C_u$
emanates from $S_{(0)}(-u)$.
For each $p\in\Vphi\cap J^+(\Si_0)$, we construct the outgoing null cone $C_{u(p)}$ with vertex $p$. We extend $u$ to $I^+(O)$ by requiring these cones to be its level sets. In particular, $C_0$ includes its vertex $O$.

We make these constructions on regular maximal-null slabs, where the null generators have no caustics or intersections away from
their vertices. There $u$ is continuous and smooth away from
$\Vphi$. The vector field
\begin{equation*}
    L:=-\grad u\qquad\mbox{ on }\;\M\setminus\Vphi
\end{equation*}
is future directed and satisfies
\begin{equation*}
    \D_L L=0.
\end{equation*}
We write $u_c(t):=u(\Si_t\cap\Vphi)$ and denote its inverse on its
range by $t_c$. The normalization along $\Vphi$ gives
\begin{equation*}
    u_c'(t)=\phi|_{\Si_t\cap\Vphi}>0,\qquad u_c(0)=0.
\end{equation*}
The nondegenerate leaves of the foliation are
\begin{equation*}
    S(t,u):=\Si_t\cap C_u,\qquad -\infty<u<u_c(t),\quad t\geq0.
\end{equation*}
We regard the points on the axis,
\begin{equation}\label{dfucttcu}
    S(t,u_c(t)):=\Si_t\cap\Vphi,\qquad
    S(t_c(u),u):=C_u\cap\Vphi,
\end{equation}
as spheres of radius $0$, with the second expression used for $u$
in the range of $u_c$.

\begin{figure}[htbp]
\centering
\begin{tikzpicture}[
    x=0.9cm,y=0.75cm,
    font=\small,
    slice/.style={black,line width=0.85pt},
    cone/.style={blue,line width=0.85pt},
    sphere/.style={draw=black,fill=red,line width=0.25pt}
]
    \path[use as bounding box] (-8.1,-0.6) rectangle (8.1,8.1);
    \coordinate (O)  at (0,0);
    \coordinate (P)  at (0,4.15);
    \coordinate (A)  at (-3.2,0);
    \coordinate (B)  at (3.5,0);
    \coordinate (SL) at (-2.4,6.55);
    \coordinate (SC) at (0,6.55);
    \coordinate (SR) at (2.4,6.55);
    % Maximal slices and the timelike axis.
    \draw[slice] (-8,0) -- (8,0);
    \draw[slice] (-8,6.55) -- (8,6.55);
    \draw[slice] (O) -- (0,8);
    % Outgoing null hypersurfaces.
    \begin{scope}
        \clip (-8,0) rectangle (8,8);
        \draw[cone] (-9,9) -- (O) -- (9,9);
        \draw[cone] (-9,13.15) -- (P) -- (9,13.15);
        \draw[cone] (A) -- (-12.2,9);
        \draw[cone] (B) -- (12.5,9);
    \end{scope}
    % Spheres represented by red points.
    \foreach \p in {O,P,A,B,SL,SC,SR}
        \filldraw[sphere] (\p) circle (1.6pt);
    % Geometric labels.
    \node[below=4pt] at (O) {$O$};
    \node[anchor=north west] at (0.15,7.95) {$\Vphi$};
    \node at (-1.85,7.65) {$I^+(O)$};
    \node at (-4.25,2.5) {$D^+(\Si_0\setminus O)$};
    \node[anchor=north] at (6.8,-0.12) {$\Si_0$};
    \node[anchor=north] at (7.25,6.3) {$\Si_t$};
    % Null-cone labels.
    \node[blue,anchor=west] at (1.65,1.35) {$C_0$};
    \node[blue,anchor=west] at (1.3,5.1) {$C_u$};
    \node[blue,anchor=east] at (-4.75,1.1) {$C_u$};
    \node[blue,anchor=west] at (5.25,1.15) {$C_u$};
    % Spherical leaves and axis points.
    \node[red,anchor=north east] at (-2.5,6.4) {$S(t,u)$};
    \node[red,anchor=north west] at (2.5,6.4) {$S(t,u)$};
    \node[red,anchor=west] at (0.18,7.05) {$S(t,u_c(t))$};
    \node[red,anchor=north west] at (0.12,4.1) {$S(t_c(u),u)$};
\end{tikzpicture}
\caption{Maximal-null foliation of the spacetime $\M$.}
\label{foliation}
\end{figure}
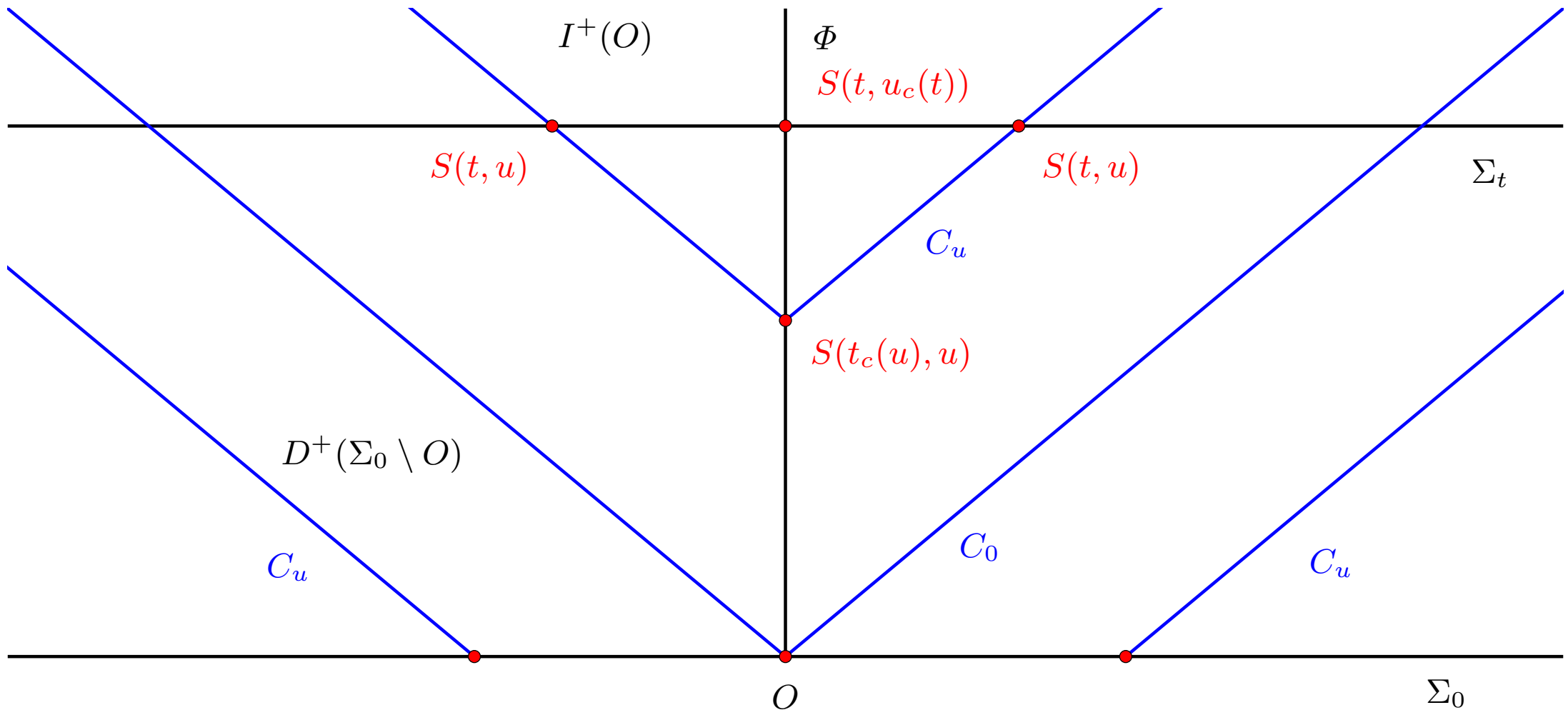

The future development is decomposed as
\begin{equation*}
    \M\cap J^+(\Si_0)=I^+(O)\cup C_0\cup D^+(\Si_0\setminus O),
\end{equation*}
where $I^+(O)$ is the chronological future of $O$ and $D^+$ denotes
the future domain of dependence. These regions satisfy
\begin{equation*}
u>0\ \mbox{ in }\;I^+(O),\qquad u=0\ \mbox{ on }\;C_0,\qquad u<0\ \mbox{ in }\;D^+(\Si_0\setminus O).
\end{equation*}
For every $t\geq0$, we define
\begin{equation}\label{dfVt}
    V_t:=J^-(\Si_t)\cap J^+(\Si_0)=\bigcup_{0\leq t'\leq t}\Si_{t'}.
\end{equation}
Here $J^\pm$ denote the causal future and past, respectively.
See Figure~\ref{foliation} for the maximal-null foliation.

On $\M\setminus\Vphi$, we define the \emph{null lapse function} $a$
and the outgoing null vectorfield $e_4$ by
\begin{equation*}
    a:=|\nab u|^{-1},\qquad e_4:=aL.
\end{equation*}
Since
\begin{equation*}
    0=\g(\D u,\D u)=|\nab u|^2-|T(u)|^2=a^{-2}-|T(u)|^2
\end{equation*}
and $T(u)>0$, we have $a=T(u)^{-1}$. The vertex normalization gives
a continuous extension of $a$ to $\Vphi$ with $a|_{\Vphi}=1$.
Finally, we set
\begin{equation*}
    N:=e_4-T.
\end{equation*}
Then
\begin{equation*}
    \g(T,N)=-aT(u)+1=0,\qquad
    \g(N,N)=2aT(u)-1=1,\qquad N=-a\nab u.
\end{equation*}
Thus, $N$ is the outward unit normal to $S(t,u)$ within $\Si_t$. We denote by $\th$ the second fundamental form of the surfaces $S(t,u)$ relative to $\Si_t$:
\begin{equation*}
    \th_{AB}:=\g(\D_A N,e_B).
\end{equation*}
We then define the null vectorfield $e_3$ by
\begin{equation*}
    e_3:=T-N.
\end{equation*}
Since $e_4=T+N$, we have
\begin{align*}
    \g(e_3,e_3)=\g(T-N,T-N)=0,
    \qquad
    \g(e_3,e_4)=\g(T-N,T+N)=-2.
\end{align*}
On a given nondegenerate $2$--sphere $S(t,u)$, we choose a local orthonormal frame $(e_1,e_2)$ and call $(e_1,e_2,e_3,e_4)$ a null frame.\footnote{The null frame is defined locally away from $\Vphi$.}
\begin{rk}
As a convention, throughout the paper, we use capital Latin letters $A,B,C,...$ to denote an index from $1$ to $2$, lowercase Latin letters $i,j,k,\ldots$ for indices from $1$ to $3$, and Greek letters $\alpha,\beta,\gamma,\ldots$ for spacetime indices, e.g., $e_A$ denotes either $e_1$ or $e_2$.
\end{rk}
The spacetime metric $\bg$ induces a Riemannian metric $\slg$ on $S(t,u)$. We use $\nabs$ to denote the Levi-Civita connection of $\slg$ on $S(t,u)$. We recall the null decomposition of the Ricci coefficients and curvature components of the null frame $(e_1,e_2,e_3,e_4)$ as follows:
\begin{equation}\label{defga}
\begin{alignedat}{3}
 \chib_{AB}&=\g(\D_A e_3,e_B),
 &\qquad \chi_{AB}&=\g(\D_A e_4,e_B),
 &\qquad \ze_A&=\tfrac12\g(\D_{e_A}e_4,e_3),\\
 \xib_A&=\frac12\g(\D_3 e_3,e_A),
 &\qquad \xi_A&=\frac12\g(\D_4 e_4,e_A),
 &\qquad \omb&=\frac14\g(\D_3 e_3,e_4),\\
 \etab_A&=\frac12\g(\D_4 e_3,e_A),
 &\qquad \eta_A&=\frac12\g(\D_3 e_4,e_A),
 &\qquad \om&=\frac14\g(\D_4 e_4,e_3),
\end{alignedat}
\end{equation}
and
\begin{equation}\label{defr}
\begin{alignedat}{3}
 \a_{AB}&=\R(e_A,e_4,e_B,e_4),
 &\qquad \b_A&=\frac12\R(e_A,e_4,e_3,e_4),
 &\qquad \rho&=\frac14\R(e_3,e_4,e_3,e_4),\\
 \aa_{AB}&=\R(e_A,e_3,e_B,e_3),
 &\qquad \bb_A&=\frac12\R(e_A,e_3,e_3,e_4),
 &\qquad \si&=\frac14{}^*\R(e_3,e_4,e_3,e_4),
\end{alignedat}
\end{equation}
where ${}^*\R$ denotes the Hodge dual of $\R$.

The second null fundamental forms $\chi, \chib$ are further decomposed in their traces $\trch$ and $\trchb$, and traceless parts $\hch$ and $\hchb$:\footnote{Throughout this paper, we use the same notation $\tr$ for both the $2$--dimensional and $3$--dimensional traces, with no confusion, as the underlying tensors are defined on $S$ and $\Si$, respectively.}
\begin{align*}
\trch&:=\slg^{AB}\chi_{AB},\qquad\quad \,\hch_{AB}:=\chi_{AB}-\frac{1}{2}\slg_{AB}\trch,\\
\trchb&:=\slg^{AB}\chib_{AB},\qquad\quad \, \hchb_{AB}:=\chib_{AB}-\frac{1}{2}\slg_{AB}\trchb.
\end{align*}
We define the horizontal covariant operator $\nabs$ as follows:
\begin{equation*}
\nabs_X Y:=\D_X Y-\frac{1}{2}\chib(X,Y)e_4-\frac{1}{2}\chi(X,Y)e_3.
\end{equation*}
We also define $\nabs_4X$ and $\nabs_3X$ as the horizontal projections of $\D_4X$ and $\D_3 X$:
\begin{align*}
\nabs_4X&:=\D_4X-\frac{1}{2}\g(X,\D_4e_3)e_4-\frac{1}{2}\g(X,\D_4e_4)e_3,\\
\nabs_3X&:=\D_3X-\frac{1}{2}\g(X,\D_3e_3)e_4-\frac{1}{2}\g(X,\D_3e_4)e_3.
\end{align*}
A spacetime tensor field $\psi$ is $S$--tangent if each contraction with $e_3$ or $e_4$ vanishes. We write $\nabs_3\psi$ and $\nabs_4\psi$ for the projections of $\D_3\psi$ and $\D_4\psi$ onto $S(t,u)$. Along $S(t,u)$, we decompose $k$ into its normal, mixed, and tangential components by setting
\begin{equation}\label{boundarydec}
    \de:=k_{NN},\qquad \eps_A:=k_{AN},\qquad \ka_{AB}:=k_{AB}.
\end{equation}
The maximal condition $\tr_g k=0$ gives $\tr_S\ka=-\de$. We denote the trace-free part of $\ka$ with respect to $\slg$ by
\[
\kah_{AB}:=\ka_{AB}-\frac12(\tr_S\ka)\slg_{AB}=\ka_{AB}+\frac12\de\,\slg_{AB}.
\]
The following identities hold for a maximal-null foliation:
\begin{equation}\label{6.6}
\begin{alignedat}{3}
 \chi_{AB}&=\th_{AB}-\ka_{AB},
 &\qquad \chib_{AB}&=-\th_{AB}-\ka_{AB},
 &\qquad \ze_A&=\eps_A,\\
 \eta_A&=\nabs_A\log a+\eps_A,
 &\qquad \etab_A&=\nabs_A\log\phi-\eps_A,
 &\qquad \xib_A&=\nabs_A\log\phi-\nabs_A\log a,\\
 2\om&=-\nab_N\log\phi+\de,
 &\qquad 2\omb&=\nab_N\log\phi+\de,
 &\qquad \xi_A&=0.
\end{alignedat}
\end{equation}
see (7.5.2b) in \cite{Ch-Kl}.
\subsection{Integral formulae}\label{ssecave}
We define the $S$-average of scalar functions.
\begin{df}\label{average}
Given any scalar function $f$, we denote its average and its mean-free part by
\begin{equation*}
    \ov{f}:=\frac{1}{|S|}\int_{S}f\,d\slg,\qquad\quad \fc:=f-\overline{f},
\end{equation*}
where $|S|$ denotes the area of $S$.
\end{df}
We recall the following integral formulae.
\begin{lem}\label{dint}
For any scalar function $f$, the following equations hold:
\begin{align*}
    \phi e_4\left(\int_{S(t,u)} f d\slg\right) &=\int_{S(t,u)} \left(\phi e_4(f) + \phi\trch f\right) d\slg,\\
    aN\left(\int_{S(t,u)} f d\slg\right) &= \int_{S(t,u)} \left(aN(f) + a\tr\th f\right) d\slg.
\end{align*}
Taking $f=1$, we obtain
\begin{equation}
    e_4r=\frac{\overline{\phi\tr\chi}}{2\phi}r,\qquad N(r)=\frac{\ov{a\tr\th}}{2a}r, \label{e3e4r}
\end{equation}
where $r$ is the \emph{area radius} defined by $4\pi r(t,u)^2=|S(t,u)|$.
\end{lem}
\begin{proof}
See Lemma 2.26 in \cite{kl-ro}.
\end{proof}
\subsection{Hodge systems in 2D--geometry}\label{ssec7.2}
\begin{df}
On an oriented Riemannian two-sphere $(S,\slg)$, let $\sfr_0(S)$,
$\sfr_1(S)$, and $\sfr_2(S)$ denote the spaces of pairs of scalar
functions, one-forms, and symmetric trace-free covariant two-tensors,
respectively. We omit $(S)$ when the sphere is clear from the context
and write $\in_{AB}$ for the area form of $\slg$. For $\xi,\eta\in\sfr_1$ and $U,V\in\sfr_2$, define the Hodge duals by
\[
 ({}^*\xi)_A:=\in_{AB}\xi^B,\qquad
 ({}^*U)_{AB}:=\in_{AC}U^C{}_B.
\]
The duals lie in $\sfr_1$ and $\sfr_2$, respectively, and satisfy
${}^*({}^*\xi)=-\xi$ and ${}^*({}^*U)=-U$.
We use the products
\begin{align*}
\xi\cdot\eta&:=\slg^{AB}\xi_A\eta_B,&\qquad \xi\wedge\eta&:=\in^{AB}\xi_A\eta_B,\\
(\xi\cdot U)_A&:=\xi^B U_{AB},&\qquad U\wedge V&:=\in^{AB}U_{AC}V^C{}_B,\\
(\xi\hot\eta)_{AB}&:=\xi_A\eta_B+\xi_B\eta_A-\slg_{AB}(\xi\cdot\eta).
\end{align*}
For $\xi\in\sfr_1$, set
\[
 \sdivs\xi:=\nabs^A\xi_A,\qquad
 \curls\xi:=\in^{AB}\nabs_A\xi_B,\qquad
 (\nabs\hot\xi)_{AB}
 :=\nabs_A\xi_B+\nabs_B\xi_A-\slg_{AB}\sdivs\xi.
\]
The Hodge operators $\sld_1:\sfr_1\to\sfr_0$ and
$\sld_2:\sfr_2\to\sfr_1$ and their formal $L^2$ adjoints
are given by
\begin{align*}
 \sld_1\xi&:=(\sdivs\xi,\curls\xi),
 & (\sld_2 U)_A&:=\nabs^B U_{AB},\\
 \bigl(\sld_1^*(f,f_*)\bigr)_A&:= -\nabs_A f+\in_A{}^B\nabs_B f_*,&\sld_2^*\xi&:=-\frac{1}{2}\nabs\hot\xi,
\end{align*}
where $(f,f_*)\in\sfr_0$. See \cite[Section~2.2]{Ch-Kl}.
\end{df}
\begin{prop}\label{Hodge2Diden}
We have the following identities:
\begin{align}
    \begin{split}\label{dddd}
        \sld_1^*\sld_1&=-\Des_1+\mathbf{K},\qquad\qquad  \sld_1\sld_1^*=-\Des_{{0}},\\
        \sld_2^*\sld_2&=-\frac{1}{2}\Des_2+\mathbf{K},\qquad\quad\; \sld_2\sld_2^*=-\frac{1}{2}(\Des_1+\mathbf{K}),
    \end{split}
\end{align}  
where $\mathbf{K}$ denotes the Gauss curvature on $S$ and $\Des_k=\nabs^A\nabs_A$ denotes the rough Laplacian on $\sfr_k$, acting componentwise when $k=0$.
\end{prop}
\begin{proof}
    See (2.2.2) in \cite{Ch-Kl}. 
\end{proof}
\begin{df}\label{dfdkb}
We define the weighted angular derivatives $\dkb$ as follows:
\begin{align*}
    \dkb U &:= r\sld_2 U,\qquad\quad\;\;\quad\quad\,\, \forall U\in \sk_2,\\
    \dkb \xi&:= r\sld_1 \xi,\qquad\qquad\quad\quad\,\, \forall \xi\in \sk_1,\\
    \dkb (f,f_*)&:= r\sld_1^*(f,f_*),\qquad\forall (f,f_*)\in \sk_0.
\end{align*}
\end{df}
\begin{prop}\label{commdkb}
We have the following commutator identities:
\begin{align*}
\dkb\sld_2U
&=\sld_1\dkb U,
&&\forall U\in\sk_2,\\
\dkb\sld_1\xi
&=\sld_1^*\dkb\xi,
&&\forall \xi\in\sk_1,\\
\dkb(2\sld_2^*\xi)
&=\sld_1^*\dkb\xi-2r\mathbf{K}\xi,
&&\forall \xi\in\sk_1,\\
\dkb\sld_1^*(f,f_*)
&=\sld_1\dkb(f,f_*),
&&\forall (f,f_*)\in\sk_0.
\end{align*}
\end{prop}
\begin{proof}
We recall from \eqref{dddd} that for any $\xi\in\sk_1$
\begin{align*}
    2\dkb \sld_2^* \xi =2r\sld_2\sld_2^*\xi=-r(\Des_1+\mathbf{K})\xi=r(\sld_1^*\sld_1-2\mathbf{K})\xi=\sld_1^*\dkb\xi-2r\mathbf{K}\xi,
\end{align*}
which implies the third identity. The other identities follow directly from Definition \ref{dfdkb}. This concludes the proof of Proposition \ref{commdkb}.
\end{proof}
\begin{df}\label{Lpnorms}
For a tensor field $f$ on a $2$--sphere $(S,\slg)$, we define
\begin{align*}
    \|f\|_{L^p(S)}:=\left(\int_S |f|^p\,d\slg\right)^{1/p},
    \qquad 1\leq p<\infty,\qquad
    \|f\|_{L^\infty(S)}:=\operatorname*{ess\,sup}_{x\in S}|f(x)|.
\end{align*}
For any integer $k\geq0$, we define the weighted Sobolev norms by
\begin{align*}
    \|f\|_{W^{k,p}(S)}&:=\left(\sum_{j=0}^k\|(r\nabs)^jf\|_{L^p(S)}^p\right)^{1/p},
    & 1&\leq p<\infty,\\
    \|f\|_{W^{k,\infty}(S)}&:=\sum_{j=0}^k\|(r\nabs)^jf\|_{L^\infty(S)},&\|f\|_{H^k(S)}&:=\|f\|_{W^{k,2}(S)}.
\end{align*}
Here all pointwise tensor norms are taken with respect to $\slg$.
\end{df}
\subsection{Null structure equations}
We choose the spacetime and sphere orientations
compatibly so that
\begin{equation*}
    \R(e_A,e_B,e_3,e_4)=2\si\in_{AB}.
\end{equation*}
All contractions and Hodge duals of horizontal tensors are taken
with respect to $\slg$ and its orientation. In the maximal-null foliation constructed above, $\xi=0$.
The null structure equations therefore take the following form;
see Proposition 7.4.1 in \cite{Ch-Kl}.
\begin{prop}\label{nulles}
In an Einstein vacuum spacetime with the maximal-null foliation
defined above, the null second fundamental forms satisfy
\begin{align*}
\nabs_4\hch+\trch\,\hch
    &=-2\om\hch-\a,\\
\nabs_4\trch+\frac12(\trch)^2
    &=-2\om\trch-|\hch|^2,\\
\nabs_3\hchb+\trchb\,\hchb
    &=\nabs\hot\xib-2\omb\hchb
      +(\eta+\etab-2\ze)\hot\xib-\aa,\\
\nabs_3\trchb+\frac12(\trchb)^2
    &=2\sdivs\xib-2\omb\trchb
      +2\xib\c(\eta+\etab-2\ze)-|\hchb|^2,\\
\nabs_4\hchb+\frac12\trch\,\hchb
    &=\nabs\hot\etab+2\om\hchb
      -\frac12\trchb\,\hch+\etab\hot\etab,\\
\nabs_4\trchb+\frac12\trch\trchb
    &=2\om\trchb-\hch\c\hchb+2\sdivs\etab
      +2|\etab|^2+2\rho,\\
\nabs_3\hch+\frac12\trchb\,\hch
    &=\nabs\hot\eta+2\omb\hch
      -\frac12\trch\,\hchb+\eta\hot\eta,\\
\nabs_3\trch+\frac12\trchb\trch
    &=2\omb\trch-\hch\c\hchb+2\sdivs\eta
      +2|\eta|^2+2\rho.
\end{align*}
The remaining Ricci coefficients satisfy
\begin{align*}
\nabs_4\eta
    &=-\chi\c(\eta-\etab)-\b,\\
\nabs_3\etab-\nabs_4\xib
    &=-4\om\xib-\chib\c(\etab-\eta)+\bb,\\
\nabs_3\ze
    &=-2\nabs\omb-\chib\c(\ze+\eta)
      +2\omb(\ze-\eta)+\chi\c\xib+2\om\xib-\bb,\\
\nabs_4\ze
    &=2\nabs\om+\chi\c(-\ze+\etab)
      +2\om(\ze+\etab)-\b,\\
\nabs_3\om+\nabs_4\omb
    &=\ze\c(\eta-\etab)-\eta\c\etab
      +4\om\omb+\rho.
\end{align*}
Moreover, the Codazzi equations are
\begin{equation}\label{codazzi}
\begin{aligned}
\sdivs\hch
    &=\frac12\nabs\trch-\ze\c\hch+\frac12\trch\,\ze-\b,\\
\sdivs\hchb
    &=\frac12\nabs\trchb+\ze\c\hchb-\frac12\trchb\,\ze+\bb.
\end{aligned}
\end{equation}
The torsion identities are
\begin{equation}\label{torsion}
\curls\eta=-\curls\etab=\si-\frac12\hch\wedge\hchb,\qquad\quad
\curls\xib=0.
\end{equation}
Finally, the Gauss equation is
\begin{equation}\label{gauss}
\K=-\frac14\trch\trchb+\frac12\hch\c\hchb-\rho.
\end{equation}
\end{prop}
\subsection{Renormalized main equations}
\begin{df}\label{renorr}
Following \cite{kl-ro,Kl-Ni}, we define
\begin{align}
\rhoc&:=\rho-\frac12\hch\c\hchb,&\sic&:=\si-\frac12\hch\wedge\hchb.\label{renorq}\\
\slu&:=\nabs\trch+\trch\,\ze,&\mu&:=-\sdivs\eta-\rhoc,\label{dfmu}
\end{align}
where $\mu$ is the mass aspect function.
\end{df}
The corrections in \eqref{renorq} cancel the $\a$ terms in the
outgoing Bianchi equations for $(\rhoc,\sic)$ and the $\aa$ terms
in their incoming equations.
\begin{prop}\label{muUtransport}
In the maximal-null foliation, we have
\begin{align}
\begin{split}\label{mu-transport}
\nabs_4\mu+\trch\,\mu
&=\frac12\trch\left(-\sdivs\etab+\rhoc-|\etab|^2\right)
+\hch\c(\nabs\hot\eta)+(\eta-\etab)\c\slu\\
&\quad-\frac14\trchb\,|\hch|^2
+2\eta\c\hch\c\etab-2\eta\c\b,
\end{split}\\
\begin{split}\label{U-transport}
\nabs_4\slu+\frac32\trch\,\slu
&=-\trch\,\b-2\om\slu-\hch\c\slu-\nabs|\hch|^2
+\trch\,\hch\c\etab-|\hch|^2(2\ze+\etab).
\end{split}
\end{align}
\end{prop}
\begin{proof}
These follow from the null structure, Bianchi and commutation
identities in Propositions 7.4.1 and 7.3.2 and Lemma 7.3.3
of \cite{Ch-Kl}.
\end{proof}
\begin{prop}\label{bianchirenor}
In the maximal-null foliation, the renormalized Bianchi
equations used below take the following form:
\begin{align*}
\nabs_3\a+\frac12\trchb\,\a
    &=-2\sld_2^*\b+4\omb\a+(\ze+4\eta)\hot\b-3\hch\left(\rhoc+\frac12\hch\c\hchb\right)
      -3{^*\hch}\left(\sic+\frac12\hch\wedge\hchb\right),\\
\nabs_4\b+2\trch\,\b
    &=\sld_2\a-2\om\b+(2\ze+\etab)\c\a,\\
\nabs_3\b+\trchb\,\b
    &=-\sld_1^*(\rhoc,-\sic)+2\omb\b+\xib\c\a+2\hch\c\bb+3\eta\left(\rhoc+\frac12\hch\c\hchb\right)
      +3{^*\eta}\left(\sic+\frac12\hch\wedge\hchb\right)\\
    &\quad+\frac12\nabs(\hch\c\hchb)
      +\frac12{^*\bigl(\nabs(\hch\wedge\hchb)\bigr)},\\
\nabs_4\rhoc+\frac32\trch\,\rhoc
    &=\sdivs\b+(\ze+2\etab)\c\b+\frac14\trchb\,|\hch|^2
      -\frac12\hch\c(\nabs\hot\etab+\etab\hot\etab),\\
\nabs_4\sic+\frac32\trch\,\sic
    &=-\curls\b-(\ze+2\etab)\c{^*\b}
      -\frac12\hch\wedge(\nabs\hot\etab+\etab\hot\etab),\\
\nabs_3\rhoc+\frac32\trchb\,\rhoc
    &=-\sdivs\bb+(\ze-2\eta)\c\bb+2\xib\c\b
      +\frac14\trch\,|\hchb|^2
      -\frac12\hchb\c(\nabs\hot\eta+\eta\hot\eta)\\
    &\quad-\frac12\hch\c
      \bigl(\nabs\hot\xib+(\eta+\etab-2\ze)\hot\xib\bigr),\\
\nabs_3\sic+\frac32\trchb\,\sic
    &=-\curls\bb+(\ze-2\eta)\c{^*\bb}-2\xib\c{^*\b}
      +\frac12\hchb\wedge(\nabs\hot\eta+\eta\hot\eta)\\
    &\quad-\frac12\hch\wedge
      \bigl(\nabs\hot\xib+(\eta+\etab-2\ze)\hot\xib\bigr),\\
\nabs_4\bb+\trch\,\bb
    &=\sld_1^*(\rhoc,\sic)+2\om\bb+2\hchb\c\b-3\etab\left(\rhoc+\frac12\hch\c\hchb\right)+3{^*\etab}\left(\sic+\frac12\hch\wedge\hchb\right)\\
    &\quad-\frac12\nabs(\hch\c\hchb)
      +\frac12{^*\bigl(\nabs(\hch\wedge\hchb)\bigr)},\\
      \nabs_3\bb+2\trchb\,\bb
    &=-\sld_2\aa-2\omb\bb+(2\ze-\eta)\c\aa-3\xib\left(\rhoc+\frac12\hch\c\hchb\right)+3{^*\xib}\left(\sic+\frac12\hch\wedge\hchb\right),\\
\nabs_4\aa+\frac12\trch\,\aa
    &=2\sld_2^*\bb+4\om\aa+(\ze-4\etab)\hot\bb-3\hchb\left(\rhoc+\frac12\hch\c\hchb\right)+3{^*\hchb}\left(\sic+\frac12\hch\wedge\hchb\right).
\end{align*}
\end{prop}
\begin{proof}
See equation~(26) of \cite{Luk-Rod17} and Lemma~4.9 of \cite{Shen23}, adapted to the present conventions.
\end{proof}

\subsection{Commutation identities}
We recall the following commutation identities, with the
curvature and Hodge conventions fixed above.
\begin{lem}\label{comm}
Let $U_{A_1\ldots A_k}$ be an $S$-tangent covariant
$k$--tensor in an Einstein vacuum spacetime. Then
\begin{align*}
[\nabs_4,\nabs_B]U_{A_1\ldots A_k}
    &=-\chi_{BC}\nabs_C U_{A_1\ldots A_k}
      +\xi_B\nabs_3 U_{A_1\ldots A_k}
      +(\ze_B+\etab_B)\nabs_4 U_{A_1\ldots A_k}\\
    &+\sum_{i=1}^k
      \Big(
          \chib_{A_iB}\xi_C-\chib_{BC}\xi_{A_i}
          +\chi_{A_iB}\etab_C-\chi_{BC}\etab_{A_i}
          +\in_{A_iC}{^*\b}_B
      \Big)
      U_{A_1\ldots C\ldots A_k},\\
[\nabs_3,\nabs_B]U_{A_1\ldots A_k}
    &=-\chib_{BC}\nabs_C U_{A_1\ldots A_k}
      +\xib_B\nabs_4 U_{A_1\ldots A_k}
      +(\eta_B-\ze_B)\nabs_3 U_{A_1\ldots A_k}\\
    &+\sum_{i=1}^k
      \Big(\chi_{A_iB}\xib_C-\chi_{BC}\xib_{A_i}+\chib_{A_iB}\eta_C-\chib_{BC}\eta_{A_i}-\in_{A_iC}{^*\bb}_B\Big)U_{A_1\ldots C\ldots A_k},\\
[\nabs_3,\nabs_4]U_{A_1\ldots A_k}
    &=-2\om\nabs_3 U_{A_1\ldots A_k}
      +2\omb\nabs_4 U_{A_1\ldots A_k}
      +2(\eta_B-\etab_B)\nabs_B U_{A_1\ldots A_k}\\
    &+2\sum_{i=1}^k
      \Big(
          \xi_{A_i}\xib_C-\xi_C\xib_{A_i}
          +\etab_{A_i}\eta_C-\etab_C\eta_{A_i}
          +\in_{A_iC}\si
      \Big)
      U_{A_1\ldots C\ldots A_k}.
\end{align*}
\end{lem}
\begin{proof}
See Lemma 7.3.3 in \cite{Ch-Kl}.
\end{proof}
\subsection{Magnetic curvature and spatial curl}
On a maximal slice, let $\in_{ijk}$ be the spatial volume form,
with $\in_{ABN}=\in_{AB}$ in the adapted orthonormal frame.
Spatial indices are raised and lowered with $g$.
\begin{prop}\label{EHiden}
The electric and magnetic parts of the spacetime curvature
satisfy:
\begin{align}
E_{ij}:=\R(T,e_i,T,e_j)&=\sRic_{ij}-k_i{}^\ell k_{\ell j},\label{initial:electric-curvature}\\
H_{ij}:={^*\R}(T,e_i,T,e_j)
&=\tfrac12{\in_i}^{lm}\R(e_l,e_m,T,e_j),\quad\;\; H_{NN}=\si,\quad\;\;
H_{AN}=\tfrac12{^*(\b-\bb)}_A,
\label{EHdecomposition}\\
\nab_i k_{jm}-\nab_j k_{im}&={\in_{ij}}^lH_{lm}.
\label{EHidentity}
\end{align}
Both $E$ and $H$ are symmetric and trace-free.
See \cite{Ch-Kl}.
\end{prop}
\begin{df}\label{dfdivcurl}
Let $(\Si,g)$ be an oriented three-dimensional Riemannian manifold
with Levi-Civita connection $\nab$ and volume form $\in_{ijk}$.
For a scalar function $f$ and a one-form $\xi$, set
\[
\De f:=\nab^i\nab_i f,\qquad\sdiv\xi:=\nab^i\xi_i,\qquad(\curl\xi)_i:={\in_i}^{lm}\nab_l\xi_m.
\]
For symmetric trace-free covariant two-tensors $U$ and $V$, set
\[
\begin{aligned}
 (\sdiv U)_i&:=\nab^jU_{ij},\qquad U\c V:=U_{ij}V^{ij},\\
 (\curl U)_{ij}
 &:=\frac{1}{2}\left(
 {\in_i}^{lm}\nab_lU_{mj}+{\in_j}^{lm}\nab_lU_{mi}
 \right).
\end{aligned}
\]
\end{df}
\section{Main theorem}\label{secmainthm}
\begin{df}\label{dfBEF}
For $t\geq0$ and $u_*\in\mathbb R\cup\{+\infty\}$, define
\begin{equation*}
 V_t(u_*):=V_t\cap\{u\leq u_*\},
\end{equation*}
where $V_t$ is given by \eqref{dfVt}. In particular, $V_t(+\infty)=V_t$. On a regular slab under consideration, we fix $t_*\geq0$ and $u_*\in\mathbb R\cup\{+\infty\}$ and write $V:=V_{t_*}(u_*)$. Unless otherwise specified, all integrals and suprema are restricted to $V$, with the suprema over $t,u$ taken only over the corresponding nonempty leaves. Throughout the paper, we abbreviate $S:=S(t,u)$.
\end{df}
\subsection{Main norms}\label{ssec8.1}
We define the curvature norm by
\begin{align}
\begin{split}\label{secRnorms}
\mr^2:={}&
\sup_u\int_{C_u\cap V}r\left(
\sum_{j=0}^1|(r\nabs)^j(\a,\b,\rhoc,\sic)|^2
+|\a_4|^2\right)d\mu_{C_u}\\
&+\sup_t\int_{\Si_t\cap V}r\left(
\sum_{j=0}^1|(r\nabs)^j(\a,\b,\rhoc,\sic,\bb)|^2
+|\a_4|^2\right)d\mu_g\\
&+\int_V\left(
\sum_{j=0}^1|(r\nabs)^j(\a,\b,\rhoc,\sic,\bb)|^2
+|\a_4|^2\right)d\mu_{\g}+\sup_{t,u}r^3\|(\a,\b,\rhoc,\sic)\|_{L^4(S)}^2,
\end{split}
\end{align}
where $\a_4:=r\left(\nabs_4\a+\tfrac52\trch\,\a+4\om\a\right)$. We define the maximal connection norm by
\begin{align}
\begin{split}\label{secknorms}
\mk^2:={}&
\sup_t\int_{\Si_t\cap V}r^{-1}
\left(
|G|^2+r^2|\nab G|^2+r^4|\nab\nabs G|^2
\right)d\mu_g\\
&+\int_Vr^{-2}
\left(
|G|^2+r^2|\nab G|^2+r^4|\nab\nabs G|^2
\right)d\mu_{\g}+\sup_{t,u}r\|G\|_{W^{1,4}(S)}^2,
\end{split}
\end{align}
where $G:=(\eps,\de,\nab\log\phi)$. For the quantities estimated by transport equations, we define
\begin{align}
\begin{split}\label{sssec8.1.3}
\mo^2:={}&
\int_{\mathbb R}du\,\sup_t\Bigg[
\|a^{1/2}r^{3/2}\nabs\mu\|_{L^2(S)}^2
+\|a^{1/2}(\eta,\hchb)\|_{W^{1,4}(S)}^2
+\|a^{1/2}r(\mu,\bb)\|_{L^4(S)}^2
\Bigg]\\
&+\sup_{t,u}\Bigg[
r\left\|\left(\trch-\frac2r,\hch\right)\right\|_{W^{1,4}(S)}^2
+r^2\|\slu\|_{H^1(S)}^2
\Bigg]\\
&+\sum_{j=0}^2\Bigg[
\sup_t\int_{\Si_t\cap V}r^{-1}
\left|(r\nabs)^j
\left(\eta,\hchb,\trch-\frac2r,\hch\right)\right|^2d\mu_g\\
&\qquad\qquad\qquad
+\int_Vr^{-2}
\left|(r\nabs)^j
\left(\eta,\hchb,\trch-\frac2r,\hch\right)\right|^2d\mu_{\g}
\Bigg]\\
&+\Bigg[
\int_{\mathbb R}du\,\sup_t\int_S ar
\left(
|(\eta,\hchb)|^4
+|r(\mu,\a,\b,\rhoc,\sic,\bb)|^4
\right)d\slg\\
&\qquad\qquad\qquad+\int_V\left(|(\eta,\hchb)|^4+|r(\mu,\a,\b,\rhoc,\sic,\bb)|^4\right)d\mu_{\g}\Bigg]^{1/2}\\
&+\Bigg[\sup_t\int_{\Si_t\cap V}r^3\left(|\eta|^6+|\hchb|^6\right)d\mu_g+\int_Vr^2\left(|\eta|^6+|\hchb|^6\right)d\mu_{\g}\Bigg]^{1/3}.
\end{split}
\end{align}
\begin{rk}
In the flux norms of $\mo$, we use both orders of taking
the supremum in $t$ and integrating in $u$.
Taking the supremum first yields a stronger mixed norm
than reversing the order. A similar choice of norms
appears in \cite[Remark~8]{Luk18}.
\end{rk}
On $\Si_0$, we define
\begin{align}
\begin{split}\label{initial:flux-norm}
\II_{(0)}^2:={}&
\sum_{j=0}^1\int_{\Si_0}r
\left(
|(r\nab)^j(\a,\b,\rhoc,\sic,\bb,\mu)|^2
+|(r\nabs)^j\slu|^2
\right)d\mu_g
+\int_{\Si_0}r|\a_4|^2d\mu_g\\
&+\sum_{j=0}^2\int_{\Si_0}r^{-1}\left|(r\nab)^j\left(\eta,\hchb,\hch,\trch-\frac{2}{r},k,\nab\log\phi\right)\right|^2\,d\mu_g\\
&+\sup_{w_1>0}\left[r^2\|\slu\|_{H^1(S_{(0)}(w_1))}^2+r\left\|\left(\hch,\trch-\frac2r\right)\right\|_{W^{1,4}(S_{(0)}(w_1))}^2\right]\\
&+\|r(\eta,\hchb)\|_{L^\infty(\Si_0)}^2+\|(\log a,\log\phi)\|_{L^\infty(\Si_0)}^2.
\end{split}
\end{align}
\subsection{Main theorem}\label{ssec8.3}
\begin{thm}\label{mainth}
Let $(\Si_0,g,k)$ be a maximal $(1,2)$--asymptotically flat vacuum initial data set of size at most $\ep_0$ in the sense of Definition~\ref{defsq}. If $\ep_0$ is sufficiently small, then its unique maximal globally hyperbolic development $(\M,\g)$ is future causally geodesically complete and admits a global maximal-null foliation $(C_u,\Si_t)$. Moreover, the norms $\mr$, $\mk$, and $\mo$ defined in \eqref{secRnorms}--\eqref{sssec8.1.3} satisfy 
\begin{equation}\label{finalesti}
\mr+\mk+\mo\les\ep_0.
\end{equation}
\end{thm}
The proof of Theorem~\ref{mainth} is given in Section~\ref{ssecproofmain} and uses the five intermediate theorems stated in Sections \ref{ssecini} and \ref{ssecmaininter}. Theorem~\ref{M0} bounds the initial norm $\II_{(0)}$, and Theorems~\ref{M1}--\ref{M3} establish the estimates for $\mr$, $\mk$, and $\mo$, respectively. Theorem~\ref{M4} establishes future causal geodesic completeness.

We set the bootstrap parameter $\ep:=\ep_0^{2/3}$ and choose $\ep_0$ sufficiently small so that $0<\ep\ll1$. Throughout the proof, $A\ll B$ means $CA<B$, where $C$ is a fixed universal constant chosen to dominate the implicit constants in the estimates below.
\subsection{Initial estimates}\label{ssecini}
\begin{thm}\label{M0}
Under the hypotheses of Theorem~\ref{mainth}, we have $\II_{(0)}\les\ep_0$.
\end{thm}
\begin{proof}
We use the initial cuts $S_{(0)}(w_1)$ constructed in Section~\ref{nullfoliation}. Write $\vr:=|x|$ for the coordinate radius and $r$ for the area radius. The construction gives $r\simeq \vr$ and $a=|\nab w|_g^{-1}=1+O(\ep_0)$, hence $\|\log a\|_{L^\infty(\Si_0)}\les\ep_0$. Using $N=a\nab w$, $\th=a\nab^2w|_{TS\times TS}$ and \eqref{6.6}, we estimate the initial connection coefficients directly from the weighted data bounds \eqref{ini} and the estimates for the initial cuts. In the coordinate-sphere region, the area formula controls $r/\vr-1$ and its first two weighted radial derivatives by sphere averages of the corresponding metric errors. This controls the replacement of $2/\vr$ by $2/r$ in the weighted trace estimate. We obtain
\[
 \sum_{j=0}^2\int_{\Si_0}r^{-1}
 \left|(r\nab)^j\left(k,\eta,\hch,\hchb,\trch-\frac{2}{r}\right)\right|^2
 \,d\mu_g\les\ep_0^2,
\]
where horizontal tensors are extended by zero in the $N$ direction for the spatial derivatives. The pointwise data bounds give
\begin{equation}\label{initial:connection-pointwise}
r\left(\|(\eta,\hchb,\hch,\trchb+2/r)\|_{W^{1,\infty}(S)}+\|\trch-2/r\|_{W^{2,\infty}(S)}+\|\eps\|_{W^{1,\infty}(S)}\right)\les\ep_0.
\end{equation}
The normalized maximal lapse satisfies $\De\phi=\phi|k|^2$, with $\phi\to1$ at infinity. The maximum principle gives $0<\phi\leq1$. Applying weighted elliptic estimates to $\De_g(\phi-1)=\phi|k|_g^2$, with $\phi-1\to0$ at infinity, and using the pointwise and weighted $L^2$ bounds for $k$ in \eqref{ini}, we obtain the following bounds for $\log\phi$:
\begin{equation*}
\|\log\phi\|_{L^\infty(\Si_0)}+\sup_{\Si_0}\langle \vr\rangle|\nab\log\phi|+\sum_{j=0}^2\left\|\vr^{j-1/2}\nab^j(\nab\log\phi)\right\|_{L^2(\Si_0)}\les\ep_0^2.
\end{equation*}
Since $r\simeq \vr$ and $|\nab r|\les1$, these are the lapse bounds required in \eqref{initial:flux-norm}.

By the Gauss--Codazzi identities \eqref{initial:electric-curvature} and \eqref{EHidentity} and the vacuum equations, the initial spacetime curvature is determined by $\sRic(g)$, $\nab k$ and quadratic terms in $k$. Thus one spatial derivative of the curvature uses at most three derivatives of $g$ and two of $k$, and \eqref{ini} gives
\[
\int_{\Si_0}r\left(|\R|^2+r^2|\nab\R|^2\right)d\mu_g\les\ep_0^2.
\]
The vacuum Bianchi identities express $\D_T\R$ in terms of spatial derivatives of the curvature components and terms of the form $k\c\R$, yielding the initial flux of $\a_4$. The bounds for $\hch$ and $\hchb$ above control the shear products in \eqref{renorq} and their first weighted spatial derivatives.

Finally, the definitions in \eqref{renorq} and \eqref{dfmu} give the weighted fluxes of $\mu$ and $\slu$. By \eqref{initial:connection-pointwise},
\[
\sup_{w_1>0}\left[r^2\|\slu\|_{H^1(S_{(0)}(w_1))}^2+r\|(\hch,\trch-2/r)\|_{W^{1,4}(S_{(0)}(w_1))}^2\right]+\|r(\eta,\hchb)\|_{L^\infty(\Si_0)}^2
\les\ep_0^2.
\]
Hence $\II_{(0)}\les\ep_0$, completing the proof of Theorem~\ref{M0}.
\end{proof}
\subsection{Main intermediate results}\label{ssecmaininter}
The following four theorems are the main intermediate results in the proof of Theorem \ref{mainth}.
\begin{thm}\label{M1}
Assume that $\II_{(0)}\les\ep_0$, $\mo\leq\ep$, $\mk\leq\ep$, and $\mr\leq\ep$. Then $\mr\les\ep_0$.
\end{thm}
In Section~\ref{secR}, we apply $r$--weighted estimates to the renormalized Bianchi and Teukolsky equations to control the curvature fluxes and bulks.
\begin{thm}\label{M2}
Assume that $\II_{(0)}\les\ep_0$, $\mo\leq\ep$, $\mk\leq\ep$,
and $\mr\les\ep_0$. Then $\mk\les\ep_0$.
\end{thm}
In Section~\ref{seck}, we use weighted elliptic estimates on the maximal slices to control the normal components of $k$ and the logarithmic lapse $\log\phi$.
\begin{thm}\label{M3}
Assume that $\II_{(0)}\les\ep_0$, $\mo\leq\ep$, $\mk\les\ep_0$, and $\mr\les\ep_0$. Then $\mo\les\ep_0$.
\end{thm}
In Section~\ref{secO}, we establish the connection coefficient flux and bulk estimates by integrating weighted transport equations along the outgoing null cones.
\begin{thm}\label{M4}
Assume that the maximal-null foliation is global and $\mo+\mk+\mr\les\ep_0$. Then the future development $(\M,\g)$ is future causally geodesically complete.
\end{thm}
In Section~\ref{sec:conclusions}, we control the energy growth along future-directed causal geodesics and show that every future-inextendible causal geodesic has infinite affine length.
\subsection{Proof of the main theorem}\label{ssecproofmain}
\paragraph{Step 1. Bootstrap}
By Theorem~\ref{M0}, the initial norm satisfies $\II_{(0)}\les\ep_0$. The initial data bounds in \eqref{ini} also give $\sRic(g),\nab k\in L^2(\Si_0)$ and a positive lower bound for $r_{\mathrm{vol}}(\Si_0,1)$. We use \cite[Theorem~2.2]{KRS} to obtain a local maximal evolution, with the lapse normalized by $\phi\to1$ at infinity, and extend the initial cuts to a local maximal-null foliation.

The norm $\II_{(0)}$ controls the initial contributions to
$\mr,\mk,\mo$. By local persistence of the weighted regularity,
we have $\mr+\mk+\mo\les\ep_0$ on a sufficiently short slab. This gives
\begin{equation}\label{B2}
    \mo\leq\ep,\qquad \mk\leq\ep,\qquad \mr\leq\ep.
\end{equation}
Let $t_*$ be the supremum of the times $t_1$ for which the maximal-null foliation exists on $V_{t_1}$ and \eqref{B2} holds there. Applying Theorems~\ref{M1}, \ref{M2}, and \ref{M3} successively gives, uniformly for $t_1<t_*$,
\begin{equation}\label{mrest}
\mr+\mk+\mo\les\ep_0.
\end{equation}
\paragraph{Step 2. Interior curvature estimates.}
Suppose that $t_*<+\infty$. We follow the interior curvature
argument of \cite[Sections~5.3--5.5]{Shen23}, taking the improved
$p=1$ bounds in \eqref{mrest} as the starting point. For the
$(\bb,\aa)$ Bianchi pair, we use the weight $p=-\delta_0$, with
$\ep\ll\delta_0\ll1$. This produces a positive bulk term of size
$\delta_0 r^{-1-\delta_0}|\aa|^2$, which absorbs the $\om r^{-\de_0}|\aa|^2$ contribution
by $r|\om|\les\ep$. With the corresponding lower-weight flux and
bulk estimates for the Bianchi hierarchy and its time-commuted
equations, the interior mean-value argument and the spatial
div--curl and Sobolev estimates yield, for any sufficiently small
$\de_B>0$,
\begin{align*}
\|(\a,\b)\|_{L^4(S(t,u))}
&\les\frac{\ep_0}{(1+t)^{3/2}},\\
\|(\rho,\si,\bb)\|_{L^4(S(t,u))}
&\les\frac{\ep_0}{(1+t)^{3/2-\de_B}\ujp^{\de_B}},&
\|\aa\|_{L^4(S(t,u))}
&\les\frac{\ep_0}{(1+t)^{1/2-\de_B}\ujp^{1+\de_B}}.
\end{align*}
Hence, we obtain
\[
\sup_{\substack{S(t,u)\\r(t,u)\leq2}}\|\R\|_{L^4(S(t,u))}\les\ep_0.
\]
By H\"older's inequality on the spheres and coarea on $\Si_t$, we conclude that
\[
\sup_{0\leq t<t_*}\int_{\Si_t\cap\{r\leq2\}}|\R|^2\,d\mu_g
\les\ep_0^2.
\]
\paragraph{Step 3. Exterior curvature estimates.}
On $\{r\geq1\}$, the slice flux in $\mr$ controls
$(\a,\b,\rhoc,\sic,\bb)$ in unweighted $L^2$. The bound
$r|\hch|\les\ep_0$ and the weighted slice estimate for $\hchb$
give the same control for $(\rho,\si)$. For $\aa$, the outgoing Bianchi equation gives
\[
 \nabs_4\aa+\tfrac12\trch\,\aa=4\om\aa+F,
 \qquad
 |F|\les|\nabs\bb|+|\Ga|\,|\bb|+|\hchb|\,|(\rho,\si)|.
\]
The curvature and quartic bulk bounds give $\|F\|_{L^2(V_t\cap\{r\geq1\})}\leq C(t_*)$, uniformly for $t<t_*$, with the renormalization terms controlled by $|\hch|\les\ep_0$ in the exterior.

Let $\vartheta$ be a smooth radial cutoff vanishing for $r\leq1$
and equal to one for $r\geq2$. Since $\bdiv e_4=\trch-2\om$,
\[
 \bdiv\bigl(\vartheta^2|\aa|^2e_4\bigr)
 =6\om\vartheta^2|\aa|^2
 +2\vartheta^2\langle\aa,F\rangle
 +e_4(\vartheta^2)|\aa|^2.
\]
We first integrate on $V_t\cap\{u\geq-U\}$, whose lateral flux
vanishes since $e_4u=0$. Set
$E_U(t)=\int_{\Si_t\cap\{u\geq-U\}}\vartheta^2|\aa|^2\,d\mu_g$.
The interior curvature bound controls the cutoff term on
$1\leq r\leq2$. Using $-\g(e_4,T)=1$, $0<\phi\leq1$, and the
exterior bound for $\om$, we obtain
\[
 E_U(t)\les E_U(0)+\int_0^t E_U(s)\,ds+C(t_*).
\]
The initial data bound $E_U(0)$ uniformly in $U$. Gr\"onwall's
inequality followed by $U\to\infty$ gives the exterior $L^2$
bound for $\aa$. Combining the interior and exterior estimates
yields
\begin{equation}\label{continuation:curvature}
 \sup_{t<t_*}\|\R\|_{L^2(\Si_t)}\leq C(t_*).
\end{equation}
\paragraph{Step 4. Continuation.}
By coarea on geodesic balls, the spatial isoperimetric inequality \eqref{normal:spatial-isoperimetry} gives $r_{\mathrm{vol}}(\Si_{t_1},1)\geq c>0$ uniformly for $t_1<t_*$. Together with \eqref{continuation:curvature}, the bounded $L^2$ curvature theorem \cite[Theorem~2.2 and Remark~2.4]{KRS} gives a common positive lifespan $\tau_*$ for the maximal evolution from $\Si_{t_1}$, with the lapse normalized by $\phi\to1$ at spatial infinity. Taking $t_1$ sufficiently close to $t_*$ extends the maximal foliation past $t_*$. On the extended maximal foliation, we extend the outgoing optical function $u$ as in \cite{Bieri} and construct the new cones from the extended axis using the normalization of Section~\ref{nullfoliation}.

Since $\ep_0\ll\ep$, the strict improvement in \eqref{mrest} and continuity of the bootstrap norms imply that \eqref{B2} continues to hold on a slightly longer slab, contradicting the definition of $t_*$. Thus $t_*=+\infty$, and \eqref{mrest} yields \eqref{finalesti}. Future causal geodesic completeness follows from Theorem~\ref{M4}. This completes the proof of Theorem~\ref{mainth}.
\section{Bootstrap assumptions and first consequences}
\label{firstboot}
Throughout this section, we work on a regular truncated maximal-null region $V=V_{t_*}(u_*)$, as in Definition~\ref{dfBEF}, and assume
\begin{equation}\label{B1}
\II_{(0)}\les\ep_0,\qquad
\mo\leq\ep,\qquad \mk\leq\ep,\qquad \mr\leq\ep,
\end{equation}
where $0<\ep_0\leq\ep\ll1$. The initial estimate follows from Theorem~\ref{M0}. All integrals and suprema below are restricted to $V$, and the estimates are uniform in the truncation parameters.
\subsection{Preliminary geometric identities and estimates}\label{2D}
By the Gauss equation \eqref{gauss}, the definition of $\rhoc$ in \eqref{renorq}, and the bootstrap sphere $L^4$ bounds, we have
\[
r\|\K-r^{-2}\|_{L^2(S)}\les\ep.
\]
\begin{prop}[Sphere Sobolev inequalities]\label{standardsobolev}
For a horizontal tensor $F$ and $2\leq p<\infty$,
\begin{align}
 r^{1/2}\|F\|_{L^\infty(S)}&\les\|F\|_{W^{1,4}(S)},\qquad
   r^{1-2/p}\|F\|_{L^p(S)}\les\|F\|_{H^1(S)},\label{bootstrap:sphere-sobolev}\\
 \|F\|_{L^2(S)}&\les r^{-1}\|F\|_{L^1(S)}+\|\nabs F\|_{L^1(S)}.\label{bootstrap:sphere-W11}
\end{align}
\end{prop}
\begin{proof}
See \cite[Lemma~4.1.3]{Kl-Ni}.
\end{proof}
\begin{prop}\label{Lpestimates}
For every $p\in(1,\infty)$, scalar function $f$, $\xi\in\sfr_1$,
and $U\in\sfr_2$, we have
\begin{align*}
 \|f-\overline f\|_{W^{2,p}(S)}&\les r^2\|\Des f\|_{L^p(S)},&
 \|\xi\|_{W^{1,p}(S)}&\les r\|\sld_1\xi\|_{L^p(S)},&
 \|U\|_{W^{1,p}(S)}&\les r\|\sld_2 U\|_{L^p(S)}.
\end{align*}
\end{prop}
\begin{proof}
See Corollary~2.3.1.1 in \cite{Ch-Kl}, with the scaled norms of Definition~\ref{Lpnorms}.
\end{proof}
\begin{cor}
Let $1<p<\infty$, let $F\in\sfr_0\cup\sfr_1\cup\sfr_2$, and let $w>0$ be a weight on $S$. Then
\begin{equation}\label{wLpestimates}
\int_S w|r\nabs F|^p\,d\slg\les\int_S w\left(|\dkb F|^p+|F|^p+|r\nabs\log w|^p|F|^p\right)d\slg,
\end{equation}
with a constant uniform in $S$ and independent of $w$.
\end{cor}
\begin{proof}
Apply Proposition~\ref{Lpestimates} to $w^{1/p}F$ and use the product rule. For scalar pairs, use the divergence-form Poisson estimate together with $\sld_1\sld_1^*=-\Des_0$.
\end{proof}
\subsection{Schematic notation}
\begin{df}\label{gammag}
We define the good and bad families by
\begin{align*}
\Gag&:=\left\{\trch-\frac2r,\,\trchb+\frac2r,\,\hch,\,\om,\,\omb,\,\eps,\,\de,\,\nab\log\phi,\,\ze,\,\etab,\,\frac{e_4r}{r}-\frac12\trch\right\},&
\Gab&:=\{\eta,\,\hchb\}.
\end{align*}
Set $\Gag^{(0)}:=\Gag$ and $\Gab^{(0)}:=\Gab$, and define
\begin{align*}
\Gag^{(1)}
 &:= (r\nabs)^{\leq1}\Gag
 \cup\left\{r\a,\,r\b,\,r\rhoc,\,r\sic,\,r\slu,\,r\K-\tfrac1r\right\},&\Gag^{(2)}
 &:= (r\nabs)^{\leq1}\Gag^{(1)}\cup\{r\a_4\},\\
\Gab^{(1)}
 &:= (r\nabs)^{\leq1}\Gab\cup\{r\mu,\,r\bb\},&
\Gab^{(2)}
 &:= (r\nabs)^{\leq1}\Gab^{(1)}.
\end{align*}
For the curvature components, set
\begin{align*}
\W&:=\{\a,\,\b,\,\rhoc,\,\sic\},
&\W^{(1)}&:=(r\nabs)^{\leq1}\W\cup\{\a_4\},\\
\Wb&:=\W\cup\{\bb\},&\Wb^{(1)}&:=(r\nabs)^{\leq1}\Wb\cup\{\a_4\}.
\end{align*}
Finally, write
\[
\Ga:=\Gag\cup\Gab,\qquad\Ga^{(i)}:=\Gag^{(i)}\cup\Gab^{(i)},\qquad 0\leq i\leq2.
\]
Using the pointwise bound for $r\Gag$, we abbreviate $r\Gag\c\Ga_l^{(i)}$ by $\Ga_l^{(i)}$ in schematic expressions, for $l\in\{g,b\}$ and $i=0,1,2$.
\end{df}
\begin{prop}\label{commutation}
On horizontal tensors, the angular commutators take the schematic form
\begin{align}
[r\nabs,\nabs_4]&=\Gag\c r\nabs+\Gag\c r\nabs_4+\Gag^{(1)},\label{bootstrap:commute4}\\
[r\nabs,\nabs_3]&=2(N(r)-1)\nabs
 +\Ga\c r\nabs+\Ga\c r\nabs_4+\Ga\c r\nabs_3+\Ga^{(1)}+r\Ga\c\Ga.
 \label{bootstrap:commute3}
\end{align}
\end{prop}
\begin{proof}
It follows from Lemma~\ref{comm} and \eqref{6.6}.
\end{proof}
\subsection{Bootstrap flux and nonlinear estimates}
\begin{lem}\label{decayGagGabGaa}
Assume \eqref{B1}. Then, we have the following bounds:
\begin{align}
 \sup_t\int_{\Si_t\cap V}r^{-1}|\Ga^{(2)}|^2\,d\mu_g
 +\int_Vr^{-2}|\Ga^{(2)}|^2\,d\mu_{\g}
 &\les\ep^2,\label{Gag-quadratic-bulk}\\
 \sup_t\int_{\Si_t\cap V}r
 |(\Gag^{(1)},\Gab,r\mu,r\bb)|^4\,d\mu_g
 +\int_V|(\Gag^{(1)},\Gab,r\mu,r\bb)|^4\,d\mu_{\g}
 &\les\ep^4,\label{Gamma-quartic-flux-bulk}\\
 \sup_t\int_{\Si_t\cap V}r^3|\Ga|^6\,d\mu_g
 +\int_Vr^2|\Ga|^6\,d\mu_{\g}
 &\les\ep^6,\label{maximal:all-sixth}\\
 \sup_t\int_{\Si_t\cap V}r|\Gab|^2|\Ga^{(1)}|^2\,d\mu_g
 +\int_V|\Gab|^2|\Ga^{(1)}|^2\,d\mu_{\g}
 &\les\ep^4,\label{wonderfulrp:bad-product-bulk}\\
 \sup_u\int_{C_u\cap V}r|\W^{(1)}|^2\,d\mu_{C_u}
 &\les\ep^2,\label{schematic-curvature-flux-bulk}\\
 \sup_{t,u}\left(r^{1/2}\|\Gag\|_{W^{1,4}(S)}+r\|\slu\|_{H^1(S)}\right)
 &\les\ep,\label{GagGabGawpointwise}\\
 \int du\,\sup_t\int_S ar
 \left(|\Gab|^4+|r(\Wb,\mu)|^4\right)d\slg
 &\les\ep^4,\label{null:zeroth-strong-quartic-bootstrap}\\
 \int du\,\sup_t\left(\|a^{1/2}\Gab\|_{W^{1,4}(S)}^2+\|a^{1/2}r(\mu,\bb)\|_{L^4(S)}^2\right)
 &\les\ep^2.\label{Gab-quadratic-flux-bulk}
\end{align}
\end{lem}
\begin{proof}
The estimates \eqref{Gab-quadratic-flux-bulk} and \eqref{null:zeroth-strong-quartic-bootstrap} follow from $\mo$, while \eqref{schematic-curvature-flux-bulk} follows from $\mr$. The sixth-power bound in $\mo$ gives
\begin{equation}\label{Gab-sixth-flux-bulk}
\sup_t\int_{\Si_t\cap V}r^3|\Gab|^6\,d\mu_g+\int_Vr^2|\Gab|^6\,d\mu_{\g}\les\ep^6.
\end{equation}
\paragraph{Sphere bounds for $\Gag$.}
We have from \eqref{e3e4r}, \eqref{gauss}, and \eqref{renorq}
\begin{equation}\label{bootstrap:gauss-coefficient}
\frac{e_4r}{r}-\frac12\trch=\frac{\overline{\phi\trch}-\phi\trch}{2\phi},\qquad r\K-\frac{1}{r}=-r\rhoc+\frac12\left(\trch-\frac2r\right)-\frac{r\trch}{4}\left(\trchb+\frac2r\right).
\end{equation}
Equations \eqref{6.6}, \eqref{dfmu}, and \eqref{bootstrap:gauss-coefficient}, together with the bounds in $\mk,\mo,\mr$ and Proposition~\ref{standardsobolev}, give \eqref{GagGabGawpointwise}. By sphere Sobolev, the sphere $L^4$ bounds in $\mr$ for $r\W$, and \eqref{bootstrap:gauss-coefficient}, we obtain
\begin{equation}\label{wonderfulrp:good-first-sphere}
r|\Gag|\les\ep,\qquad\sup_{t,u}r^{1/2}\|\Gag^{(1)}\|_{L^4(S)}\les\ep.
\end{equation}
\paragraph{$L^2$, $L^4$, and $L^6$ bounds.} The bootstrap assumption \eqref{B1}, \eqref{bootstrap:gauss-coefficient}, and
\begin{equation}\label{bootstrap:mass-and-trace}
 \mu=-\sdivs\eta-\rhoc,\qquad
 r\slu=r\nabs(\trch-2/r)+(r\trch)\eps
\end{equation}
give \eqref{Gag-quadratic-bulk}. The $L^6$ flux and bulk bounds \eqref{maximal:all-sixth} follow from \eqref{Gab-sixth-flux-bulk}, \eqref{Gag-quadratic-bulk}, and $r|\Gag|\les\ep$.

We now prove \eqref{Gamma-quartic-flux-bulk}. Indeed, the spacetime $L^4$ bounds for $(\Gab,r\Wb,r\mu)$ follow from $\mo\leq\ep$. The corresponding spatial $L^4$ fluxes follow from \eqref{null:zeroth-strong-quartic-bootstrap} and $d\mu_g=a\,du\,d\slg$. The zeroth-order $\Gag$ terms follow from $r|\Gag|\les\ep$ and \eqref{Gag-quadratic-bulk}. For any positive weight $w$ on $S$, integration by parts gives
\[
\int_S w|\nabs F|^4\,d\slg=-\int_S\left\langle F,\nabs^A\bigl(w|\nabs F|^2\nabs_AF\bigr)\right\rangle d\slg\les\|F\|_{L^\infty(S)}\int_S w\left(|\nabs F|^2|\nabs^2F|+|\nabs\log w|\,|\nabs F|^3\right)d\slg.
\]
Applying Young's inequality and multiplying by $r^4$, we obtain
\begin{equation}\label{bootstrap:angular-interpolation}
 \int_S w|r\nabs F|^4\,d\slg\les \|rF\|_{L^\infty(S)}^2\int_S wr^{-2}|(r\nabs)^2F|^2\,d\slg+\|rF\|_{L^\infty(S)}^4\int_S w|\nabs\log w|^4\,d\slg.
\end{equation}
By \eqref{6.6},
\begin{equation}\label{bootstrap:angular-densities}
 \nabs\log a=\eta-\eps,\qquad
 \nabs\log(\phi a)=\eta+\etab,\qquad
 \nabs\log\phi=\etab+\eps.
\end{equation}
Thus, for $w=a$ on slices or $w=\phi a$ in the spacetime bulk, we have
$|\nabs\log w|\les|\Gag|+|\eta|$. Taking $F\in\Gag$ in
\eqref{bootstrap:angular-interpolation} and using $r|\Gag|\les\ep$,
we get
\[
\int_S w|r\nabs\Gag|^4\,d\slg\les \ep^2\int_S wr^{-2}|(r\nabs)^2\Gag|^2\,d\slg+\ep^4\int_S w\left(|\Gag|^4+|\eta|^4\right)\,d\slg.
\]
Integrating in $r\,du$ on each slice and in $dt\,du$ over the spacetime region, we control the first term by \eqref{Gag-quadratic-bulk} and the second by the zeroth-order $L^4$ flux and bulk bounds above. This gives the $L^4$ flux and bulk bounds for $r\nabs\Gag$. The corresponding bounds for $r\slu$ and $r\K-1/r$ follow from \eqref{bootstrap:mass-and-trace} and \eqref{bootstrap:gauss-coefficient}, respectively. Together with the bounds already obtained, these prove \eqref{Gamma-quartic-flux-bulk}.
\paragraph{Products involving $\Gab$.}
Let $w=a$ on spatial slices and $w=\phi a$ in the spacetime bulk. Integration by parts gives
\[
\int_S w\left(|\Gab|^2|\nabs\Gab|^2+\frac12|\nabs(|\Gab|^2)|^2\right)d\slg
=-\int_S w|\Gab|^2\langle\Gab,\Des\Gab\rangle\,d\slg-\frac12\int_S |\Gab|^2
\langle\nabs w,\nabs(|\Gab|^2)\rangle\,d\slg.
\]
By Young's inequality, we absorb part of the last term into the L.H.S. and bound the remainder by $C\int_S w|\nabs\log w|^2|\Gab|^4\,d\slg$. Recalling from \eqref{bootstrap:angular-densities} that $|\nabs\log w|\les|\Gab|+\ep/r$, we infer
\begin{equation}\label{bootstrap:bad-product-interpolation}
 \int_S w|\Gab|^2|\nabs\Gab|^2\,d\slg\les\int_S w\left(|\Gab|^3|\nabs^2\Gab|+|\Gab|^6+\ep^2r^{-2}|\Gab|^4
 \right)d\slg.
\end{equation}
Integrate against $r^3\,du$ on each slice and against $r^2\,dt\,du$ over the spacetime region. By Cauchy--Schwarz, \eqref{Gag-quadratic-bulk} and \eqref{Gab-sixth-flux-bulk} bound the contribution of $|\Gab|^3|\nabs^2\Gab|$ by $C\ep^4$. The remaining terms are bounded by $C\ep^6$ directly from \eqref{Gab-sixth-flux-bulk} and \eqref{Gamma-quartic-flux-bulk}. Hence
\begin{equation}\label{maximal:bad-product-derivative-bound}
\sup_t\int_{\Si_t\cap V}r^3|\Gab|^2|\nabs\Gab|^2\,d\mu_g+\int_Vr^2|\Gab|^2|\nabs\Gab|^2\,d\mu_{\g}
\les\ep^4.
\end{equation}
The other products in \eqref{wonderfulrp:bad-product-bulk} follow from H\"older's inequality and \eqref{Gamma-quartic-flux-bulk}.
\end{proof}
\subsection{Weighted geometry and trace estimates}
\begin{lem}\label{bootstrap:weighted-estimates}
Assume \eqref{B1}. The measures and the area radius satisfy
\begin{align}
 d\mu_g
 &= a\,du\,d\slg,\qquad d\mu_{C_u}=\phi\,dt\,d\slg,\qquad\quad\;\, d\mu_{\g}=a\,du\,d\mu_{C_u},
 \label{bootstrap:coarea-exact}\\
 \phi/\ov{\phi}&=1+O(\ep),\qquad e_4r=1+O(\ep),\qquad N(r)>0,
 \label{bootstrap:radius-geometry}\\
 d\mu_{C_u}
 &\simeq dr\,d\slg,\qquad\quad\; d\mu_{\g}\simeq a\,du\,dr\,d\slg.
 \label{wonderfulrp:coarea}
\end{align}
At fixed $u$, the suprema in $r$ and $t$ are equivalent. The weighted sphere estimates are
\begin{align}
 \int du\,dr\,r^{-2}\|a^{1/2}\Gab\|_{H^2(S)}^2
 &\les\ep^2,
 \label{wonderfulrp:weighted-trace-energy}\\
 \int du\,dr\,r^{-1}\left(
 \|a^{1/2}\Ga^{(1)}\|_{L^4(S)}^2
 +\|a^{1/2}\Gab\|_{W^{1,4}(S)}^2\right)
 &\les\ep^2,
 \label{wonderfulrp:weighted-mixed}\\
 \int du\,\sup_t r\|a^{1/2}\Gab\|_{L^\infty(S)}^2
 +\int du\,dr\,\|a^{1/2}\Gab\|_{L^\infty(S)}^2
 &\les\ep^2.
 \label{wonderfulrp:coefficient-size}
\end{align}
The derivatives in these Sobolev norms act on the entire weighted field. In particular, for any tensor $F$,
\begin{equation}\label{bootstrap:coefficient-product}
\int_V|\Gab|^2|F|^2\,d\mu_{\g}\les\ep^2\min\left\{\sup_u\int_{C_u\cap V}r^{-1}|F|^2\,d\mu_{C_u},\; \sup_{t,u}\|F\|_{L^2(S)}^2\right\}.
\end{equation}
For $w=a$ or $w=\phi a$ and a scalar $f$, we have
\begin{align}
 \|w^{1/2}-\overline{w^{1/2}}\|_{L^\infty(S)}
 &\les r^{1/2}\|w^{1/2}\eta\|_{L^4(S)}
 +\ep\,\overline{w^{1/2}},
 \label{normal:sqrt-lapse-oscillation}\\
 \|w^{1/2}\|_{L^2(S)}|\overline f|
 &\les\|w^{1/2}f\|_{L^2(S)}
 +r^{1/2}\|w^{1/2}\eta\|_{L^4(S)}\|f\|_{L^2(S)},
 \label{bootstrap:weighted-mean}\\
 \|w^{1/2}(f-\overline f)\|_{L^2(S)}
 &\les r\|w^{1/2}\nabs f\|_{L^2(S)}
 +r^{1/2}\|w^{1/2}\eta\|_{L^4(S)}\|f-\overline f\|_{L^2(S)},
 \label{bootstrap:weighted-poincare}\\
 a|N(r)-1|
 &\les\ep a+\sqrt{ar}\,\|a^{1/2}\eta\|_{L^4(S)}
 +r\|a^{1/2}\eta\|_{L^4(S)}^2.
 \label{normal:radius-defect}
\end{align}
\end{lem}
\begin{proof}
The coarea identities give \eqref{bootstrap:coarea-exact}. Since $r|\nabs\log\phi|\les\ep$, sphere Poincar\'e and Sobolev imply $\phi/\overline\phi=1+O(\ep)$. By \eqref{e3e4r}, we have
\[
 e_4r=\frac{r\overline{\phi\trch}}{2\phi},
 \qquad
 aN(r)=\frac r2\overline{a\tr\th}.
\]
As $r\trch=2+O(\ep)$ and $r\tr\th=2+O(\ep)$, we obtain \eqref{bootstrap:radius-geometry}. Since $dr=\phi(e_4r)\,dt$ along $C_u$,
\[
 d\mu_{C_u}=(e_4r)^{-1}dr\,d\slg,
\]
which proves \eqref{wonderfulrp:coarea}.
\paragraph{Weighted sphere estimates.}
The identity $\nabs\log a=\eta-\eps\in\Ga$ gives
\[
 \big|(r\nabs)^{\leq2}(a^{1/2}\Gab)\big|
 \les a^{1/2}\left(
 |\Ga^{(2)}|+r|\Gab|\,|\Ga^{(1)}|+r^2|\Gab|\,|\Ga|^2
 \right).
\]
Consequently, by \eqref{wonderfulrp:coarea},
\begin{align*}
 \int du\,dr\,r^{-2}\|a^{1/2}\Gab\|_{H^2(S)}^2\les\int_V\left(
 r^{-2}|\Ga^{(2)}|^2+|\Gab|^2|\Ga^{(1)}|^2+r^2|\Ga|^6
 \right)d\mu_{\g}\les\ep^2,
\end{align*}
using \eqref{Gag-quadratic-bulk}, \eqref{wonderfulrp:bad-product-bulk}, and \eqref{maximal:all-sixth}. This proves \eqref{wonderfulrp:weighted-trace-energy}. The sphere embedding $H^1(S)\subset L^4(S)$ and the product rule give
\begin{equation}\label{bootstrap:weighted-H1}
\int du\,dr\,r^{-1}\|a^{1/2}\Ga^{(1)}\|_{L^4(S)}^2\les\int_V\left(
r^{-2}|\Ga^{(2)}|^2+(|\eta|^2+|\Gag|^2)|\Ga^{(1)}|^2\right)d\mu_{\g}\les\ep^2.
\end{equation}
The last product is controlled by \eqref{wonderfulrp:bad-product-bulk} for $\eta$, and by $r|\Gag|\les\ep$. Applying the same embedding to $r\nabs(a^{1/2}\Gab)$ and using \eqref{wonderfulrp:weighted-trace-energy} proves \eqref{wonderfulrp:weighted-mixed}. The Sobolev embeddings $H^2(S)\subset  W^{1,4}(S)\subset L^\infty(S)$ then give \eqref{wonderfulrp:coefficient-size}, using \eqref{Gab-quadratic-flux-bulk} for the first term.

To prove \eqref{bootstrap:coefficient-product}, we use \eqref{wonderfulrp:coarea} to write
\[
\int_V|\Gab|^2|F|^2\,d\mu_{\g}\simeq\int du\,dr\int_S a|\Gab|^2|F|^2\,d\slg.
\]
Bounding $a|\Gab|^2$ by its sphere supremum and using the first term of \eqref{wonderfulrp:coefficient-size} gives
\[
 \int_V|\Gab|^2|F|^2\,d\mu_{\g}
 \les\ep^2\sup_u\int_{C_u\cap V}r^{-1}|F|^2\,d\mu_{C_u}.
\]
Using the second term of \eqref{wonderfulrp:coefficient-size} gives
\[
 \int_V|\Gab|^2|F|^2\,d\mu_{\g}
 \les\ep^2\sup_{t,u}\|F\|_{L^2(S)}^2.
\]
Taking the minimum proves \eqref{bootstrap:coefficient-product}.
\paragraph{Lapse oscillation and the radius defect.}
For $w=a$ or $w=\phi a$, \eqref{bootstrap:angular-densities} gives
\[
 |\nabs w^{1/2}|\les w^{1/2}(|\eta|+|\Gag|).
\]
Sphere Poincar\'e and Sobolev yield
\begin{align*}
 \|w^{1/2}-\overline{w^{1/2}}\|_{L^\infty(S)}
 &\les r^{1/2}\|w^{1/2}\eta\|_{L^4(S)}
 +\ep\|w^{1/2}\|_{L^\infty(S)}\\
 &\les r^{1/2}\|w^{1/2}\eta\|_{L^4(S)}
 +\ep\,\overline{w^{1/2}}
 +\ep\|w^{1/2}-\overline{w^{1/2}}\|_{L^\infty(S)},
\end{align*}
which proves \eqref{normal:sqrt-lapse-oscillation}. Fix a sufficiently small universal constant $c>0$. On spheres where
\[
r^{1/2}\|w^{1/2}\eta\|_{L^4(S)}\leq c\,\ov{w^{1/2}},
\]
we have $w^{1/2}\simeq\overline{w^{1/2}}$. Cauchy--Schwarz and sphere Poincar\'e then give
\[
 \|w^{1/2}\|_{L^2(S)}|\overline f|
 \les\|w^{1/2}f\|_{L^2(S)},\qquad
 \|w^{1/2}(f-\overline f)\|_{L^2(S)}
 \les r\|w^{1/2}\nabs f\|_{L^2(S)}.
\]
On the remaining spheres, \eqref{normal:sqrt-lapse-oscillation} implies
\[
\|w^{1/2}\|_{L^\infty(S)}\les r^{1/2}\|w^{1/2}\eta\|_{L^4(S)}.
\]
Consequently,
\begin{align*}
 \|w^{1/2}\|_{L^2(S)}|\ov{f}|&\les r^{1/2}\|w^{1/2}\eta\|_{L^4(S)}\|f\|_{L^2(S)},\\
 \|w^{1/2}(f-\ov{f})\|_{L^2(S)}&\les r^{1/2}\|w^{1/2}\eta\|_{L^4(S)}\|f-\ov{f}\|_{L^2(S)}.
\end{align*}
These prove \eqref{bootstrap:weighted-mean} and \eqref{bootstrap:weighted-poincare}. For \eqref{normal:radius-defect}, taking $w=a$ in \eqref{normal:sqrt-lapse-oscillation} gives
\[
\|a^{1/2}-\ov{a^{1/2}}\|_{L^\infty(S)}\les r^{1/2}\|a^{1/2}\eta\|_{L^4(S)}+\ep a^{1/2}.
\]
Since
\[
\ov a=(\ov{a^{1/2}})^2+\ov{(a^{1/2}-\ov{a^{1/2}})^2},
\]
we obtain
\begin{align*}
|a-\ov a|\les a^{1/2}\left|a^{1/2}-\ov{a^{1/2}}\right|+\left\|a^{1/2}-\ov{a^{1/2}}\right\|_{L^\infty(S)}^2\les\ep a+\sqrt{ar}\,\|a^{1/2}\eta\|_{L^4(S)}+r\|a^{1/2}\eta\|_{L^4(S)}^2.
\end{align*}
Finally, \eqref{e3e4r} and $r\tr\th=2+O(\ep)$ imply $|aN(r)-\overline a|\les\ep\overline a$. Combining the two estimates proves \eqref{normal:radius-defect}.
\end{proof}
\begin{lem}\label{bootstrap:sobolev-traces}
Assume \eqref{B1}. Let $\Omega_{t,u}$ denote the exterior of a cut $S(t,u)\subset\Si_t\cap V$ on a complete maximal slice. For decaying tensors $Y$ with a metric connection and $p\geq0$, we have
\begin{align}
\|Y\|_{L^6(\Om_{t,u})}+\|Y\|_{L^4(S(t,u))}&\les\|\nab Y\|_{L^2(\Om_{t,u})},\label{normal:exterior-sobolev-trace}\\
\|r^{p/2}Y\|_{L^6(\Omega_{t,u})}^2&\les\int_{\Omega_{t,u}}r^p|\nab Y|^2\,d\mu_g,\label{normal:weighted-spatial-sobolev}\\
 r^p\|Y\|_{L^4(S(t,u))}^2
 &\les\int_{\Omega_{t,u}}r^p|\nab Y|^2\,d\mu_g.
 \label{normal:weighted-L4-trace}
\end{align}
When $\Si_t\subset V$, the spatial Sobolev estimates also hold on the whole slice, and we have the uniform isoperimetric inequality
\begin{equation}\label{normal:spatial-isoperimetry}
\Vol_g(E)^{2/3}\les\Area_g(\partial E),\qquad E\Subset\Si_t.
\end{equation}
For $R>0$ and a horizontal tensor $Y$, the cuts $S\subset C_u\cap V$ with $R\leq r\leq2R$ satisfy
\begin{equation}\label{bootstrap:cone-trace}
\begin{aligned}
\sup_S\|Y\|_{L^4(S)}^2&\les\int_{C_u\cap V\cap\left\{\frac{R}{2}\leq r\leq2R\right\}}\left(r^{-2}|Y|^2+|\nabs Y|^2+|\nabs_4Y|^2\right)d\mu_{C_u}+\mathbf{1}_{\left\{\frac{R}{2}\leq r_0\leq2R\right\}}\|Y\|_{L^4(S_0)}^2.
\end{aligned}
\end{equation}
Here $S_0=C_u\cap\Si_0$ and $r_0=r|_{S_0}$. The initial term is omitted for cones issuing from the axis.
\end{lem}
\begin{proof}
The spatial estimates are standard Sobolev, trace, and isoperimetric
inequalities and their weighted variants. See
\cite[Proposition~3.2.1 and Corollaries~3.2.1.1--3.2.1.2]{Ch-Kl}.
The cone trace estimate follows by localizing
\cite[Proposition~4.10]{Shen24} with a radial cutoff.
\end{proof}
\subsection{Geometric coefficient identities}
\begin{lem}\label{bootstrap:coefficient-families}
Assume \eqref{B1}. We have the schematic inclusions
\begin{align}
e_4r-\tfrac r2\trch,\ e_4r-1,\ r\tr\th-2 &\in r\Gag,\label{bootstrap:good-geometric-coefficients}\\
 \th-r^{-1}\slg,\ \kah,\ \xib,\ \nabs\log a,\ \nabs\log(\phi a)
 &\in\Ga,\label{bootstrap:general-geometric-coefficients}
\end{align}
with the corresponding bounds at each derivative order. Moreover,
\begin{equation}\label{bootstrap:incoming-radius}
 e_3r-\tfrac r2\trchb+2(N(r)-1)\in r\Gag.
\end{equation}
For any tensor $F$, the following estimates hold:
\begin{align}
 \int_{\Si_t\cap V}|N(r)-1|\,|F|^2\,d\mu_g
 &\les\ep\left(
 \int_{\Si_t\cap V}|F|^2\,d\mu_g
 +\sup_u r\|F\|_{L^2(S(t,u))}^2
 \right),
 \label{bootstrap:normal-radius-multiplier-flux}\\
 \int_V|N(r)-1|\,|F|^2\,d\mu_{\g}
 &\les\ep\left(
 \int_V|F|^2\,d\mu_{\g}
 +\sup_u\int_{C_u\cap V}r|F|^2\,d\mu_{C_u}
 \right).
 \label{bootstrap:normal-radius-multiplier}
\end{align}
In \eqref{bootstrap:normal-radius-multiplier-flux}, the supremum is over the cuts $S(t,u)\subset V$.
\end{lem}
\begin{proof}
The inclusions \eqref{bootstrap:good-geometric-coefficients} and \eqref{bootstrap:general-geometric-coefficients} follow from Definition~\ref{gammag}, \eqref{6.6}, and \eqref{bootstrap:angular-densities}. Since $e_3r=e_4r-2N(r)$, the expression in \eqref{bootstrap:incoming-radius} equals $(e_4r-1)-\tfrac r2(\trchb+2/r)$ and hence belongs to $r\Gag$. By \eqref{normal:radius-defect} and Young's inequality,
\[
a|N(r)-1|\les\ep a+\ep^{-1}r\|a^{1/2}\eta\|_{L^4(S)}^2.
\]
Using $d\mu_g=a\,du\,d\slg$ on each slice and $d\mu_{\g}=a\,du\,d\mu_{C_u}$ in spacetime, we obtain
\begin{align*}
\int_{\Si_t\cap V}|N(r)-1|\,|F|^2\,d\mu_g&\les\ep\int_{\Si_t\cap V}|F|^2\,d\mu_g+\ep^{-1}\left(\int du\,\sup_{t'}\|a^{1/2}\eta\|_{L^4(S(t',u))}^2\right)
\sup_u r\|F\|_{L^2(S(t,u))}^2,\\
\int_V|N(r)-1|\,|F|^2\,d\mu_{\g}&\les\ep\int_V|F|^2\,d\mu_{\g}+\ep^{-1}\left(\int du\,\sup_{t'}\|a^{1/2}\eta\|_{L^4(S(t',u))}^2\right)\sup_u\int_{C_u\cap V}r|F|^2\,d\mu_{C_u}.
\end{align*}
Applying \eqref{Gab-quadratic-flux-bulk}, we deduce both estimates.
\end{proof}
\subsection{Evolution lemma and higher integrability estimates}
\begin{lem}\label{evolution}
Let $X$ be horizontal and locally bounded near the axis. Suppose
\begin{equation}\label{evolution:equation}
 \nabs_4X+\la\trch\,X=F,\qquad \la\geq\tfrac12.
\end{equation}
For each $p\in[2,12]$, we have
\begin{align}
\begin{split}\label{evolution:Lp}
 &\int du\,\sup_{t'}\int_S ar^{p-3}|X|^p\,d\slg
 +\int_Vr^{p-4}|X|^p\,d\mu_{\g}\\
 &\qquad\les\int_{\Si_0\cap V}r^{p-3}|X|^p\,d\mu_g
 +\int_Vr^{p-3}|X|^{p-1}|F|\,d\mu_{\g}.
\end{split}
\end{align}
\end{lem}
\begin{proof}
Fix $u$ and write $S=S(t',u)$. By metric compatibility and \eqref{evolution:equation},
\[
e_4(|X|^p)=p|X|^{p-2}\langle X,F\rangle-p\la\trch\,|X|^p.
\]
By Lemma~\ref{dint}, $\frac{d}{dt'}\int_S f\,d\slg=\int_S\phi(e_4f+\trch\,f)\,d\slg$. Applying this identity to $f=ar^{p-3}|X|^p$ and using $e_4\log a=-2\om$, we obtain
\begin{align}
\begin{split}\label{evolution:identity}
&\frac{d}{dt'}\int_S ar^{p-3}|X|^p\,d\slg+\int_S\phi ar^{p-3}\left((p\la-1)\trch-(p-3)\frac{e_4r}{r}+2\om\right)|X|^p\,d\slg\\
&=p\int_S\phi ar^{p-3}|X|^{p-2}\langle X,F\rangle\,d\slg.
\end{split}
\end{align}
By \eqref{wonderfulrp:good-first-sphere} and \eqref{bootstrap:radius-geometry}, we have $r\trch=2+O(\ep)$, $e_4r=1+O(\ep)$, and $r\om=O(\ep)$. Thus
\[
r\left((p\la-1)\trch-(p-3)\frac{e_4r}{r}+2\om\right)=p(2\la-1)+1+O\bigl(\ep(p\la+1)\bigr)\geq\frac{1}{2}.
\]
We integrate \eqref{evolution:identity} from the initial cut for $u<0$ and from the vertex for $u\geq 0$. At a regular vertex, $a$ and $X$ are locally bounded, so the sphere term is $O(r^{p-1})$ and tends to zero. On each cone, we take the supremum over terminal cuts to bound the sphere term and integrate to the last cut to bound the bulk term. Integrating both bounds in $u$ and using \eqref{bootstrap:coarea-exact} yields \eqref{evolution:Lp}.
\end{proof}
We now derive the $L^{12}$ flux and bulk bounds of $\eta$. These bounds will be used to recover the quartic estimates for $r\nabs\Gab$ in Lemma \ref{null:augmented-quartic-bootstrap}.
\begin{prop}\label{null:eta12-bootstrap}
Assume \eqref{B1}. Then we have
\begin{equation}\label{null:eta12-weak}
\sup_{t'}\int_{\Si_{t'}\cap V}r^9|\eta|^{12}\,d\mu_g+\int_Vr^8|\eta|^{12}\,d\mu_{\g}\les\ep^{12}.
\end{equation}
\end{prop}
\begin{proof}
By Proposition~\ref{nulles},
\[
\nabs_4\eta+\tfrac12\trch\,\eta=-\b-\hch\c\eta+\tfrac12\trch\,\etab+\hch\c\etab.
\]
Using the schematic families of Definition~\ref{gammag}, we apply Lemma~\ref{evolution} with $X=\eta$, $\la=\tfrac12$, and $p=12$ to obtain
\begin{align*}
&\int du\,\sup_{t'}\int_S ar^9|\eta|^{12}\,d\slg+\int_Vr^8|\eta|^{12}\,d\mu_{\g}\\
&\qquad\les\int_{\Si_0\cap V}r^9|\eta|^{12}\,d\mu_g+\int_Vr^9|\eta|^{11}
\big(|\b|+|\Gag|\,|\eta|+r^{-1}|\Gag|\big)d\mu_{\g}.
\end{align*}
The initial integral is bounded by
\[
\|r\eta\|_{L^\infty(\Si_0)}^{10}\int_{\Si_0\cap V}r^{-1}|\eta|^2\,d\mu_g\les\ep_0^{12}.
\]
By \eqref{wonderfulrp:good-first-sphere}, $r|\Gag|\les\ep$. Hence, we obtain
\begin{align*}
r^9|\eta|^{11}\left(|\Gag|\,|\eta|+r^{-1}|\Gag|\right)\les\ep r^8|\eta|^{12}+\ep r^7|\eta|^{11}.
\end{align*}
After absorbing the first contribution into the bulk, it remains to estimate
\[
\int_Vr^9|\eta|^{11}|\b|\,d\mu_{\g}+\ep\int_Vr^7|\eta|^{11}\,d\mu_{\g}.
\]
For the curvature term, H\"older on $S$ followed by \eqref{bootstrap:sphere-sobolev} gives
\begin{equation}\label{null:eta12-sphere-holder}
\begin{aligned}
 \int_S a|\eta|^{11}|\b|\,d\slg
 &\leq\|a^{1/12}\eta\|_{L^{12}(S)}^{54/5}
 \|a^{1/2}\eta\|_{L^6(S)}^{1/5}\|\b\|_{L^{15}(S)}\\
 &\les r^{-1}\|a^{1/12}\eta\|_{L^{12}(S)}^{54/5}\|a^{1/2}\eta\|_{H^1(S)}^{1/5}\|\b\|_{H^1(S)}.
\end{aligned}
\end{equation}
Since $\|a^{1/12}\eta\|_{L^{12}(S)}^{54/5}
=(\int_S a|\eta|^{12}\,d\slg)^{1/2+2/5}$, we obtain
\begin{align}
\begin{split}\label{null:eta12-sphere-source}
 r^9\int_S a|\eta|^{11}|\b|\,d\slg
 &\les\left(\sup_{t'}\int_{S(t',u)}ar^9|\eta|^{12}\,d\slg\right)^{1/2}
 \left(\int_S ar^8|\eta|^{12}\,d\slg\right)^{2/5}\\
 &\quad\times
 \left(r^{-2}\|a^{1/2}\eta\|_{H^1(S)}^2\right)^{1/10}
 \left(r\|\b\|_{H^1(S)}^2\right)^{1/2}.
\end{split}
\end{align}
For fixed $u$, we take the sphere flux factor outside the $r$ integral and apply H\"older with exponents $(5/2,10,2)$ to the remaining factors. This gives
\begin{align*}
 \int dr\,r^9\int_S a|\eta|^{11}|\b|\,d\slg
 &\les \left(\sup_{t'}\int_{S(t',u)}ar^9|\eta|^{12}\,d\slg\right)^{1/2}
 \left(\int dr\int_S ar^8|\eta|^{12}\,d\slg\right)^{2/5}\\
 &\quad\times
 \left(\int dr\,r^{-2}\|a^{1/2}\eta\|_{H^1(S)}^2\right)^{1/10}
 \left(\int dr\,r\|\b\|_{H^1(S)}^2\right)^{1/2}.
\end{align*}
By \eqref{wonderfulrp:coarea} and \eqref{schematic-curvature-flux-bulk}, the last factor is bounded by $C\ep$. We now apply H\"older in $u$ with exponents $(2,5/2,10)$ to the first three factors. Using \eqref{wonderfulrp:coarea}, we find
\begin{align*}
 \int_Vr^9|\eta|^{11}|\b|\,d\mu_{\g}\les \ep
 \left(\int du\,\sup_{t'}\int_S ar^9|\eta|^{12}\,d\slg\right)^{1/2}
 \left(\int_Vr^8|\eta|^{12}\,d\mu_{\g}\right)^{2/5}
 \left(\int du\,dr\,r^{-2}\|a^{1/2}\eta\|_{H^1(S)}^2\right)^{1/10}.
\end{align*}
The last factor is bounded by $C\ep^{1/5}$ by \eqref{wonderfulrp:weighted-trace-energy}. Young's inequality with exponents $(2,5/2,10)$ then gives
\begin{equation}\label{null:eta12-beta-source}
 \int_Vr^9|\eta|^{11}|\b|\,d\mu_{\g}
 \leq \tau\left(
 \int du\,\sup_{t'}\int_S ar^9|\eta|^{12}\,d\slg
 +\int_Vr^8|\eta|^{12}\,d\mu_{\g}\right)
 +C_\tau\ep^{12}.
\end{equation}
For the second term, H\"older's inequality and \eqref{Gab-sixth-flux-bulk}, followed by Young's inequality, give
\begin{equation}\label{null:eta12-good-source}
 \ep\int_Vr^7|\eta|^{11}\,d\mu_{\g}
 \les\ep^2\left(\int_Vr^8|\eta|^{12}\,d\mu_{\g}\right)^{5/6}
 \leq\tau\int_Vr^8|\eta|^{12}\,d\mu_{\g}+C_\tau\ep^{12}.
\end{equation}
Choosing $\tau$ sufficiently small and absorbing the flux and bulk terms yields
\[
\int du\,\sup_{t'}\int_S ar^9|\eta|^{12}\,d\slg+\int_Vr^8|\eta|^{12}\,d\mu_{\g}\les\ep_0^{12}+\ep^{12}\les\ep^{12}.
\]
Taking the supremum over $t'$ in the slice coarea identity \eqref{bootstrap:coarea-exact} gives \eqref{null:eta12-weak}.
\end{proof}
\begin{lem}\label{null:augmented-quartic-bootstrap}
Assume \eqref{B1}. Then we have
\begin{equation}\label{null:full-quartic-bootstrap}
\sup_t\int_{\Si_t\cap V}r|\Ga^{(1)}|^4\,d\mu_g+\int_V|\Ga^{(1)}|^4\,d\mu_{\g}\les\ep^4.
\end{equation}
\end{lem}
\begin{proof}
By Lemma~\ref{decayGagGabGaa}, it remains to estimate $r\nabs\Gab$. The Codazzi and torsion equations \eqref{codazzi} and \eqref{torsion}, together with \eqref{dfmu}, give
\[
\sdivs\eta=-\mu-\rhoc,\qquad\curls\eta=\sic,\qquad\sdivs\hchb=\tfrac12\nabs\trchb+\ze\c\hchb-\tfrac12\trchb\,\ze+\bb.
\]
Since $r|\Gag|\les\ep$, we obtain
\[
r\big|(\sdivs\eta,\curls\eta,\sdivs\hchb)\big|\les|(\Gag^{(1)},\Gab,r\mu,r\bb)|.
\]
Their quartic flux and bulk are therefore controlled by Lemma~\ref{decayGagGabGaa}. We apply \eqref{wLpestimates} with $p=4$, taking $w=a$ for the slice flux and $w=\phi a$ for the spacetime bulk. By \eqref{bootstrap:angular-densities},
\[
|r\nabs\log w|^4|\Gab|^4\les r^4|\eta|^4|\Gab|^4+\ep^4|\Gab|^4.
\]
The second term is controlled by Lemma~\ref{decayGagGabGaa}. For the first term, H\"older gives
\begin{align*}
&\sup_t\int_{\Si_t\cap V}r^5|\eta|^4|\Gab|^4\,d\mu_g+\int_Vr^4|\eta|^4|\Gab|^4\,d\mu_{\g}\\
&\quad\leq\left(\sup_t\int_{\Si_t\cap V}r^9|\eta|^{12}\,d\mu_g\right)^{1/3}\left(\sup_t\int_{\Si_t\cap V}r^3|\Gab|^6\,d\mu_g\right)^{2/3}\\
&\qquad+\left(\int_Vr^8|\eta|^{12}\,d\mu_{\g}\right)^{1/3}\left(\int_Vr^2|\Gab|^6\,d\mu_{\g}\right)^{2/3}\les\ep^8,
\end{align*}
where the last bound follows from Proposition~\ref{null:eta12-bootstrap} and Lemma~\ref{decayGagGabGaa}. Integrating the weighted Hodge estimates using \eqref{bootstrap:coarea-exact}, we obtain the required quartic flux and bulk bounds for $r\nabs\Gab$.
\end{proof}
\section{\texorpdfstring{Hyperbolic estimates (Theorem \ref{M1})}{Curvature estimates}}\label{secR}
In this section, we prove Theorem~\ref{M1} by deriving and closing the curvature flux and bulk estimates on $V=V_t(u)$. To obtain the underlying energy identities, we use the $r^p$--weighted method of \cite{Da-Ro} at $p=1$, in the form adapted to the Bianchi equations in \cite{holzegel,Shen22,Shen23,Shen24}.
\subsection{General Bianchi pairs}
\begin{figure}[htbp]
    \centering
    \begin{tikzpicture}[
        x=0.95cm,y=0.95cm,
        font=\small,line cap=round,line join=round
    ]
    \foreach \offset/\interior in {0cm/0,7.5cm/1}{
        \begin{scope}[xshift=\offset]
            \coordinate (O) at (0,0);
            \coordinate (Izero) at (6.9,0);
            \coordinate (It) at (4.5,2.4);
            \coordinate (St) at (0,2.4);
            \ifnum\interior=0\relax
                \coordinate (B) at (0.7,0);
                \coordinate (S) at (3.1,2.4);
                \fill[yellow!35]
                    (B) -- (Izero) -- (It) -- (S) -- cycle;
                \draw[blue,thick] (B) -- (S);
                \fill[red] (B) circle[radius=1.7pt];
                \fill[red] (St) circle[radius=1.7pt];
                \node[below=4pt] at (B) {$S(0,u)$};
                \node[blue,anchor=west] at (2.1,1.1) {$C_u$};
                \node at (4.4,1.05) {$V_R$};
                \node at (3.45,3.35) {$u<0$};
            \else
                \coordinate (P) at (0,1);
                \coordinate (D) at (0.5,1.5);
                \coordinate (S) at (1.4,2.4);
                \fill[gray!20]
                    (O) -- (0.5,0) -- (D) -- (P) -- cycle;
                \fill[yellow!35]
                    (0.5,0) -- (Izero) -- (It)
                    -- (S) -- (D) -- cycle;
                \draw[blue,thick,dashed] (P) -- (D);
                \draw[blue,thick] (D) -- (S);
                \draw[violet,thick,dashed] (0.5,0) -- (D);
                \draw[violet,thin] (0.5,0.65) -- (0.8,0.65);
                \node[violet,anchor=west]
                    at (0.85,0.65) {$\mathcal T_\delta$};
                \node[rotate=90,font=\scriptsize]
                    at (0.23,0.55) {$r<\delta$};
                \draw[<->,thin] (0,-0.4) -- (0.5,-0.4)
                    node[midway,below=2pt] {$\delta$};
                \fill[red] (P) circle[radius=1.7pt];
                \fill[red] (D) circle[radius=1.7pt];
                \node[blue,anchor=east] at (0.85,2.05) {$C_u$};
                \node at (3.6,1.05) {$V_{\delta,R}$};
                \node at (3.45,3.35) {$u>0$};
            \fi
            \draw[black,thick] (O) -- (0,3.1);
            \draw[black,thick] (O) -- (Izero);
            \draw[black,thick] (St) -- (It);
            \draw[blue,thick] (O) -- (2.8,2.8);
            \draw[red,thick,dashed] (Izero) -- (It);
            \foreach \Q in {O,S,Izero,It}
                \fill[red] (\Q) circle[radius=1.7pt];
            \node[below left=2pt] at (O) {$O$};
            \node[anchor=west] at (0.1,3.05) {$\Vphi$};
            \node[above=5pt,red] at (S) {$S(t,u)$};
            \node[below=5pt] at (5.25,0) {$\Si_0$};
            \node[anchor=west] at (4.65,2.4) {$\Si_t$};
            \node[blue,anchor=east] at (1.45,1.55) {$C_0$};
            \node[red,anchor=west] at (5.85,1.2) {$\mathcal N_R$};
        \end{scope}
    }
    \end{tikzpicture}
    \caption{The truncated integration domains for $u<0$ (left) and $u>0$ (right). In the right panel, the gray strip represents the excised cylinder $\{r<\delta\}$ near the axis, whose boundary within $V_R$ is $\mathcal T_\delta$.}
    \label{Vt}
\end{figure}
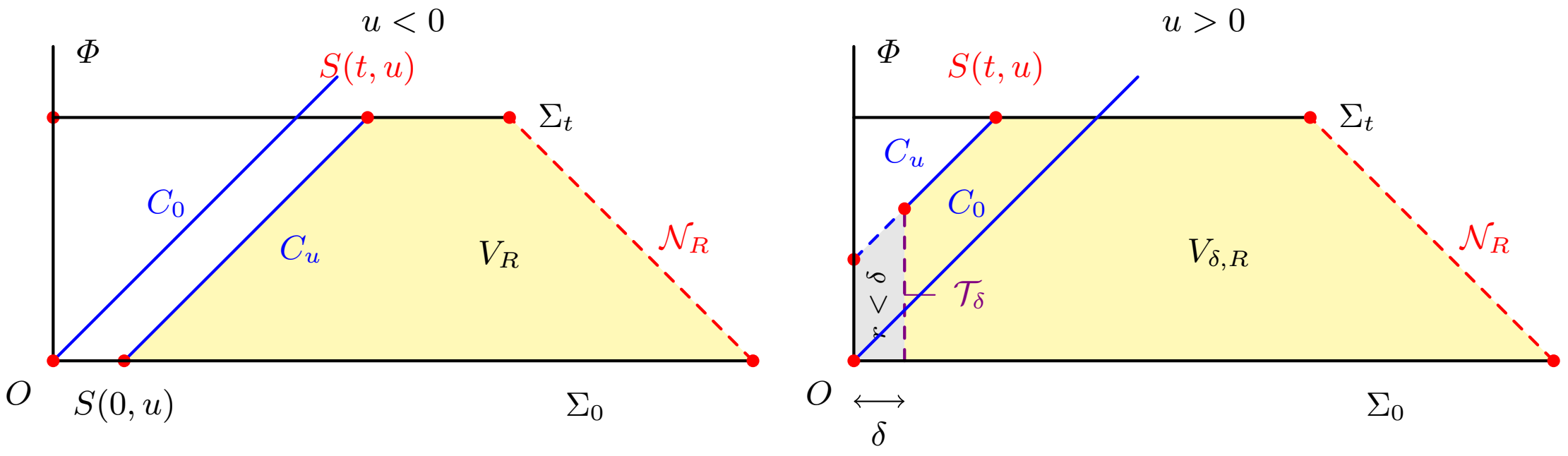

The following lemma gives the $r$-weighted divergence identities for Bianchi pairs.
\begin{lem}\label{keypoint}
Let $k\in\{1,2\}$ and $a_{(1)},a_{(2)}\in\mathbb R$.

\smallskip
\noindent\textup{(i)}
If $\psi_{(1)},h_{(1)}\in\sk_k$ and $\psi_{(2)},h_{(2)}\in\sk_{k-1}$ satisfy
\begin{equation}\label{bianchi1}
\begin{split}
    \nabs_3\psi_{(1)}+a_{(1)}\trchb\,\psi_{(1)}
    &=-k\sld_k^*\psi_{(2)}+h_{(1)},\\
    \nabs_4\psi_{(2)}+a_{(2)}\trch\,\psi_{(2)}
    &=\sld_k\psi_{(1)}+h_{(2)},
\end{split}
\end{equation}
then
\begin{equation}\label{div}
\begin{split}
    &\quad\bdiv(r|\psi_{(1)}|^2e_3)+k\bdiv(r|\psi_{(2)}|^2e_4)+\left(2a_{(1)}-\frac32\right)r\trchb|\psi_{(1)}|^2
    +k\left(2a_{(2)}-\frac32\right)r\trch|\psi_{(2)}|^2\\
    &=2kr\,\sdivs(\psi_{(1)}\cdot\psi_{(2)})
    +2r\psi_{(1)}\cdot h_{(1)}+2kr\psi_{(2)}\cdot h_{(2)}
    -2r\omb|\psi_{(1)}|^2-2kr\om|\psi_{(2)}|^2\\
    &\quad+\left(e_3r-\frac r2\trchb\right)|\psi_{(1)}|^2
    +k\left(e_4r-\frac r2\trch\right)|\psi_{(2)}|^2.
\end{split}
\end{equation}

\smallskip
\noindent\textup{(ii)}
If $\psi_{(1)},h_{(1)}\in\sk_{k-1}$ and $\psi_{(2)},h_{(2)}\in\sk_k$ satisfy
\begin{equation}\label{bianchi2}
\begin{split}
    \nabs_3\psi_{(1)}+a_{(1)}\trchb\,\psi_{(1)}
    &=\sld_k\psi_{(2)}+h_{(1)},\\
    \nabs_4\psi_{(2)}+a_{(2)}\trch\,\psi_{(2)}
    &=-k\sld_k^*\psi_{(1)}+h_{(2)},
\end{split}
\end{equation}
then
\begin{equation}\label{div2}
\begin{split}
    &\quad k\bdiv(r|\psi_{(1)}|^2e_3)+\bdiv(r|\psi_{(2)}|^2e_4)
    +k\left(2a_{(1)}-\frac32\right)r\trchb|\psi_{(1)}|^2
    +\left(2a_{(2)}-\frac32\right)r\trch|\psi_{(2)}|^2\\
    &=2kr\,\sdivs(\psi_{(1)}\cdot\psi_{(2)})
    +2kr\psi_{(1)}\cdot h_{(1)}+2r\psi_{(2)}\cdot h_{(2)}
    -2kr\omb|\psi_{(1)}|^2-2r\om|\psi_{(2)}|^2\\
    &\quad+k\left(e_3r-\frac r2\trchb\right)|\psi_{(1)}|^2
    +\left(e_4r-\frac r2\trch\right)|\psi_{(2)}|^2.
\end{split}
\end{equation}
\end{lem}
\begin{proof}
See Lemma~4.2 in \cite{Shen22}, with $p=1$.
\end{proof}
\begin{rk}
The Bianchi equations can be written in the form \eqref{bianchi1} or \eqref{bianchi2}. In particular,
\begin{itemize}
    \item $(\a_4,\as)$ satisfies \eqref{bianchi1} with $k=2$, $a_{(1)}=0$, $a_{(2)}=\frac52$;
    \item $(\a,\b)$ satisfies \eqref{bianchi1} with $k=2$, $a_{(1)}=\frac12$, $a_{(2)}=2$;
    \item $(\b,(\rhoc,-\sic))$ satisfies \eqref{bianchi1} with $k=1$, $a_{(1)}=1$, $a_{(2)}=\frac32$;
    \item $((\rhoc,\sic),-\bb)$ satisfies \eqref{bianchi2} with $k=1$, $a_{(1)}=\frac32$, $a_{(2)}=1$.
\end{itemize}
\end{rk}
\begin{prop}\label{keyintegral}
Under \eqref{B1}, let $(\psi_{(1)},\psi_{(2)})$ satisfy \eqref{bianchi1} or \eqref{bianchi2} on $V=V_t(u)$. Set $(c_1,c_2)=(1,k)$ in \eqref{bianchi1} and $(c_1,c_2)=(k,1)$ in \eqref{bianchi2}. Assume that the $(\psi_{(1)},\psi_{(2)})$ are locally bounded near the axis and that the displayed integrals are finite. Then
\begin{align}
\begin{split}\label{keyintegral:estimate}
&2c_1\int_{C_u\cap V}r|\psi_{(1)}|^2\,d\mu_{C_u}
+\int_{\Si_t\cap V}r\bigl(c_1|\psi_{(1)}|^2+c_2|\psi_{(2)}|^2\bigr)\,d\mu_g\\
&\quad+\int_V\Bigl(c_1(3-4a_{(1)})|\psi_{(1)}|^2+c_2(4a_{(2)}-3)|\psi_{(2)}|^2\Bigr)d\mu_{\g}\\
&\leq\int_{\Si_0\cap V}r\bigl(c_1|\psi_{(1)}|^2+c_2|\psi_{(2)}|^2\bigr)\,d\mu_g+2\int_Vr\bigl(c_1|\psi_{(1)}|\,|h_{(1)}|+c_2|\psi_{(2)}|\,|h_{(2)}|\bigr)\,d\mu_{\g}+\mathcal E(V),
\end{split}
\end{align}
where
\begin{align*}
\mathcal E(V)\les\ep\left[\sup_{u'\leq u}\int_{C_{u'}\cap V}r|\psi_{(1)}|^2\,d\mu_{C_{u'}}+\int_V\bigl(|\psi_{(1)}|^2+|\psi_{(2)}|^2\bigr)\,d\mu_{\g}\right]+\int_Vr|\Gab|\,|\psi_{(1)}|\,|\psi_{(2)}|\,d\mu_{\g}.
\end{align*}
\end{prop}
\begin{proof}
Let $V_R\nearrow V$ be the exhaustion with outer incoming null boundary $\mathcal N_R$, and set
\[
V_{\de,R}:=V_R\cap\{r\geq\de\},\qquad \mathcal T_\de:=V_R\cap\{r=\de\}.
\]
Figure~\ref{Vt} illustrates the truncated domains. In \eqref{div} and \eqref{div2}, use
\[
\bdiv\bigl(r(\psi_{(1)}\cdot\psi_{(2)})\bigr)=r\sdivs(\psi_{(1)}\cdot\psi_{(2)})+r(\eta+\etab)\c(\psi_{(1)}\cdot\psi_{(2)}).
\]
Moving the spacetime divergence to the left gives the current
\[
J=r\Big(c_1|\psi_{(1)}|^2e_3+c_2|\psi_{(2)}|^2e_4-2k(\psi_{(1)}\cdot\psi_{(2)})\Big),
\]
where the horizontal contraction is identified with its metric dual. Since $c_1c_2=k$ and $|\psi_{(1)}\cdot\psi_{(2)}|\leq k^{-1/2}|\psi_{(1)}|\,|\psi_{(2)}|$, the current is future causal and its flux through $\mathcal N_R$ is nonnegative. The horizontal part of $J$ is tangent to the remaining boundary hypersurfaces. Stokes therefore gives
\begin{align}
\begin{split}\label{linear:stokes}
\int_{V_{\delta,R}}\bdiv J\,d\mu_{\g}
&\geq2c_1\int_{C_u\cap V_{\delta,R}}r|\psi_{(1)}|^2\,d\mu_{C_u}
+\int_{\Si_t\cap V_{\delta,R}}
 r\bigl(c_1|\psi_{(1)}|^2+c_2|\psi_{(2)}|^2\bigr)\,d\mu_g\\
&\quad-\int_{\Si_0\cap V_{\delta,R}}
 r\bigl(c_1|\psi_{(1)}|^2+c_2|\psi_{(2)}|^2\bigr)\,d\mu_g\\
&\quad-\delta\int_0^t\int_{\Si_{t'}\cap\mathcal T_\delta}
 \frac{\phi}{N(r)}
 \bigl(c_1e_3(r)|\psi_{(1)}|^2+c_2e_4(r)|\psi_{(2)}|^2\bigr)
 \,d\slg\,dt'.
\end{split}
\end{align}
For $u<0$, the cylinder is absent for sufficiently small $\delta$. Otherwise, regularity of the foliation near the axis, local boundedness of the fields, and the area bound $4\pi\delta^2$ give
\begin{equation}\label{linear:axis-limit}
 \delta\int_0^t\int_{\Si_{t'}\cap\mathcal T_\delta}
 \frac{\phi}{N(r)}
 \bigl(c_1|e_3r|\,|\psi_{(1)}|^2+c_2|e_4r|\,|\psi_{(2)}|^2\bigr)
 \,d\slg\,dt'=O_{t,R}(\delta^3)\longrightarrow0.
\end{equation}
We integrate \eqref{div} or \eqref{div2} using \eqref{linear:stokes}, and let $\delta\to0$ and then $R\to\infty$. The expansions $r\trchb=-2+O(\ep)$ and $r\trch=2+O(\ep)$ give the signed bulk coefficients in \eqref{keyintegral:estimate}. The errors in these expansions and the terms with $\om,\omb$ are bounded by $C\ep\int_V(|\psi_{(1)}|^2+|\psi_{(2)}|^2)\,d\mu_{\g}$. By \eqref{bootstrap:good-geometric-coefficients} and \eqref{bootstrap:incoming-radius}, the terms with derivatives of $r$ are bounded pointwise by $C\ep(|\psi_{(1)}|^2+|\psi_{(2)}|^2)+2|N(r)-1|\,|\psi_{(1)}|^2$. Applying \eqref{bootstrap:normal-radius-multiplier}, we obtain
\begin{align}
\begin{split}\label{linear:radius-errors}
 &\int_V\left[\left|e_3r-\tfrac r2\trchb\right||\psi_{(1)}|^2
 +\left|e_4r-\tfrac r2\trch\right||\psi_{(2)}|^2\right]d\mu_{\g}\\
 &\les\ep\left[\sup_{u'\leq u}\int_{C_{u'}\cap V}r|\psi_{(1)}|^2\,d\mu_{C_{u'}}
 +\int_V\bigl(|\psi_{(1)}|^2+|\psi_{(2)}|^2\bigr)d\mu_{\g}\right].
\end{split}
\end{align}
The remaining contribution from the horizontal divergence is bounded by
\[
 C\int_V r|\eta+\etab|\,|\psi_{(1)}|\,|\psi_{(2)}|\,d\mu_{\g}.
\]
Using $\eta\in\Gab$ and $r|\etab|\les\ep$ gives the stated bound for $\mathcal E(V)$. Taking absolute values of the source pairings with their coefficients $2c_1$ and $2c_2$ gives \eqref{keyintegral:estimate}.
\end{proof}
\subsection{Nonlinear error estimates}
\begin{df}\label{defschematicerror}
We introduce the schematic error family
\begin{equation*}
\err[\Ga]:=r^{-1}\Ga\c\Ga^{(1)}+\Ga\c\Ga\c\Ga,
\end{equation*}
and its augmented first-order family
\begin{equation*}
\err[\Ga]^{(1)}:=r^{-1}\Ga\c\Ga^{(2)}+r^{-1}\Ga^{(1)}\c\Ga^{(1)}+\Ga\c\Ga\c\Ga^{(1)}.
\end{equation*}
We set $\err[\Ga]^{(0)}=\err[\Ga]$ and use the families of Definition \ref{gammag}.
\end{df}
We estimate the nonlinear source terms in the Bianchi and Teukolsky equations.
\begin{thm}\label{Junkman}
For $q\in\{0,1\}$, define
\begin{equation}\label{deferr}
 \NN^q[\psi_{(1)},\psi_{(2)}](V)
 :=\int_V r\left|\psi_{(1)}^{(q)}\c\psi_{(2)}^{(q)}\right|\,d\mu_{\g}.
\end{equation}
Here $\c$ denotes a contraction appearing in the curvature equations. The superscripts $(q)$ are understood as in Definition~\ref{gammag}, with the schematic product rule applied to products. Under the bootstrap assumption \eqref{B1}, we have
\begin{align}
 \begin{aligned}[b]
 &\NN^1[\W,\Gab\c\Wb](V)
  +\NN^1[\W,r^{-1}\Ga\c\Ga^{(1)}](V)\\
 &+\NN^1[\Wb,\Gab\c\W](V)
  +\NN^1[\Wb,\Gab\c\slu](V)+\NN^1[\Wb,r^{-1}\Gag\c\Ga^{(1)}](V)
 \end{aligned}
 &\les\ep^3,\label{wonderfulrp:estimate}\\
 \NN^1[\W,\Ga\c\Ga\c\Ga](V)
 +\NN^1[\Wb,\Gag\c\Gag\c\Ga](V)
 &\les\ep^4,\label{wonderfulrp:cubicestimate}\\
 \NN^0[\W^{(1)},r\Gab\c\Gab\c\Gab\c\Gab](V)
 &\les\ep^5.\label{wonderfulrp:Gab4}
\end{align}
\end{thm}
\begin{proof}
The proof is divided into three steps.
\paragraph{Step 1. Preliminary product estimates.}
We use Lemmas~\ref{decayGagGabGaa} and \ref{bootstrap:weighted-estimates}. Since $r\Wb^{(1)}\subseteq\Ga^{(2)}$, estimate \eqref{Gag-quadratic-bulk} gives $\|\W^{(1)}\|_{L^2(V)}\leq\|\Wb^{(1)}\|_{L^2(V)}\les\ep$. Applying \eqref{bootstrap:coefficient-product} to $r\W^{(1)}$ and using \eqref{schematic-curvature-flux-bulk} and Cauchy--Schwarz gives
\begin{equation}\label{wonderfulrp:separated-products}
 \int_Vr|\W^{(1)}|\,|\Gab|\,|F|\,d\mu_{\g}
 \les\left(\int_Vr^2|\W^{(1)}|^2|\Gab|^2\,d\mu_{\g}\right)^{1/2}
 \|F\|_{L^2(V)}
 \les\ep^2\|F\|_{L^2(V)}.
\end{equation}
Taking $F=r|\Gab|^3$ and using \eqref{maximal:all-sixth}, we obtain \eqref{wonderfulrp:Gab4} from
\begin{equation}\label{wonderfulrp:quartic-bound}
 \NN^0[\W^{(1)},r\Gab\c\Gab\c\Gab\c\Gab](V)
 \les\int_Vr^2|\W^{(1)}|\,|\Gab|^4\,d\mu_{\g}
 \les\ep^2\|r|\Gab|^3\|_{L^2(V)}
 \les\ep^5.
\end{equation}
We recall from \eqref{null:full-quartic-bootstrap} the $L^4$ bulk bound $\|\Ga^{(1)}\|_{L^4(V)}\les\ep$. H\"older gives
\begin{equation}\label{wonderfulrp:good-mixed-product}
 \int_V|\Wb^{(1)}|\,|\Ga^{(1)}|^2\,d\mu_{\g}
 \leq\|\Wb^{(1)}\|_{L^2(V)}\|\Ga^{(1)}\|_{L^4(V)}^2
 \les\ep^3,
\end{equation}
\paragraph{Step 2. The quadratic sources.}
The product rule gives
\begin{align*}
\NN^1[\W,\Gab\c\Wb](V)\les\int_Vr|\W^{(1)}|\,|\Gab|\,|\Wb^{(1)}|\,d\mu_{\g}+\int_V|\W^{(1)}|\,|\Gab^{(1)}|\,|r\Wb|\,d\mu_{\g}\les\ep^3.
\end{align*}
The first integral uses \eqref{wonderfulrp:separated-products} with $F=|\Wb^{(1)}|$, and the second uses \eqref{wonderfulrp:good-mixed-product} and $r\Wb\subseteq\Ga^{(1)}$.
The product rule gives
\begin{align*}
 \NN^1[\W,r^{-1}\Ga\c\Ga^{(1)}](V)\les\int_V|\W^{(1)}|\,|\Ga^{(1)}|^2\,d\mu_{\g}+\int_V|\W^{(1)}|\,|\Ga|\,|\Ga^{(2)}|\,d\mu_{\g}.
\end{align*}
The first integral is bounded by \eqref{wonderfulrp:good-mixed-product}. For the second, split $\Ga=\Gag\cup\Gab$ and use $r|\Gag|\les\ep$, \eqref{Gag-quadratic-bulk}, and \eqref{wonderfulrp:separated-products} to obtain
\begin{align*}
 \int_V|\W^{(1)}|\,|\Ga|\,|\Ga^{(2)}|\,d\mu_{\g}
 &\les\ep\|\W^{(1)}\|_{L^2(V)}\|r^{-1}\Ga^{(2)}\|_{L^2(V)}+\int_Vr|\W^{(1)}|\,|\Gab|\,|r^{-1}\Ga^{(2)}|\,d\mu_{\g}
 \les\ep^3.
\end{align*}
We also have from $r\W\subseteq r\Wb\subseteq \Ga^{(1)}$, \eqref{wonderfulrp:separated-products} and \eqref{wonderfulrp:good-mixed-product}
\begin{align*}
 \NN^1[\Wb,\Gab\c\W](V)
 &\les\int_Vr|\Wb^{(1)}|\,|\Gab|\,|\W^{(1)}|\,d\mu_{\g}+\int_V|\Wb^{(1)}|\,|\Gab^{(1)}|\,|r\W|\,d\mu_{\g}
 \les\ep^3.
\end{align*}
The product rule gives
\begin{align*}
 \NN^1[\Wb,\Gab\c\slu](V)\les\int_V|\Wb^{(1)}|\,|\Gab^{(1)}|\,|r\slu|\,d\mu_{\g}+\int_Vr|\Wb^{(1)}|\,|\Gab|\,|(r\nabs)^{\leq1}\slu|\,d\mu_{\g}.
\end{align*}
The first integral is bounded by \eqref{wonderfulrp:good-mixed-product}, since $r\slu\in\Gag^{(1)}$. By Cauchy--Schwarz, \eqref{bootstrap:coefficient-product}, and \eqref{GagGabGawpointwise}, the second integral satisfies
\begin{align*}
 \int_Vr|\Wb^{(1)}|\,|\Gab|\,|(r\nabs)^{\leq1}\slu|\,d\mu_{\g}
 &\leq\|\Wb^{(1)}\|_{L^2(V)}
 \big\||\Gab|\,|(r\nabs)^{\leq1}(r\slu)|\big\|_{L^2(V)}\\
 &\les\ep\|\Wb^{(1)}\|_{L^2(V)}\sup_{t,u}r\|\slu\|_{H^1(S(t,u))}
 \les\ep^3.
\end{align*}
Finally, the product rule gives
\begin{align*}
 \NN^1[\Wb,r^{-1}\Gag\c\Ga^{(1)}](V)
 &\les\int_V|\Wb^{(1)}|\,|\Gag^{(1)}|\,|\Ga^{(1)}|\,d\mu_{\g}
 +\ep\|\Wb^{(1)}\|_{L^2(V)}\|r^{-1}\Ga^{(2)}\|_{L^2(V)}
 \les\ep^3.
\end{align*}
Here the first term uses \eqref{wonderfulrp:good-mixed-product}, and the second uses $r|\Gag|\les\ep$ and \eqref{Gag-quadratic-bulk}. This proves \eqref{wonderfulrp:estimate}.
\paragraph{Step 3. The cubic sources.}
The product rule and $r|\Gag|\les\ep$ give
\begin{equation}\label{wonderfulrp:cubic-product}
 r|(\Ga\c\Ga\c\Ga)^{(1)}|
 \les\ep^2r^{-1}|\Ga^{(1)}|
 +\ep|\Gab|\,|\Ga^{(1)}|
 +r|\Gab|^2|\Ga^{(1)}|.
\end{equation}
The first two terms have $L^2(V)$ norm bounded by $C\ep^3$ by \eqref{Gag-quadratic-bulk} and \eqref{null:full-quartic-bootstrap}. Pairing them with $\W^{(1)}$ in $L^2(V)$ gives $C\ep^4$. For the last term, use \eqref{wonderfulrp:separated-products} with $F=|\Gab|\,|\Ga^{(1)}|$ and apply \eqref{null:full-quartic-bootstrap} again to obtain
\[
\int_Vr|\W^{(1)}|\,|\Gab|^2|\Ga^{(1)}|\,d\mu_{\g}\les\ep^2\|\Gab\|_{L^4(V)}\|\Ga^{(1)}\|_{L^4(V)}\les\ep^4.
\]
For the second cubic source, the product rule and $r|\Gag|\les\ep$ give
\[
r|(\Gag\c\Gag\c\Ga)^{(1)}|\les\ep^2r^{-1}|\Ga^{(1)}|+\ep|\Ga^{(1)}|^2.
\]
Both terms have $L^2(V)$ norm bounded by $C\ep^3$ by \eqref{Gag-quadratic-bulk} and \eqref{null:full-quartic-bootstrap}. Pairing with $\Wb^{(1)}$ in $L^2(V)$ proves \eqref{wonderfulrp:cubicestimate}.
\end{proof}
\subsection{Proof of Theorem~\ref{M1}}\label{exteriorest}
\paragraph{Step 1. The commuted equations.}
Set
\[
\a_4:=r\left(\nabs_4\a+\tfrac52\trch\,\a+4\om\a\right),\qquad\quad \as:=r\sld_2\a.
\]
The Bianchi equations and the commutator identities give the Teukolsky pair below. See also Lemma~5.5 in \cite{Shen23}.
\begin{align*}
 \nabs_3\a_4
 &=-2\sld_2^*\as+\frac4r\a-\frac{2(N(r)-1)}r\a_4+\err[\Ga]^{(1)},\\
 \nabs_4\as+\tfrac52\trch\,\as
 &=\sld_2\a_4+\err[\Ga]^{(1)}.
\end{align*}
For the remaining pairs, we use the renormalized Bianchi equations of Proposition~\ref{bianchirenor} and the Codazzi equation \eqref{codazzi}. Commuting with $\dkb$ using Propositions \ref{commutation} and \ref{commdkb}, we obtain the following equations for $i=0,1$.
\begin{align*}
 \nabs_3(\dkb^i\a)+\tfrac12\trchb\,\dkb^i\a
 &=-(2-i)\sld_{2-i}^*(\dkb^i\b)+2ir^{-1}\b+\err[\Ga]^{(i)}-\frac{2i(N(r)-1)}r\dkb^i\a,\\
 \nabs_4(\dkb^i\b)+2\trch\,\dkb^i\b
 &=\sld_{2-i}(\dkb^i\a)+\err[\Ga]^{(i)},\\
 \nabs_3(\dkb^i\b)+\trchb\,\dkb^i\b
 &=-r^{-1}\dkb^{i+1}(\rhoc,-\sic)+\err[\Ga]^{(i)}
   -\frac{2i(N(r)-1)}r\dkb^i\b,\\
 \nabs_4\bigl(\dkb^i(\rhoc,-\sic)\bigr)+\tfrac32\trch\,\dkb^i(\rhoc,-\sic)
 &=r^{-1}\dkb^{i+1}\b+\err[\Ga]^{(i)},\\
 \nabs_3\bigl(\dkb^i(\rhoc,\sic)\bigr)+\tfrac32\trchb\,\dkb^i(\rhoc,\sic)
 &=-r^{-1}\dkb^{i+1}\bb+\err[\Ga]^{(i)}
   -\frac{2i(N(r)-1)}r\dkb^i(\rhoc,\sic)+ir\,\Gab\c\Gab\c\Gab\c\Gab,\\
 \nabs_4(\dkb^i\bb)+\trch\,\dkb^i\bb
 &=r^{-1}\dkb^{i+1}(\rhoc,\sic)+\bigl(\Gab\c\W+\Gab\c\slu\bigr)^{(i)}+\bigl(r^{-1}\Gag\c\Ga^{(1)}+\Gag\c\Gag\c\Ga\bigr)^{(i)}.
\end{align*}
Here we split $2ir\K\b$ into $2ir^{-1}\b$ and a remainder in $\err[\Ga]^{(1)}$, using $r\K-1/r\in\Gag^{(1)}$. The outgoing radius corrections are included in the terms involving $\Gag$ by \eqref{bootstrap:good-geometric-coefficients}.
\paragraph{Step 2. Weighted curvature estimates.}
We apply Proposition~\ref{keyintegral} to the four Bianchi pairs above. In each case, \eqref{Gag-quadratic-bulk}, \eqref{schematic-curvature-flux-bulk}, and \eqref{wonderfulrp:separated-products} give $\mathcal E(V)\les\ep^3$. The angular integration-by-parts error in the combined scalar bulk satisfies the same bound by \eqref{bootstrap:angular-densities} and \eqref{wonderfulrp:separated-products}.
For the explicit $N(r)-1$ terms, the multiplier weight $r$ cancels the factor $r^{-1}$. By \eqref{bootstrap:normal-radius-multiplier},
\[
 \int_V|N(r)-1|\,|\W^{(1)}|^2\,d\mu_{\g}
 \les\ep\left(
 \int_V|\W^{(1)}|^2\,d\mu_{\g}
 +\sup_{u'}\int_{C_{u'}\cap V}r|\W^{(1)}|^2\,d\mu_{C_{u'}}
 \right)\les\ep^3.
\]
The nonlinear sources contribute at most $C\ep^3$ by Theorem~\ref{Junkman}.
We sum over $i=0,1$ for each of the last three pairs and combine the four weighted estimates with coefficients $1$, $10$, $2$, and $1$, respectively. Keeping the term $2|\b|^2$ on the right-hand side and dropping excess positive terms gives
\begin{align}
\begin{split}\label{curvature:combined-bianchi}
 &\int_{C_u\cap V}r\left(
 |\a_4|^2+\sum_{i=0}^1
 \bigl(|\dkb^i\a|^2+|\dkb^i\b|^2+|\dkb^i(\rhoc,\sic)|^2\bigr)
 \right)d\mu_{C_u}\\
 &\quad+\int_{\Si_t\cap V}r\Bigg(
 |\a_4|^2+\sum_{i=0}^1\Bigl(
 |\dkb^i\a|^2+|\dkb^i\b|^2+|\dkb^i\bb|^2
 +|\dkb^i(\rhoc,\sic)|^2\Bigr)\Bigg)d\mu_g\\
 &\quad+\int_V\Bigg(
 3|\a_4|^2+14|\as|^2
 +\sum_{i=0}^1\Bigl(
 10|\dkb^i\a|^2+48|\dkb^i\b|^2
 +3|\dkb^i(\rhoc,\sic)|^2+|\dkb^i\bb|^2
 \Bigr)\Bigg)d\mu_{\g}\\
 &\leq C\ep_0^2+C\ep^3
 +\int_V\left(40|\dkb\a|\,|\b|+2|\b|^2+8|\a_4|\,|\a|\right)d\mu_{\g}.
\end{split}
\end{align}
By Young's inequality,
\[
40|\dkb\a|\,|\b|+2|\b|^2+8|\a_4|\,|\a|\leq 9|\dkb\a|^2+47|\b|^2+2|\a_4|^2+8|\a|^2.
\]
Since $\as=\dkb\a$, absorption leaves the positive bulk coefficients $15$, $1$, $1$, and $2$ for $|\dkb\a|^2$, $|\b|^2$, $|\a_4|^2$, and $|\a|^2$, respectively.
Applying \eqref{wLpestimates} with $p=2$ and $w=a$ or $w=\phi a$, and bounding the weight errors by \eqref{bootstrap:angular-densities} and \eqref{Gamma-quartic-flux-bulk}, we obtain
\begin{align*}
\int_{\Si_{t'}\cap V}r|r\nabs(\a,\b,\rhoc,\sic,\bb)|^2\,d\mu_g&\les\int_{\Si_{t'}\cap V}r|\dkb^{\leq1}(\a,\b,\rhoc,\sic,\bb)|^2\,d\mu_g+\ep^4,\\
\int_V|r\nabs(\a,\b,\rhoc,\sic,\bb)|^2\,d\mu_{\g}&\les\int_V|\dkb^{\leq1}(\a,\b,\rhoc,\sic,\bb)|^2\,d\mu_{\g}+\ep^4.
\end{align*}
On the outgoing cones, the unweighted sphere Hodge estimate and \eqref{wonderfulrp:coarea} give the corresponding angular derivative bound for $\W$. Taking suprema over the subregions $V_{t'}(u')\subseteq V$ implies
\begin{align}\label{Fluxest}
\sup_u\int_{C_u\cap V}r|\W^{(1)}|^2\,d\mu_{C_u}+\sup_t\int_{\Si_t\cap V}r|\Wb^{(1)}|^2\,d\mu_g+\int_V|\Wb^{(1)}|^2\,d\mu_\g\les\ep_0^2.
\end{align}
\paragraph{Step 3. The sphere estimate.}
The outgoing Bianchi equations in Proposition~\ref{bianchirenor} and the definition of $\a_4$ give
\begin{equation}\label{curvature:outgoing-derivative}
r|\nabs_4\W|\les|\W^{(1)}|+r|\Gag|\,|\W|+|\Gag|\,|\Gag^{(1)}|+r|\Gag|^3.
\end{equation}
Let $S\subset C_u\cap V$ of radius $R$. Combining the cone trace estimate \eqref{bootstrap:cone-trace}, the outgoing flux in \eqref{Fluxest}, and \eqref{curvature:outgoing-derivative}, we infer
\begin{align}\label{curvature:annular-trace-energy}
R^3\|\W\|_{L^4(S)}^2&\les\II_{(0)}^2+\int_{C_u\cap V\cap\{R/2\leq r\leq2R\}}r\bigl(|\W^{(1)}|^2+|r\nabs_4\W|^2\bigr)\,d\mu_{C_u}\les\ep_0^2,
\end{align}
where we used \eqref{GagGabGawpointwise} for the lower-order terms. Taking the supremum over the cuts gives
\begin{equation}\label{curvature:sphere-L4}
\sup_{S(t,u)\subset V}r^3\|\W\|_{L^4(S)}^2\les\ep_0^2.
\end{equation}
Combining \eqref{Fluxest} and \eqref{curvature:sphere-L4}, we obtain $\mr\les\ep_0$. This completes the proof of Theorem~\ref{M1}.
\section{\texorpdfstring{3D--elliptic estimates (Theorem~\ref{M2})}{Maximal connection estimates}}\label{seck}
We work under the hypotheses of Theorem~\ref{M2} on the complete maximal slices of $V_{t_*}:=\bigcup_{0\leq t\leq t_*}\Si_t$. All constants are independent of $t_*$. We set
\[
\vphi:=\log\phi,\qquad G:=(\eps,\de,\nab\vphi),\qquad i_Nk:=N\c k.
\]
We write $Y=zN^\flat+y$ with $(z,y)=(\de,\eps)$ for $Y=i_Nk$ and $(z,y)=(N\vphi,\nabs\vphi)$ for $Y=\nab\vphi$. We write $\nabs_N$ for the projected derivative in the $N$ direction.
\subsection{Basic equations}
\begin{prop}
The div--curl system for $i_Nk$ is given schematically by
\begin{align}
\begin{split}\label{normal:divcurl-sources}
 \sdiv(i_Nk)&=-r^{-1}\de+\Ga\c\Ga,\\
 (\curl(i_Nk))_N&=\sic,\\
 (\curl(i_Nk))_A&=r^{-1}({}^*\eps)_A+\tfrac12{}^*(\b-\bb)_A+\Ga\c\Ga.
\end{split}
\end{align}
We also have the lapse equation:
\begin{equation}\label{lapse:divcurl}
\sdiv(\nab\vphi)=\De\vphi=|k|^2-|\nab\vphi|^2,\qquad\curl(\nab\vphi)=0.
\end{equation}
On each cut $S(t,u)$, the components of $k$ satisfy
\begin{equation}\label{keyelliptic}
\sdivs\kah+\tfrac12\nabs\de+\tfrac12\tr\th\,\eps=\tfrac12(\b-\bb)+\thh\c\eps.
\end{equation}
\end{prop}
\begin{proof}
By the product rule, we obtain
\[
 \sdiv(i_Nk)=k_{ij}\nab^iN^j,\qquad
 (\curl(i_Nk))_i=H_{iN}+{\in_i}^{jm}k_{m\ell}\nab_jN^\ell.
\]
Substituting $(\nab_A N)_B=\th_{AB}$ and $(\nab_N N)_A=-\nabs_A\log a$, the linear part of $\th$ gives $-r^{-1}\de$ in the divergence and $r^{-1}{}^*\eps$ in the tangential curl. The remaining terms involving $\nab N$ in these two equations are of the form $\Ga\c\Ga$. The magnetic terms follow from \eqref{EHdecomposition}, and the normal curl equals $\sic$ by \eqref{6.6} and \eqref{renorq}.

Equation~\eqref{lapse:divcurl} follows from \eqref{lapseequation} with $\vphi=\log\phi$. Adding the two equations in \eqref{codazzi} and using \eqref{6.6} gives \eqref{keyelliptic}.
\end{proof}
For a spatial one-form $Y=zN^\flat+y$, taking the trace and the antisymmetric part of $\nab Y$ gives
\begin{equation}\label{normal:split-divcurl}
\begin{aligned}
\sdiv Y&=Nz+\sdivs y+\tr\th\,z+\nabs\log a\c y,\\
{\in_A}^{B}(\curl Y)_B&=(\nabs_N y)_A-\nabs_A z+\th_A{}^B y_B-z\nabs_A\log a,\\
(\curl Y)_N&=\curls y.
\end{aligned}
\end{equation}
\subsection{Lapse normalization and radius errors}
\begin{prop}[Lower bound for the lapse]\label{lapse:uniform-lower-bound}
On $V_{t_*}$ we have $1-C\ep\leq\phi\leq1$.
\end{prop}
\begin{proof}
The maximum principle gives $0<\phi\leq1$. Fix $t$ and denote the cut of area radius $s$ by $S_s$ and the region it encloses by $B_s$. By \eqref{lapse:divcurl} and the divergence theorem,
\[
 \int_0^\infty|\overline{N\vphi}(s)|\,ds
 \les\int_0^\infty s^{-2}\int_{B_s}|\Ga|^2d\mu_g\,ds
 =\int_{\Si_t}r^{-1}|\Ga|^2d\mu_g\les\ep^2.
\]
By \eqref{bootstrap:radius-geometry}, $\|\vphi-\overline\vphi\|_{L^\infty(S)}\les\ep$. Lemma~\ref{dint} gives
\[
-\pr_u r=\frac r2\overline{a\tr\th}\simeq\ov{a},\qquad\quad
\pr_r\ov\vphi=\frac{\ov{aN\vphi}+\ov{a\tr\th(\vphi-\ov\vphi)}}{-\pr_u r}.
\]
In particular, $|\partial_r\overline\vphi|\les\ep/r$ by \eqref{wonderfulrp:good-first-sphere}. Fix a sufficiently small universal $c_0>0$ and set
\begin{equation}\label{dfEt}
E_t=\left\{r>0:r^{1/2}\|a^{1/2}\eta\|_{L^4(S_r)}>c_0\overline{a^{1/2}}\right\}.
\end{equation}
By \eqref{normal:sqrt-lapse-oscillation}, $a\simeq\ov{a}$ for $r\notin E_t$, and $\sup_{S_r}a\les r\|a^{1/2}\eta\|_{L^4(S_r)}^2$ for $r\in E_t$. Thus $|dr|\simeq\overline a\,|du|$ and \eqref{Gab-quadratic-flux-bulk} give
\begin{equation}\label{Etdrrest}
 \int_{E_t}|\partial_r\overline\vphi|\,dr
 \les\ep\int_{E_t}\frac{dr}{r}\les\ep\int du\,\|a^{1/2}\eta\|_{L^4(S(t,u))}^2\les\ep^3.
\end{equation}
For $r\notin E_t$, we have $a\simeq\overline a$. We write the numerator of $\partial_r\overline\vphi$ as
\[
\begin{aligned}
 \overline{aN\vphi}+\overline{a\tr\th(\vphi-\overline\vphi)}
 &=\overline a\,\overline{N\vphi}
 +\overline{(a-\overline a)N\vphi}
 +\frac2r\overline{(a-\overline a)(\vphi-\overline\vphi)}+\overline{a(\tr\th-2/r)(\vphi-\overline\vphi)}.
\end{aligned}
\]
By sphere Poincar\'e and \eqref{bootstrap:angular-densities},
\[
 \|a-\ov{a}\|_{L^2(S_r)}\les r\ov{a}\|\Ga\|_{L^2(S_r)},\qquad
 \|\vphi-\ov{\vphi}\|_{L^2(S_r)}\les r\|\Gag\|_{L^2(S_r)}.
\]
Since $|S_r|=4\pi r^2$, Cauchy--Schwarz bounds the last three terms by $C\overline a\,r^{-1}\|\Ga\|_{L^2(S_r)}^2$. Dividing by $-\partial_u r\simeq\overline a$ gives
\[
|\partial_r\overline\vphi|\les|\overline{N\vphi}|+r^{-1}\|\Ga\|_{L^2(S_r)}^2.
\]
On these cuts, $|dr|\simeq a\,|du|$, so coarea and \eqref{Gag-quadratic-bulk} give
\[
 \int_{r\notin E_t}r^{-1}\|\Ga\|_{L^2(S_r)}^2\,dr
 \les\int_{\Si_t}r^{-1}|\Ga|^2\,d\mu_g\les\ep^2.
\]
Combining the estimates on $E_t$ and its complement gives $\int_0^\infty|\partial_r\overline\vphi|\,dr\les\ep^2$, which implies $|\overline\vphi|\les\ep^2$. Combining with $r|\nabs\vphi|\les\ep$, we deduce $|\vphi|\les\ep$.
\end{proof}
\begin{prop}[Lower-order radius terms]\label{normal:radius-errors}
We have
\begin{equation}\label{normal:radius-flux-error}
\begin{aligned}
 \sup_t\int_{\Si_t}|N(r)-1|
 \left(r^{-1}|G|^2+|G|\,|\nabs G|\right)d\mu_g&\les\ep^3,\\
 \int_{V_{t_*}}|N(r)-1|r^{-2}|G|^2d\mu_{\g}&\les\ep^3.
\end{aligned}
\end{equation}
\end{prop}
\begin{proof}
Since $r^{-1}|G|^2+|G|\,|\nabs G|\les r^{-1}|\Gag^{(1)}|^2$, the first estimate in \eqref{normal:radius-flux-error} follows by applying \eqref{bootstrap:normal-radius-multiplier-flux} componentwise to $r^{-1/2}\Gag^{(1)}$ and using \eqref{Gag-quadratic-bulk} and \eqref{wonderfulrp:good-first-sphere}. Next, applying \eqref{bootstrap:normal-radius-multiplier} with $F=r^{-1}G$, we obtain
\begin{align*}
\int_{V_{t_*}}|N(r)-1|r^{-2}|G|^2d\mu_{\g}\les\ep\left(\int_{V_{t_*}}r^{-2}|G|^2d\mu_\g+\sup_u\int_{C_u\cap V_{t_*}}r^{-1}|G|^2d\mu_{C_u}\right).
\end{align*}
The first term is bounded by \eqref{Gag-quadratic-bulk}. By H\"older, \eqref{normal:exterior-sobolev-trace}, and $\phi\simeq1$, the second term satisfies
\[
\sup_u\int_{C_u\cap V_{t_*}}r^{-1}|G|^2d\mu_{C_u}\les\int_0^{t_*}\sup_u\|G\|_{L^4(S(t,u))}^2dt
\les\int_{V_{t_*}}|\nab G|^2d\mu_{\g}\les\ep^2.
\]
This gives the second estimate in \eqref{normal:radius-flux-error}.
\end{proof}
\subsection{Hardy inequalities}
All integrations by parts on the complete slices are justified by first introducing annular cutoffs. The weighted $L^2$ bootstrap bounds for $\Ga^{(2)}$, together with the radius estimates above, ensure that the corresponding boundary terms vanish along suitable sequences of inner and outer radii. In particular, no boundary contributions remain when the cutoffs are removed.
\begin{lem}[Angular spectral gap]\label{normal:angular-gap}
Let $z$ be a scalar, $y\in\sk_1$, and $U\in\sk_2$, satisfying
\begin{equation}\label{normal:abstract-angular-constraint}
 \sdivs U+\tfrac12\nabs z+r^{-1}y=B.
\end{equation}
Then, we have
\begin{equation}\label{normal:quantitative-angular-gap}
 \int_S\left(|\nabs z-r^{-1}y|^2+|\nabs y+r^{-1}z\slg|^2\right)d\slg
 \geq (2-C\ep)r^{-2}\|(z,y)\|_{L^2(S)}^2-C\|B\|_{L^2(S)}^2.
\end{equation}
Moreover, for any $y\in\sk_1$, we have
\begin{equation}\label{normal:oneform-gap}
 \int_S|\nabs y|^2d\slg
 \geq(1-C\ep)r^{-2}\int_S|y|^2d\slg.
\end{equation}
\end{lem}
\begin{proof}
By scaling, it suffices to take $r=1$. By Gauss--Bonnet, we can solve
\begin{equation*}
    \Des f=\K-1,\qquad\quad \ov{f}=0.
\end{equation*}
Proposition~\ref{Lpestimates} gives
$$
\|f\|_{H^2(S)}+\|f\|_{L^\infty(S)}\les\ep.
$$
By \cite[Lemma~2.1]{KlainermanSzeftelGCM2}, the metric $\widetilde g:=e^{2f}\slg$ has Gauss curvature $\K_{\widetilde{g}}=e^{-2f}=1+O(\ep)$ and area radius $\widetilde{r}=1+O(\ep)$. Applying \cite[Corollary~3.8]{KlainermanSzeftelGCM2} to $\widetilde g$ gives a parametrization $\Phi:\mathbb S^2\to S$ in which $\Phi^*\widetilde g=\widetilde r^2e^{2u}\slg_0$, where $\slg_0$ is the unit round metric and $\|u\|_{L^\infty}\les\ep$. Thus
\[
\Phi^*\slg=e^{2\psi}\slg_0,\qquad\psi=u-f\circ\Phi+\log\widetilde r,\qquad\|\psi\|_{L^\infty}\les\ep.
\]
The curvature equation
$\Des_{\slg_0}\psi=1-(\K\circ\Phi)e^{2\psi}$ and the bound for $\K-1$ then give $\|\psi\|_{H^2(\slg_0)}\les\ep$. We now work in this parametrization and omit the pullback notation.

Let $x^j$ be the restrictions of the Cartesian coordinates to $\mathbb S^2$, and let $X_j=\grad_{\slg_0}x^j$. Writing $\nabs^0$ for the round connection, we have $(\nabs^0)^2x^j=-x^j\slg_0$ and hence $\mathcal L_{X_j}\slg=2(X_j\psi-x^j)\slg$. Testing \eqref{normal:abstract-angular-constraint} against $X_j$ and integrating by parts, we obtain
\[
\begin{aligned}
\int_S B(X_j)\,d\slg&=-\frac12\int_S U\c\mathcal L_{X_j}\slg\,d\slg-\frac12\int_S z\,\sdivs X_j\,d\slg+\int_S y(X_j)\,d\slg\\
&=\int_S\bigl(zx^j+y(X_j)-zX_j\psi\bigr)\,d\slg.
\end{aligned}
\]
The second equality uses trace-freeness of $U$ and $\sdivs X_j=2(X_j\psi-x^j)$.

Define the $\mathbb R^3$-valued field $W$ on $\mathbb S^2$ by $W^j=zx^j+y(X_j)$. Since $d\slg=e^{2\psi}d\mu_{\slg_0}$, the preceding identity and the bounds for $\psi$ give
\[
|\overline W_0|\les\|B\|_{L^2(S)}+\ep\|(z,y)\|_{L^2(S)},
\]
where $\overline W_0$ is the round mean of $W$. For any tangent vector field $Y$ on $\mathbb S^2$, we have
\[
 Y(W^j)=(Yz-y(Y))x^j+(\nabs_Y^0y)(X_j)+z\slg_0(Y,X_j).
\]
Applying the round Poincar\'e inequality componentwise gives
\[
\int_{\mathbb S^2}\left(|\nabs^0z-y|^2+|\nabs^0y+z\slg_0|^2\right)d\mu_{\slg_0}\geq 2\|(z,y)\|_{L^2(\slg_0)}^2-8\pi|\ov{W}_0|^2.
\]
The difference of the two connections satisfies
\[
\|\nabs y-\nabs^0y\|_{L^2(\slg_0)}\les\|\nabs^0\psi\|_{L^4(\slg_0)}\|y\|_{L^4(\slg_0)}\les\ep\left(\|\nabs^0y+z\slg_0\|_{L^2(\slg_0)}+\|(z,y)\|_{L^2(\slg_0)}\right).
\]
By Young's inequality and $\|\psi\|_{L^\infty}\les\ep$, we obtain
\[
\int_S\left(|dz-y|_{\slg}^2+|\nabs y+z\slg|_{\slg}^2\right)d\slg\geq (1-C\ep)\int_{\mathbb S^2}\left(|dz-y|_{\slg_0}^2+|\nabs^0y+z\slg_0|_{\slg_0}^2\right)d\mu_{\slg_0}-C\ep\|(z,y)\|_{L^2(\slg_0)}^2.
\]
The squared $L^2$ norms for $\slg$ and $\slg_0$ differ by factors $1+O(\ep)$. Combining these estimates and restoring $r$ gives \eqref{normal:quantitative-angular-gap}. Finally, the first eigenvalue of the rough Laplacian on one-forms on the unit round sphere is $1$, so
\[
\int_{\mathbb S^2}|\nabs^0 y|^2d\mu_{\slg_0}\geq\int_{\mathbb S^2}|y|^2d\mu_{\slg_0}.
\]
Comparing the connections and measures using $\|\psi\|_{L^\infty}+\|\psi\|_{H^2}\les\ep$ gives \eqref{normal:oneform-gap}.
\end{proof}
\begin{prop}[Hardy inequalities]\label{normal:hardy-flux}
For $Y=i_Nk$ or $Y=\nab\vphi$, put $F=r^{1/2}Y$. Then
\begin{align}
 \frac94\int_{\Si_t}r^{-2}|F|^2d\mu_g
 &\leq\int_{\Si_t}|\nab F|^2d\mu_g
 +C\int_{\Si_t}r|(\b,\bb)|^2d\mu_g+C\ep^3,
 &&Y=i_Nk,\label{normal:flux-hardy}\\
 \frac94\int_{\Si_t}r^{-2}|F|^2d\mu_g
 &\leq\int_{\Si_t}|\nab F|^2d\mu_g+C\ep^3,
 &&Y=\nab\vphi.\label{normal:lapse-flux-hardy}
\end{align}
Moreover,
\begin{align}
\frac94\int_{V_{t_*}}r^{-2}|i_Nk|^2d\mu_{\g}&\leq\int_{V_{t_*}}|\nab(i_Nk)|^2d\mu_{\g}+C\int_{V_{t_*}}|(\b,\bb)|^2d\mu_{\g}+C\ep^3,\label{normal:bulk-hardy}\\
\frac14\int_{V_{t_*}}r^{-2}|Y|^2d\mu_{\g}&\leq\int_{V_{t_*}}|\nab Y|^2d\mu_{\g}+C\ep^3,\qquad Y=i_Nk,\nab\vphi.\label{normal:radial-bulk-hardy}
\end{align}
\end{prop}
\begin{proof}
We first work with annular cutoffs and take limits after the integrations by parts. For a compactly supported spatial one-form $F$ away from the axis, integrating $\sdiv(r^{-1}|F|^2N)$ gives
\begin{equation}\label{normal:unweighted-hardy}
\begin{aligned}
 \int_{\Si_t}|\nab_NF|^2d\mu_g
 &=\int_{\Si_t}\left|\nab_NF+\frac{F}{2r}\right|^2d\mu_g
 +\frac14\int_{\Si_t}r^{-2}|F|^2d\mu_g\\
 &\quad-\frac12\int_{\Si_t}\left(
 \frac{N(r)-1}{r^2}-\frac{1}{r}\left(\tr\th-\frac{2}{r}\right)\right)|F|^2d\mu_g.
\end{aligned}
\end{equation}
For $F=r^{1/2}Y$, the last integral in \eqref{normal:unweighted-hardy} is $O(\ep^3)$ by \eqref{normal:radius-flux-error}, \eqref{Gag-quadratic-bulk}, and $r|\Gag|\les\ep$. Thus
\[
\frac14\int_{\Si_t}r^{-2}|F|^2d\mu_g\leq\int_{\Si_t}|\nab_NF|^2d\mu_g+C\ep^3.
\]
For the bulk, integrate $\sdiv(\phi r^{-1}|Y|^2N)$ on each slice and then in time. The additional coefficient is $N\vphi/r$, and the same bounds give
\[
\frac14\int_{V_{t_*}}r^{-2}|Y|^2d\mu_{\g}\leq\int_{V_{t_*}}|\nab_NY|^2d\mu_{\g}+C\ep^3.
\]
This implies \eqref{normal:radial-bulk-hardy}.

For $Y=zN^\flat+y$, the angular derivatives satisfy
\begin{equation}\label{normal:angular-components}
 (\nab_AY)_N=\nabs_Az-\th_A{}^By_B,\qquad(\nab_AY)_B=\nabs_Ay_B+z\th_{AB}.
\end{equation}
For $Y=i_Nk$, apply Lemma~\ref{normal:angular-gap} to $\varpi^{1/2}(\de,\eps,\kah)$, taking $\varpi=a$ on a slice and $\varpi=\phi a$ for the bulk. By \eqref{keyelliptic} and \eqref{bootstrap:angular-densities}, its source is schematically $\varpi^{1/2}(\tfrac12(\b-\bb)+\Ga\c\Ga)$. Replacing $r^{-1}\slg$ by $\th$ and differentiating the weight gives
\begin{equation}\label{normal:angular-hardy}
(2-C\ep)r^{-2}\int_S\varpi|Y|^2d\slg\leq\int_S\varpi\sum_A|\nab_AY|^2d\slg+C\int_S\varpi\left(|(\b,\bb)|^2+\ep^{-1}|\Ga|^4\right)d\slg.
\end{equation}
The $C\ep$ loss and the quartic errors contribute $O(\ep^3)$ after integration, by \eqref{Gag-quadratic-bulk} and \eqref{null:full-quartic-bootstrap}. Multiplying by $r$ and integrating in $u$, we obtain
\begin{equation}\label{normal:angular-flux-coercivity}
2\int_{\Si_t}r^{-2}|F|^2d\mu_g\leq\int_{\Si_t}\sum_A|\nab_AF|^2d\mu_g+C\int_{\Si_t}r|(\b,\bb)|^2d\mu_g+C\ep^3,
\end{equation}
since $\nab_AF=r^{1/2}\nab_AY$. Adding the radial $1/4$ proves \eqref{normal:flux-hardy}. Integrating \eqref{normal:angular-hardy} with $\varpi=\phi a$ in $u,t$ and adding the radial bulk estimate proves \eqref{normal:bulk-hardy}.

For $Y=\nab\vphi=zN^\flat+y$, the mixed curl equation in \eqref{normal:split-divcurl} gives
\[
y\c\nabs z=\tfrac12N(|y|^2)+\th(y,y)-z\,y\c\nabs\log a.
\]
Integrating the normal derivative by parts, we find
\[
\int_{\Si_t}y\c\nabs z\,d\mu_g=\int_{\Si_t}\left[\left(\th-\tfrac12(\tr\th)\slg\right)(y,y)-z\,y\c\nabs\log a\right]d\mu_g=O(\ep^3).
\]
Sphere integration by parts with $d\mu_g=a\,du\,d\slg$ also gives
\[
\int_{\Si_t}z\sdivs y\,d\mu_g=-\int_{\Si_t}\left(y\c\nabs z+z\,y\c\nabs\log a\right)d\mu_g=O(\ep^3).
\]
Applying \eqref{normal:oneform-gap} to $a^{1/2}y$ and using $\nabs\log a=\eta-\eps$, we obtain
\[
\int_{\Si_t}r|\nabs y|^2d\mu_g\geq(1-C\ep)\int_{\Si_t}r^{-1}|y|^2d\mu_g
-C\ep^{-1}\int_{\Si_t}r|\Ga|^2|y|^2d\mu_g\geq\int_{\Si_t}r^{-1}|y|^2d\mu_g-C\ep^3,
\]
by \eqref{Gag-quadratic-bulk} and \eqref{null:full-quartic-bootstrap}. Expanding \eqref{normal:angular-components} therefore yields
\[
\begin{aligned}
\int_{\Si_t}r\sum_A|\nab_AY|^2d\mu_g&=\int_{\Si_t}\left[r(|\nabs z|^2+|\nabs y|^2)+r^{-1}(|y|^2+2|z|^2)-2y\c\nabs z+2z\sdivs y\right]d\mu_g+O(\ep^3)\\
&\geq 2\int_{\Si_t}r^{-1}|Y|^2d\mu_g-C\ep^3.
\end{aligned}
\]
The geometric errors are controlled by $r|Y|\les\ep$ and the quadratic and quartic bootstrap estimates. Hence
\begin{equation}\label{normal:lapse-angular-coercivity}
 2\int_{\Si_t}r^{-2}|F|^2d\mu_g
 \leq\int_{\Si_t}\sum_A|\nab_AF|^2d\mu_g+C\ep^3.
\end{equation}
Adding the radial $1/4$ proves \eqref{normal:lapse-flux-hardy}.
\end{proof}
\subsection{Unweighted div--curl estimates}\label{ssec:maximal-divcurl}
\begin{prop}\label{normal:spatial-ricci}\label{normal:divcurl-estimate}
Let $F$ have components $F_i$ or $F_{iA}$. The operators $\sdiv$ and $\curl$ act on $i$, and the derivatives on $A$ are projected. In the spatial index write $z_F=F(N)$, $y_F=F|_{TS}$, and set
\begin{equation}\label{dfqFUF}
q_F:=|y_F|^2+2|z_F|^2,\qquad U_F:=\left(\tr\th-\frac2r\right)|z_F|^2+(\th-r^{-1}\slg)_{AB}\langle y_F^A,y_F^B\rangle.
\end{equation}
In a local orthonormal frame $\{e_j\}_{j=1,2,3}$, define $J(F)$ by
\[
g(J(F),X):=\langle F(e_j),(\nab_{e_j}F)(X)\rangle-\langle F(X),\sdiv F\rangle+U_Fg(N,X).
\]
The unweighted Hodge identity is
\begin{equation}\label{Hodgeid}
|\nab F|^2-|\sdiv F|^2-|\curl F|^2=\sdiv J(F)-\frac{N(r)-1}{r^2}q_F-\mathcal E_F-\mathcal K_F,
\end{equation}
where
\[
|\mathcal E_F|\les|\Ga|\,|F|\,|\nab F|+\left(r^{-1}|\Ga^{(1)}|+|\Ga|^2\right)|F|^2.
\]
The term $\mathcal K_F$ vanishes for $F_i$. For $F_{iA}$, it is bounded by
\[
 |\mathcal K_F|\les\left(r^{-2}+r^{-1}|\Ga^{(1)}|+|\Ga|^2\right)|F|^2.
\]
The normal component of the current is
\begin{equation}\label{normal:hodge-normal-current}
 g(J(F),N)=\mathcal T_F-r^{-1}q_F,\qquad
 \mathcal T_F:=\langle y_F,\nabs z_F\rangle-\langle z_F,\sdivs y_F\rangle.
\end{equation}
Moreover, the integration identity holds:
\begin{equation}\label{normal:hodge-multiplier}
\begin{aligned}
&\int_{\Si_t}w\left(|\nab F|^2-|\sdiv F|^2-|\curl F|^2\right)d\mu_g\\
&\qquad=-\int_{\Si_t}w\left(\frac{N(r)-1}{r^2}q_F+\mathcal E_F+\mathcal K_F\right)d\mu_g-\int_{\Si_t}\nab w\c J(F)\,d\mu_g.
\end{aligned}
\end{equation}
\end{prop}
\begin{proof}
The spatial Gauss and Codazzi identities and the normal curvature equation give \cite[Equations~(2.92), (2.93), (2.103), and (3.43)]{Gourgoulhon07}
\begin{equation}\label{normal:ricci-decomposition}
\begin{aligned}
\sRic_{NN}&=-N(\tr\th)-|\th|^2-\Des\log a-|\nabs\log a|^2,\\
\sRic_{NA}&=(\sdivs\th)_A-\nabs_A\tr\th,\\
\sRic_{AB}&=\K\slg_{AB}-\nabs_N\th_{AB}-(\tr\th)\th_{AB}-\nabs_A\nabs_B\log a-(\nabs_A\log a)(\nabs_B\log a).
\end{aligned}
\end{equation}
Using $N(r^{-1})=-r^{-2}N(r)$ and separating the Minkowski parts of $\th$ and $\K$, we obtain
\begin{equation}\label{normal:ricci-remainder}
\begin{aligned}
\sRic(F,F)&=\frac{N(r)-1}{r^2}q_F-N\left(\tr\th-\frac2r\right)|z_F|^2-\bigl(\nabs_N(\th-r^{-1}\slg)\bigr)_{AB}\langle y_F^A,y_F^B\rangle\\
&\quad+O\left(\left(r^{-1}|\Ga^{(1)}|+|\Ga|^2\right)|F|^2\right)\\
&=\frac{N(r)-1}{r^2}q_F-\sdiv(U_FN)+\mathcal E_F.
\end{aligned}
\end{equation}
The first error uses the angular derivatives of $\th-r^{-1}\slg$ and $\nabs\log a$, and \eqref{bootstrap:gauss-coefficient}. In the last line, the product rule gives terms bounded by $|\Ga|\,|F|(|\nab F|+|\Ga|\,|F|)$, and $(\sdiv N)U_F$ is bounded by $(r^{-1}|\Ga|+|\Ga|^2)|F|^2$.

In the standard pointwise Hodge identity \cite[Lemma~4.4.1]{Ch-Kl}, the current is $J(F)-U_FN$ and the Ricci term is $-\sRic(F,F)$. For $F_{iA}$, commuting derivatives on $A$ adds $-\mathcal K_F$. The angular curvature is $\K$, and the mixed curvature is bounded by $|\nabs\th|+|\th|\,|\nabs\log a|$. Thus the curvature of the projected connection satisfies
\begin{equation}\label{normal:coefficient-curvature}
 |R^{\mathrm{ang}}|\les r^{-2}+r^{-1}|\Ga^{(1)}|+|\Ga|^2.
\end{equation}
Substitution of \eqref{normal:ricci-remainder} gives \eqref{Hodgeid}. For the normal component of the current, the $e_j=N$ term cancels the normal derivative in $\sdiv F$. The remaining terms give
\[
g(J(F),N)=\langle y_F,\nabs z_F\rangle-\langle z_F,\sdivs y_F\rangle-\th_{AB}\langle y_F^A,y_F^B\rangle-(\tr\th)|z_F|^2+U_F=\mathcal T_F-r^{-1}q_F.
\]
For a scalar $f$, the terms with derivatives of $f$ cancel in $J(fF)$. The Ricci contraction scales quadratically, so \eqref{normal:ricci-remainder} also gives
\begin{equation}\label{normal:rescaling-identities}
J(fF)=f^2J(F),\qquad\mathcal E_{fF}=f^2\mathcal E_F+N(f^2)U_F.
\end{equation}
Finally, integration of \eqref{Hodgeid} against a scalar weight $w$ gives \eqref{normal:hodge-multiplier}.
\end{proof}
\begin{cor}\label{normal:geometric-errors}
For $Y=i_Nk$ or $Y=\nab\vphi$,
\begin{equation}\label{normal:weighted-hodge-bulk}
\left|\int_{V_{t_*}}\left(|\nab Y|^2-|\sdiv Y|^2-|\curl Y|^2\right)d\mu_{\g}\right|\les\ep^3.
\end{equation}
\end{cor}
\begin{proof}
Take $w=\phi$ in \eqref{normal:hodge-multiplier} and integrate in time using $d\mu_{\g}=\phi\,dt\,d\mu_g$. Since $r|Y|\les\ep$ and $|J(Y)|\les|Y|\,|\nab Y|+|\Ga|\,|Y|^2$, we have
\[
|\mathcal E_Y|+|\nab\vphi|\,|J(Y)|\les\ep\left(|\nab Y|^2+r^{-2}|\Ga^{(1)}|^2\right).
\]
Use \eqref{Gag-quadratic-bulk} and \eqref{secknorms} for these terms and \eqref{normal:radius-flux-error} for the radius term.
\end{proof}
\begin{prop}\label{normal:rescaled-divcurl}
For $F$ as in Proposition~\ref{normal:spatial-ricci} and $p=2,3$, we have
\begin{equation}\label{normal:commuted-hodge-estimate}
 \int_{\Si_t}r^p|\nab F|^2d\mu_g\les\int_{\Si_t}r^p\left(|\sdiv F|^2+|\curl F|^2\right)d\mu_g+\int_{\Si_t}\left(r^{p-2}+r^{p-1}|\Ga^{(1)}|+r^p|\Ga|^2\right)|F|^2d\mu_g.
\end{equation}
\end{prop}
\begin{proof}
Apply \eqref{bootstrap:normal-radius-multiplier-flux} to $r^{(p-2)/2}F$. By H\"older and \eqref{normal:weighted-L4-trace}, the cut term satisfies
\[
r^{p-1}\|F\|_{L^2(S)}^2\les r^p\|F\|_{L^4(S)}^2\les\int_{\Omega_{t,u}}r^p|\nab F|^2d\mu_g.
\]
Taking the supremum over the cuts gives
\begin{equation}\label{normal:high-radius-quadratic}
\int_{\Si_t}r^{p-2}|N(r)-1|\,|F|^2d\mu_g\les\ep\int_{\Si_t}\left(r^p|\nab F|^2+r^{p-2}|F|^2\right)d\mu_g.
\end{equation}
Let $E_t$ be defined in \eqref{dfEt}. On $r\notin E_t$, \eqref{e3e4r} and $a\simeq\ov a$ give $N(r)\les1$. Since $|\mathcal T_F|\les |F|(|\nab F|+(r^{-1}+|\Ga|)|F|)$, Young's inequality gives
\[
\begin{aligned}
\int_{\{r\notin E_t\}}r^{p-1}N(r)|\mathcal T_F|\,d\mu_g&\les\int_{\Si_t}\left[r^{p-1}|F|\,|\nab F|+(r^{p-2}+r^{p-1}|\Ga|)|F|^2\right]d\mu_g\\
&\leq\de_0\int_{\Si_t}r^p|\nab F|^2d\mu_g+C_{\de_0}\int_{\Si_t}(r^{p-2}+r^p|\Ga|^2)|F|^2d\mu_g.
\end{aligned}
\]
The same calculation bounds $\int_{\Si_t}r^{p-1}|\mathcal T_F|\,d\mu_g$. Thus, for the estimate with $N(r)-1$, we can include the term with coefficient $-1$ here, allocating $\de_0/2$ to each contribution.

To estimate the contribution on $E_t$, integrate \eqref{Hodgeid} over the exterior $\Omega_{t,u}$. Its inner boundary has outward normal $-N$, so by \eqref{normal:hodge-normal-current},
\[
\int_S\mathcal T_F\,d\slg=r(t,u)^{-1}\int_S q_F\,d\slg-\int_{\Omega_{t,u}}\left(|\nab F|^2-|\sdiv F|^2-|\curl F|^2+\frac{N(r)-1}{r^2}q_F+\mathcal E_F+\mathcal K_F\right)d\mu_g.
\]
Multiply by $r(t,u)^p$. The boundary term is bounded by $Cr(t,u)^{p-1}\|F\|_{L^2(S)}^2$, which is controlled by the trace estimate above. In the volume integral, $r(t,u)^p\leq r^p$. Thus the div--curl squares are bounded by $Cr^p|\nab F|^2$, the radius term is controlled by \eqref{normal:high-radius-quadratic}, and the mixed term from $\mathcal E_F$ satisfies
\[
r^p|\Ga|\,|F|\,|\nab F|\les r^p|\nab F|^2+r^p|\Ga|^2|F|^2.
\]
The remaining terms in $\mathcal E_F$ and $\mathcal K_F$ have the weights displayed below. We obtain
\begin{equation}\label{normal:unweighted-boundary-trace}
 r(t,u)^p\left|\int_S\mathcal T_F\,d\slg\right|\les\int_{\Si_t}r^p|\nab F|^2d\mu_g+\int_{\Si_t}\left(
 r^{p-2}+r^{p-1}|\Ga^{(1)}|+r^p|\Ga|^2
 \right)|F|^2d\mu_g.
\end{equation}
Since $aN(r)=-\partial_u r$, changing the coarea variable from $u$ to $r$ gives $N(r)d\mu_g=dr\,d\slg$ as positive measures. Taking the sphere integral first, we have
\[
\left|\int_{\{r\in E_t\}}r^{p-1}N(r)\mathcal T_F\,d\mu_g\right|=\left|\int_{E_t}r^{p-1}\left(\int_{S_r}\mathcal T_F\,d\slg\right)dr\right|\leq\left(\int_{E_t}\frac{dr}{r}\right)\sup_u r(t,u)^p\left|\int_{S(t,u)}\mathcal T_F\,d\slg\right|.
\]
The first factor satisfies $\int_{E_t}r^{-1}dr\les\ep^2$ by \eqref{Etdrrest}. Using \eqref{normal:unweighted-boundary-trace} for the second factor and adding the preceding Young bound gives
\begin{equation}\label{normal:high-radius-form-bound}
\begin{aligned}
 \left|\int_{\Si_t}r^{p-1}N(r)\mathcal T_F\,d\mu_g\right|
 &\leq(\de_0+C\ep^2)\int_{\Si_t}r^p|\nab F|^2d\mu_g\\
 &\quad+C_{\de_0}\int_{\Si_t}
 \left(r^{p-2}+r^{p-1}|\Ga^{(1)}|+r^p|\Ga|^2\right)|F|^2d\mu_g.
\end{aligned}
\end{equation}
The same bound holds with $N(r)-1$ in place of $N(r)$.

Now take $w=r^p$ in \eqref{normal:hodge-multiplier}. By \eqref{normal:hodge-normal-current}, the radius term and the normal current contribute
\[
\int_{\Si_t}r^{p-2}\bigl(p+(p-1)(N(r)-1)\bigr)q_F\,d\mu_g-p\int_{\Si_t}r^{p-1}N(r)\mathcal T_F\,d\mu_g.
\]
Use \eqref{normal:high-radius-quadratic} and \eqref{normal:high-radius-form-bound} for the terms containing $N(r)-1$ and $N(r)\mathcal T_F$. In the geometric remainder, Young's inequality bounds $r^p|\Ga|\,|F|\,|\nab F|$ by $\de_0 r^p|\nab F|^2+C_{\de_0}r^p|\Ga|^2|F|^2$. The other terms already have the weights on the right of \eqref{normal:commuted-hodge-estimate}. Choosing $\de_0$ and then $\ep$ sufficiently small absorbs the derivative terms.
\end{proof}
\subsection{First-order estimates}
\begin{prop}[First-order estimates]\label{normal:epsilon-delta-estimates}
We have
\begin{align}
 \sup_t\int_{\Si_t}\left(r^{-1}|G|^2+r|\nab G|^2\right)d\mu_g
 &\les\sup_t\int_{\Si_t}r|(\b,\sic,\bb)|^2d\mu_g+\ep^3\les\ep_0^2,
 \label{normal:epsilon-delta-flux}\\
 \int_{V_{t_*}}\left(r^{-2}|G|^2+|\nab G|^2\right)d\mu_{\g}
 &\les\int_{V_{t_*}}|(\b,\sic,\bb)|^2d\mu_{\g}+\ep^3\les\ep_0^2.
 \label{normal:epsilon-delta-bulk}
\end{align}
\end{prop}
\begin{proof}
For $F=r^{1/2}i_Nk$, the product rule and \eqref{normal:divcurl-sources} give
\begin{align*}
&\bigl(\sdiv F,{\in_A}^{B}(\curl F)_B,(\curl F)_N\bigr)-\frac{N(r)-1}{2r}(F_N,F_A,0)\\
&\qquad=-\frac1{2r}(F_N,F_A,0)+r^{1/2}\bigl(0,-\tfrac12(\b-\bb)_A,\sic\bigr)+O(r^{1/2}|\Ga|^2).
\end{align*}
For $F=r^{1/2}\nab\vphi$, \eqref{lapse:divcurl} gives
\begin{equation*}
\bigl(\sdiv F,{\in_A}^{B}(\curl F)_B,(\curl F)_N\bigr)-\frac{N(r)-1}{2r}(F_N,F_A,0)=\frac1{2r}(F_N,F_A,0)+O(r^{1/2}|\Ga|^2).
\end{equation*}
Insert \eqref{normal:unweighted-hardy} into the unweighted integral of \eqref{Hodgeid}. Expanding gives
\[
\left|\nab_NF+\frac{F}{2r}\right|^2-\left|\nab_NF+\frac{2-N(r)}{2r}F\right|^2
=\frac{N(r)-1}{r}\left(\langle\nab_NF,F\rangle+\frac{|F|^2}{2r}\right)-\frac{(N(r)-1)^2}{4r^2}|F|^2.
\]
Expanding the corrected div--curl square gives
\[
\begin{aligned}
&|\sdiv F|^2+|\curl F|^2-\left|\bigl(\sdiv F,{\in_A}^{B}(\curl F)_B,(\curl F)_N\bigr)-\frac{N(r)-1}{2r}(F_N,F_A,0)\right|^2\\
&\qquad=\frac{N(r)-1}{r}\left(F_N\sdiv F+F^A{\in_A}^{B}(\curl F)_B\right)-\frac{(N(r)-1)^2}{4r^2}\bigl(|z_F|^2+|y_F|^2\bigr).
\end{aligned}
\]
The radial square enters with a plus sign and the div--curl square with a minus sign. The quadratic corrections therefore cancel pointwise:
\[
 -\frac{(N(r)-1)^2}{4r^2}|F|^2
 +\frac{(N(r)-1)^2}{4r^2}\bigl(|z_F|^2+|y_F|^2\bigr)=0.
\]
For the linear corrections, \eqref{normal:split-divcurl} and
\eqref{normal:hodge-normal-current} give
\[
\langle\nab_NF,F\rangle-F_N\sdiv F-F^A{\in_A}^{B}(\curl F)_B=g(J(F),N)-U_F=\mathcal T_F-r^{-1}q_F-U_F.
\]
Thus, combining \eqref{Hodgeid} and \eqref{normal:unweighted-hardy} with the two square expansions yields
\[
\begin{aligned}
 &\int_{\Si_t}\Bigg[
 \sum_A|\nab_AF|^2
 +\left|\nab_NF+\frac{2-N(r)}{2r}F\right|^2
 +\frac{|F|^2}{4r^2}\\
 &\hspace{25mm}
 -\left|
 \bigl(\sdiv F,{\in_A}^{B}(\curl F)_B,(\curl F)_N\bigr)
 -\frac{N(r)-1}{2r}(F_N,F_A,0)
 \right|^2
 \Bigg]d\mu_g\\
 &=\int_{\Si_t}\Bigg[
 -\frac{N(r)-1}{r^2}q_F-\mathcal E_F
 +\frac{N(r)-1}{2r^2}|F|^2
 -\frac{\tr\th-2/r}{2r}|F|^2\\
 &\hspace{25mm}
 -\frac{N(r)-1}{r}
 \left(\mathcal T_F-r^{-1}q_F-U_F+\frac{|F|^2}{2r}\right)
 \Bigg]d\mu_g\\
 &=-\int_{\Si_t}\left[
 \frac{N(r)-1}{r}\mathcal T_F
 +r\mathcal E_{r^{-1/2}F}
 +r^{-1}U_F+\frac{\tr\th-2/r}{2r}|F|^2
 \right]d\mu_g
 =O(\ep^3),
\end{aligned}
\]
where we used \eqref{normal:rescaling-identities} in the form
$\mathcal E_F-\frac{N(r)-1}{r}U_F=r\mathcal E_{r^{-1/2}F}+r^{-1}U_F$.
Since $\mathcal T_F=r\mathcal T_{r^{-1/2}F}$, the first term in the last integral is controlled by \eqref{normal:radius-flux-error}. The other terms are controlled by $r|r^{-1/2}F|\les\ep$, \eqref{Gag-quadratic-bulk}, and \eqref{secknorms}. We perform the cancellations before removing the annular cutoffs. By Young's inequality and \eqref{null:full-quartic-bootstrap}, the corrected sources give
\begin{equation}\label{normal:slice-hodge-reduced}
\begin{aligned}
 &\int_{\Si_t}\left(
 \sum_A|\nab_AF|^2+\left|\nab_NF+\frac{2-N(r)}{2r}F\right|^2
 +\frac{|F|^2}{4r^2}\right)d\mu_g\\
 &\qquad\leq\frac12\int_{\Si_t}r^{-2}|F|^2d\mu_g
 +C\int_{\Si_t}r|(\b,\sic,\bb)|^2d\mu_g+C\ep^3.
\end{aligned}
\end{equation}
By \eqref{normal:angular-flux-coercivity} and \eqref{normal:lapse-angular-coercivity}, the angular term and the retained radial Hardy term control $\int_{\Si_t}r^{-2}|F|^2d\mu_g$ with coefficient $9/4$. This absorbs the coefficient $1/2$ on the right. Since
\[
\nab_NF+\frac{2-N(r)}{2r}F=r^{1/2}\left(\nab_NY+r^{-1}Y\right),
\]
the flux estimates for $i_Nk$ and $\nab\vphi$ follow by the product rule. For the bulk, \eqref{normal:weighted-hodge-bulk} with the sources in \eqref{normal:divcurl-sources} and \eqref{lapse:divcurl} gives
\[
\begin{aligned}
\int_{V_{t_*}}|\nab(i_Nk)|^2d\mu_{\g}&\leq\tfrac32\int_{V_{t_*}}r^{-2}|i_Nk|^2d\mu_{\g}+C\int_{V_{t_*}}|(\b,\sic,\bb)|^2d\mu_{\g}+C\ep^3,\\
\int_{V_{t_*}}|\nab(\nab\vphi)|^2d\mu_{\g}&\les\ep^3.
\end{aligned}
\]
Apply \eqref{normal:bulk-hardy} to $i_Nk$ and \eqref{normal:radial-bulk-hardy} to $\nab\vphi$. Finally, $|G|^2=|i_Nk|^2+|\nab\vphi|^2$. The component derivatives of $(\de,\eps)$ differ from those of $i_Nk$ by terms bounded by $(r^{-1}+|\Ga|)|i_Nk|$. Their squares are controlled by the field estimates and \eqref{null:full-quartic-bootstrap}. Taking the supremum in $t$ for the flux and using $\mr\les\ep_0$ and $\ep^3=\ep_0^2$ gives both conclusions.
\end{proof}
\subsection{Second-order estimates}\label{ssec:maximal-angular-estimate}
For each field $Y=zN^\flat+y$, set
\[
 F=(\nabs z)N^\flat+\nabs y.
\]
\begin{prop}\label{normal:commutation}
The commuted system satisfies
\begin{equation}\label{maximal:commuted-system}
|\sdiv F|+|\curl F|\les|\nabs(\b,\sic,\bb)|+r^{-1}|F|+r^{-2}|Y|+r^{-1}|\Ga|\,|\Ga^{(1)}|+|\Ga|^3.
\end{equation}
\end{prop}
\begin{proof}
Differentiate \eqref{normal:split-divcurl} angularly. In the divergence and mixed curl equations, use
\begin{equation}\label{maximal:normal-angular-commutators}
\begin{aligned}
 {}[\nabs_N,\nabs_A]z&=-\th_A{}^B\nabs_Bz+(\nabs_A\log a)Nz,\\
 [\nabs_N,\nabs_A]y_B
 &=-\th_A{}^C\nabs_Cy_B+(\nabs_A\log a)\nabs_Ny_B
 +\bigl(r^{-1}\Ga^{(1)}+\Ga\c\Ga\bigr)\c y.
\end{aligned}
\end{equation}
These follow from $\nab_AN=\th_A{}^Be_B$, $(\nab_NN)_A=-\nabs_A\log a$, and the mixed curvature bound used in \eqref{normal:coefficient-curvature}. Commuting angular derivatives contributes $\K\c y$. The undifferentiated equations bound $(Nz,\nabs_Ny)$ by $C(|F|+r^{-1}|Y|+|(\b,\bb)|+|\Ga|^2)$. Differentiating the sources in \eqref{normal:divcurl-sources} and \eqref{lapse:divcurl}, and using $r(\b,\sic,\bb)\subseteq\Ga^{(1)}$ and \eqref{bootstrap:gauss-coefficient}, gives \eqref{maximal:commuted-system}.
\end{proof}
\begin{prop}[Mixed second derivatives]\label{normal:mixed-second-reduction}
For $(z,y)=(\de,\eps)$ and $(N\vphi,\nabs\vphi)$,
\begin{equation}\label{maximal:component-second-estimate}
\begin{aligned}
&\sup_t\int_{\Si_t}r^3\left(|\nabs^2(z,y)|^2+|\nabs_N\nabs(z,y)|^2\right)d\mu_g\\
&\quad+\int_{V_{t_*}}r^2\left(|\nabs^2(z,y)|^2+|\nabs_N\nabs(z,y)|^2\right)d\mu_{\g}
 \les\ep_0^2.
\end{aligned}
\end{equation}
\end{prop}
\begin{proof}
In \eqref{maximal:commuted-system}, the terms $r^{-1}F$ and $r^{-2}Y$ use the first-order estimates, and the curvature derivatives use $\mr\les\ep_0$. The products use \eqref{null:full-quartic-bootstrap} and \eqref{maximal:all-sixth}. Hence
\begin{equation}\label{maximal:commuted-source-estimate}
\begin{aligned}
 &\sup_t\int_{\Si_t}\left[r|F|^2+r^3\left(|\sdiv F|^2+|\curl F|^2\right)\right]d\mu_g\\
 &\quad+\int_{V_{t_*}}\left[|F|^2+r^2\left(|\sdiv F|^2+|\curl F|^2\right)\right]d\mu_{\g}
 \les\ep_0^2.
\end{aligned}
\end{equation}
Since $r|F|\les|\Ga^{(1)}|$, the remaining lower-order terms obey
\begin{equation}\label{normal:second-order-product-bounds}
\sup_t\int_{\Si_t}r^3\left(r^{-1}|\Ga^{(1)}|+|\Ga|^2\right)|F|^2d\mu_g+\int_{V_{t_*}}r^2\left(r^{-1}|\Ga^{(1)}|+|\Ga|^2\right)|F|^2d\mu_{\g}
 \les\ep^3.
\end{equation}
Here the cubic term uses Cauchy--Schwarz between \eqref{Gag-quadratic-bulk} and \eqref{null:full-quartic-bootstrap}, and the quartic term uses \eqref{null:full-quartic-bootstrap}. Apply Proposition~\ref{normal:rescaled-divcurl} with $p=3$ on the slices and with $p=2$ before integrating in time, using $\phi\simeq1$. By \eqref{maximal:commuted-source-estimate}, \eqref{normal:second-order-product-bounds}, and $\ep^3=\ep_0^2$,
\[
\sup_t\int_{\Si_t}r^3|\nab F|^2d\mu_g+\int_{V_{t_*}}r^2|\nab F|^2d\mu_{\g}\les\ep_0^2.
\]
The terms from passing between $\nab F$ and the projected derivatives are bounded by $(r^{-1}+|\Ga|)|F|$. Their weighted squares are covered by the same two estimates, proving \eqref{maximal:component-second-estimate}.
\end{proof}
Combining Propositions \ref{normal:epsilon-delta-estimates} and \ref{normal:mixed-second-reduction}, we obtain
\begin{equation}\label{normal:complete-mixed-estimate}
\begin{aligned}
 &\sup_t\int_{\Si_t}r^{-1}\left(|G|^2+r^2|\nab G|^2+r^4|\nab\nabs G|^2\right)d\mu_g\\
 &\quad+\int_{V_{t_*}}r^{-2}\left(|G|^2+r^2|\nab G|^2+r^4|\nab\nabs G|^2\right)d\mu_{\g}
 \les\ep_0^2.
\end{aligned}
\end{equation}
Apply \eqref{normal:weighted-L4-trace} with $p=1$ to $G$ and with $p=3$ to $\nabs G$. The sphere Sobolev inequality then gives
\begin{equation}\label{maximal:improved-W14}
 \sup_{t,u}\left(r^{1/2}\|G\|_{W^{1,4}(S)}+r\|G\|_{L^\infty(S)}\right)\les\ep_0.
\end{equation}
By \eqref{secknorms}, the two bounds give $\mk\les\ep_0$.
\section{\texorpdfstring{Null transport and sphere elliptic estimates (Theorem~\ref{M3})}{Null connection estimates}}\label{secO}
We work on $V=V_{t_*}$ under the hypotheses of Theorem~\ref{M3}. In the cone integrals, $S=S(t',u)$ denotes a sphere in $V$, and the supremum in $t'$ is taken at fixed $u$.
\subsection{Quadratic estimates for \texorpdfstring{$(\Wb,\slu,\mu)$}{}}
\label{ssec:transport-main-equations}\label{ssec:noncritical-transport}
We apply Lemma~\ref{evolution} to rescaled fields. By \eqref{bootstrap:good-geometric-coefficients}, the rescaling of \eqref{evolution:equation} is
\begin{equation}\label{transport:rescale}
\nabs_4(rX)+\left(\la-\tfrac12\right)\trch\,(rX)=rF+\Gag\c(rX).
\end{equation}
The definition of $\a_4$ gives
\begin{equation}\label{transport:curvature-equations}
 \nabs_4(r\a)+2\trch\,(r\a)=\a_4+\Gag\c(r\a).
\end{equation}
Using \eqref{mu-transport}, \eqref{U-transport}, and \eqref{transport:rescale}, we obtain
\begin{equation}\label{transport:connection-equations}
\begin{aligned}
\nabs_4(r\slu)+\trch\,(r\slu)&=-2\b+\Gag\c\Gag^{(1)}+r\Gag\c\Gag\c\Gag,\\
\nabs_4(r\mu)+\tfrac12\trch\,(r\mu)&=-\sdivs\etab+\rhoc+r\Gab\c(\b,\slu)+\Gag\c\Ga^{(1)}+r\Gag\c\Gag\c\Ga.
\end{aligned}
\end{equation}
Proposition \ref{bianchirenor} gives
\begin{equation}\label{transport:rescaled-betab}
 \nabs_4(r\bb)+\tfrac12\trch\,(r\bb)
 =r\sld_1^*(\rhoc,\sic)+r\Gab\c(\W,\slu)
   +\Gag\c\Ga^{(1)}+r\Gag\c\Gag\c\Ga.
\end{equation}
We will also use the outgoing Codazzi equation in the form
\begin{equation}\label{transport:codazzi}
 \sdivs\hch=\tfrac12\slu-\b-\hch\c\eps,
\end{equation}
which follows from \eqref{codazzi}, \eqref{6.6}, and \eqref{dfmu}.
\begin{prop}\label{transport:noncritical-estimates}
We have
\begin{equation}\label{transport:noncritical-quadratic}
\begin{aligned}
 &\int du\,\sup_{t'}\int_S ar\left(|\Wb|^2+|(r\nabs)^{\leq1}(\slu,\mu)|^2\right)d\slg\\
 &\quad+\int_V\left(|\Wb|^2+|(r\nabs)^{\leq1}(\slu,\mu)|^2
 \right)d\mu_{\g}+\sup_{t,u}r^2\|\slu\|_{H^1(S)}^2\les\ep_0^2.
\end{aligned}
\end{equation}
\end{prop}
\begin{proof}
Combining Proposition~\ref{bianchirenor} with \eqref{transport:curvature-equations}--\eqref{transport:rescaled-betab}, we have the unified transport form
\begin{align*}
    \nabs_4X+\lambda_X\trch\,X=F_X,
\end{align*}
with the coefficients $\lambda_X$ given by
\begin{equation}\label{transport:unified-quadratic-system}
\begin{array}{c|cccccc}
 X&r\a&r\b&r(\rhoc,-\sic)&r\slu&r\mu&r\bb\\ \hline
 \lambda_X&2&\frac32&1&1&\frac12&\frac12
\end{array}.
\end{equation}
The corresponding sources are schematically given by
\begin{equation}\label{transport:unified-quadratic-sources}
\begin{aligned}
F_{r\a}&=\a_4+r\Gag\c\W,\\
F_{r\b}&=r\sld_2\a+r\Gag\c\W,\\
F_{r(\rhoc,-\sic)}&=r\sld_1\b+\Gag\c\Gag^{(1)}+r\Gag\c\Gag\c\Gag,\\
F_{r\slu}&=-2\b+\Gag\c\Gag^{(1)}+r\Gag\c\Gag\c\Gag,\\
F_{r\mu}&=-\sdivs\etab+\rhoc+r\Gab\c(\b,\slu)+\Gag\c\Ga^{(1)}+r\Gag\c\Gag\c\Ga,\\
F_{r\bb}&=r\sld_1^*(\rhoc,\sic)+r\Gab\c(\W,\slu)+\Gag\c\Ga^{(1)}+r\Gag\c\Gag\c\Ga.
\end{aligned}
\end{equation}
The linear sources have $L^2(V)$ norm bounded by $C(\mr+\mk)$, where the term $\sdivs\etab$ is controlled by $\mk$ using $\etab=\nabs\log\phi-\eps$ from \eqref{6.6}. The other terms satisfy
\begin{equation}\label{transport:unified-nonlinear-bulk}
\begin{aligned}
\|r\Gag\c\W\|_{L^2(V)}&\les\ep\mr,\\
\|\Gag\c\Ga^{(1)}\|_{L^2(V)}+\|r\Gag\c\Gag\c\Ga\|_{L^2(V)}&\les\ep^2,\\
\|r\Gab\c(\W,\slu)\|_{L^2(V)}&\les\|\Gab\|_{L^4(V)}\|r(\W,\slu)\|_{L^4(V)}\les\ep^2.
\end{aligned}
\end{equation}
Here we used $r|\Gag|\les\ep$,
\eqref{Gag-quadratic-bulk}, and
\eqref{Gamma-quartic-flux-bulk}, together with
$r(\W,\slu)\subseteq\Gag^{(1)}$. In \eqref{evolution:Lp} with $p=2$, Young's inequality bounds the source integral by
\[
\de_0\int_Vr^{-2}|X|^2d\mu_{\g}+C_{\de_0}\int_V|F_X|^2d\mu_{\g}.
\]
Applying Lemma \ref{evolution} to $X$ and absorbing the first term, we obtain
\begin{equation}\label{transport:unified-quadratic-bound}
 \sum_X\left[
 \int du\,\sup_{t'}\int_S ar^{-1}|X|^2d\slg
 +\int_Vr^{-2}|X|^2d\mu_{\g}\right]
 \les\II_{(0)}^2+\mr^2+\mk^2+\ep^4\les\ep_0^2.
\end{equation}
\paragraph{Estimate for $\slu$.}
Set $Y_{\slu}=(r\nabs)(r\slu)$. By \eqref{U-transport} and \eqref{bootstrap:commute4}, its equation takes the schematic form
\begin{equation}\label{transport:unified-commuted-U}
\nabs_4Y_{\slu}+\trch\,Y_{\slu}=-2(r\nabs)\b+\Gag\c Y_{\slu}+r\Gag\c\W^{(1)}+\hch\c(r\nabs)^2\hch+\Gag^{(1)}\c\Gag^{(1)}.
\end{equation}
By \eqref{null:full-quartic-bootstrap}, \eqref{Gag-quadratic-bulk}, and $r|\Gag|\les\ep$, we have
\[
\|\hch\c(r\nabs)^2\hch\|_{L^2(V)}+\|\Gag^{(1)}\c\Gag^{(1)}\|_{L^2(V)}+\|r\Gag\c\W^{(1)}\|_{L^2(V)}\les\ep^2.
\]
The term $\Gag\c Y_{\slu}$ contributes a small multiple of the transport bulk. Thus \eqref{evolution:Lp} gives
\begin{equation}\label{transport:unified-commuted-U-bound}
 \int du\,\sup_{t'}\int_S ar^{-1}|Y_{\slu}|^2d\slg
 +\int_Vr^{-2}|Y_{\slu}|^2d\mu_{\g}
 \les\II_{(0)}^2+\mr^2+\ep^4\les\ep_0^2.
\end{equation}
By \eqref{transport:codazzi} and Proposition~\ref{Lpestimates},
\begin{equation}\label{transport:unified-U-Hodge}
 r^2\|\nabs^2\hch\|_{L^2(S)}
 \les\|Y_{\slu}\|_{L^2(S)}+r\|\W^{(1)}\|_{L^2(S)}+\ep.
\end{equation}
Using this bound for $\hch\c(r\nabs)^2\hch$, and
\eqref{wonderfulrp:good-first-sphere} for the remaining products, we have
\[
\begin{aligned}
 \|F_{r\slu}+2\b\|_{L^2(S)}
 &\les\frac{\ep^2}{r},\\
 \|F_{Y_{\slu}}+2(r\nabs)\b\|_{L^2(S)}
 &\les\frac{\ep}{r}\|Y_{\slu}\|_{L^2(S)}
 +\ep\|\W^{(1)}\|_{L^2(S)}+\frac{\ep^2}{r},
\end{aligned}
\]
where $F_{r\slu}$ and $F_{Y_{\slu}}$ are the sources in
\eqref{transport:connection-equations} and \eqref{transport:unified-commuted-U}. By Lemma~\ref{dint} and $e_4r=1+O(\ep)$, the two transport equations give
\begin{equation}\label{transport:unified-U-sphere-ODE}
\frac{d}{dr}\|r\slu\|_{H^1(S)}+\frac{1-C\ep}{r}\|r\slu\|_{H^1(S)}\les\|\W^{(1)}\|_{L^2(S)}+\frac{\ep^2}{r}.
\end{equation}
Since $1-C\ep\geq\tfrac12$, integration with the factor $r^{1/2}$ and Cauchy--Schwarz give
\begin{equation}\label{sluH2est}
r\|\slu\|_{H^1(S)}\les\II_{(0)}+\ep^2+\left(\int_{C_u\cap V}r|\W^{(1)}|^2d\mu_{C_u}\right)^{1/2}\les\ep_0.
\end{equation}
The initial-cut term uses $\II_{(0)}\les\ep_0$, and the vertex term vanishes.
\paragraph{Estimate for $\mu$.}
Set $Y_\mu=(r\nabs)(r\mu)$. Commuting \eqref{transport:connection-equations} by \eqref{bootstrap:commute4} gives
\begin{equation}\label{transport:unified-commuted-mu}
\begin{aligned}
 \nabs_4Y_\mu+\tfrac12\trch\,Y_\mu&=-(r\nabs)\sdivs\etab+r\nabs\rhoc+\Gag\c Y_\mu\\
 &\quad+\Ga^{(1)}\c\Ga^{(1)}+\Gag\c\Ga^{(2)}+\Gab\c(r\nabs)(r\b,r\slu).
\end{aligned}
\end{equation}
The linear sources have $L^2(V)$ norm bounded by $C(\mk+\mr)$ using \eqref{normal:complete-mixed-estimate} for the term involving $\etab$. The first three products are handled as in \eqref{transport:unified-commuted-U}. For the last product use \eqref{bootstrap:coefficient-product}, taking the cone flux for $(r\nabs)(r\b)$ and the sphere supremum for $Y_{\slu}$. This gives
\begin{equation}\label{transport:mu-bad-products}
 \int_V|\Gab|^2|(r\nabs)(r\b,r\slu)|^2d\mu_{\g}
 \les\ep^2\left(\mr^2+\sup_{t,u}r^2\|\slu\|_{H^1(S)}^2\right)
 \les\ep^2\ep_0^2.
\end{equation}
Applying \eqref{evolution:Lp} to $Y_\mu$ with $p=2$ yields the first angular estimate for $\mu$. Together with \eqref{transport:unified-quadratic-bound} and \eqref{sluH2est}, this proves Proposition \ref{transport:noncritical-estimates}.
\end{proof}
\subsection{Angular estimates for \texorpdfstring{$\Ga$}{}}
\label{ssec:null-hodge-traces}
The incoming Codazzi and torsion equations, together with \eqref{6.6}, \eqref{renorq}, and \eqref{dfmu}, give
\begin{align}
\sdivs\hchb&=\bb+\nabs\de-\tfrac12\slu+(\trch-\de)\eps+\eps\c\hchb,\label{transport:incoming-codazzi}\\
\sdivs\eta&=-\mu-\rhoc,\qquad\qquad\curls\eta=\sic.\label{transport:eta-hodge}
\end{align}
\begin{prop}\label{null:connection-hodge-estimates}
We have
\begin{equation}\label{null:good-sphere-estimates}
 \sup_{t,u}
 r^{1/2}\left\|
 \left(\trch-\frac{2}{r},\hch,\trchb+\frac{2}{r}\right)
 \right\|_{W^{1,4}(S)}\les\ep_0
\end{equation}
and
\begin{equation}\label{null:augmented-quadratic-bulk}
 \sup_t\int_{\Si_t}r^{-1}|\Ga^{(2)}|^2d\mu_g
 +\int_Vr^{-2}|\Ga^{(2)}|^2d\mu_{\g}\les\ep_0^2.
\end{equation}
\end{prop}
\begin{proof}
By \eqref{bootstrap:mass-and-trace},
\begin{equation}\label{null:trace-gradient}
\nabs\trch=\slu-\trch\,\eps.
\end{equation}
By \eqref{null:trace-gradient} and \eqref{transport:codazzi}, the sphere bounds for $\slu$, $\b$, and $\eps$ give
\[
r^{1/2}\|(\trchc,\hch)\|_{W^{1,4}(S)}\les\ep_0,\qquad\quad \trchc:=\trch-\ov{\trch}.
\]
Averaging \eqref{gauss} and using Gauss--Bonnet, \eqref{renorq}, and $\trchb=2\de-\trch$ gives
\begin{equation}\label{null:mean-trace-Gauss}
\left(\overline{\trch}-\frac2r\right)\left(\overline{\trch}+\frac2r-2\overline\de\right)=4\overline{\rhoc}+\frac4r\overline\de-\overline{\trchc^2}+2\overline{\de\,\trchc}.
\end{equation}
The second factor is comparable to $r^{-1}$, and the right-hand side
is $O(\ep_0r^{-2})$ by the sphere bounds for $\rhoc$, $\trchc$, and $\de$. Thus $r|\overline{\trch}-2/r|\les\ep_0$. Combining with $\trchb+2/r=2\de-(\trch-2/r)$ gives \eqref{null:good-sphere-estimates}.

Take $w=a$ on the slices and $w=\phi a$ in spacetime. Apply the $H^2(S)$ Hodge estimates to $w^{1/2}(\eta,\hch,\hchb)$ using \eqref{transport:codazzi},
\eqref{transport:incoming-codazzi}, and \eqref{transport:eta-hodge}.
Propositions~\ref{Lpestimates} and~\ref{Hodge2Diden}, together with \eqref{bootstrap:angular-densities}, give
\[
\begin{aligned}
 r^{-2}\left\|w^{1/2}(r\nabs)^{\leq2}
     (\eta,\hch,\hchb)\right\|_{L^2(S)}^2
 &\les\left\|w^{1/2}(r\nabs)^{\leq1}
     (\mu,\rhoc,\sic,\slu,\b,\bb)\right\|_{L^2(S)}^2\\
 &\quad+r^{-2}\left\|w^{1/2}(r\nabs)^{\leq2}
     (\eps,\de)\right\|_{L^2(S)}^2+\int_S w\left(|\Ga^{(1)}|^4+r^2|\Ga|^6\right)d\slg.
\end{aligned}
\]
Multiply by $r$ and integrate in $u$ for the slice bound, and integrate
in $(t,u)$ for the spacetime bound. For $(\mu,\slu)$ use
Proposition~\ref{transport:noncritical-estimates}, and for $(\eps,\de)$
use \eqref{normal:complete-mixed-estimate}.
The last term contributes $O(\ep^4+\ep^6)$ by
\eqref{null:full-quartic-bootstrap} and \eqref{maximal:all-sixth}.
This controls the contributions of $(\eta,\hch,\hchb)$ to
\eqref{null:augmented-quadratic-bulk}, and also gives the same bound for
$\kah=-\tfrac12(\hch+\hchb)$.

For the expansions, \eqref{bootstrap:weighted-poincare} and
\eqref{bootstrap:weighted-mean} applied to \eqref{null:trace-gradient}
and \eqref{null:mean-trace-Gauss} give
\[
\begin{aligned}
 \|w^{1/2}\trchc\|_{L^2(S)}
 &\les\|w^{1/2}(r\slu,\eps)\|_{L^2(S)}
      +\ep_0 r^{1/2}\|w^{1/2}\eta\|_{L^4(S)},\\
 \left\|w^{1/2}\left(\overline{\trch}-\tfrac2r\right)\right\|_{L^2(S)}
 &\les\|w^{1/2}(r\rhoc,\de)\|_{L^2(S)}
      +\ep_0\|w^{1/2}\trchc\|_{L^2(S)}
      +\ep_0 r^{1/2}\|w^{1/2}\eta\|_{L^4(S)}.
\end{aligned}
\]
The squared weight errors contribute $O(\ep_0^2\ep^2)$ by
\eqref{Gab-quadratic-flux-bulk} and \eqref{wonderfulrp:weighted-mixed}.
For the angular derivatives, use
\[
 (r\nabs)\trchc=r\slu-(r\trch)\eps,\qquad
 (r\nabs)^2\trchc
 =(r\nabs)(r\slu)-(r\trch)(r\nabs)\eps
 -r\bigl((r\nabs)\trchc\bigr)\otimes\eps.
\]
The quadratic estimates for $\slu,\rhoc,\eps,\de$ then control
$\trch-2/r$ and $\trchb+2/r$ through two angular derivatives.

The maximal coefficients in $\Gag$ follow from $\mk$ and \eqref{6.6}.
By \eqref{e3e4r},
\[
 \frac{e_4r}{r}-\frac{\trch}{2}
 =\frac{\overline{\phi\trch}-\phi\trch}{2\phi}.
\]
Since $\phi\simeq1$ and
$\|\phi\trch-\overline{\phi\trch}\|_{L^2(S)}\les\ep_0$,
\eqref{bootstrap:weighted-poincare} gives
\[
 \left\|w^{1/2}\left(\frac{e_4r}{r}-\frac{\trch}{2}\right)\right\|_{L^2(S)}
 \les\|w^{1/2}(r\nabs\trch,\nabs\log\phi)\|_{L^2(S)}
     +\ep_0 r^{1/2}\|w^{1/2}\eta\|_{L^4(S)}.
\]
Differentiating the identity bounds its angular derivatives by those
of $\trch$ and $\log\phi$. The same identity and the sphere bounds give
$r|\Gag|\les\ep_0$.

The terms involving $r\mu$ and $r\slu$ use
Proposition~\ref{transport:noncritical-estimates}, and $r\K-r^{-1}$
is controlled by \eqref{bootstrap:gauss-coefficient}.
This proves \eqref{null:augmented-quadratic-bulk}.
\end{proof}
\subsection{Estimates for strong fluxes}
\begin{prop}
We have
\begin{equation}\label{null:incoming-shear-strong-W14}
\int du\,\sup_{t'}\left[\|a^{1/2}(\eta,\hch,\hchb)\|_{W^{1,4}(S)}^2+\|a^{1/2}r(\mu,\bb)\|_{L^4(S)}^2\right]\les\ep_0^2
\end{equation}
and
\begin{equation}\label{null:lower-strong-fluxes}
\int du\,\sup_{t'}\int_S ar^{-1}|(r\nabs)^{\leq1}(\eta,\hch,\hchb)|^2d\slg\les\ep_0^2.
\end{equation}
The outgoing curvature components satisfy
\begin{equation}\label{transport:noncritical-quartic}
\int du\,\sup_{t'}\int_S ar|r\W|^4d\slg+\int_V|r\W|^4d\mu_{\g}\les\ep_0^4.
\end{equation}
\end{prop}
\begin{proof}
By Proposition~\ref{nulles},
\begin{equation}
\begin{aligned}\label{transport:bad-structure-equations}
\nabs_4\eta+\tfrac12\trch\,\eta&=-\b+\tfrac12\trch\,\etab-\hch\c(\eta-\etab),\\
\nabs_4\hchb+\tfrac12\trch\,\hchb&=\nabs\hot\etab-\tfrac12\trchb\,\hch+2\om\hchb+\etab\hot\etab,\\
\nabs_4\hch+\trch\,\hch&=-\a-2\om\hch.
\end{aligned}
\end{equation}
Apply \eqref{bootstrap:cone-trace} at each cut $S_R$ of radius $R$ to
\[
B\in\left\{a^{1/2}(\eta,\hch,\hchb),\quad (r\nabs)\bigl(a^{1/2}(\eta,\hch,\hchb)\bigr),\quad a^{1/2}r(\Wb,\mu)\right\}.
\]
Taking the supremum along the cone and integrating in $dR/R$ give
\[
\begin{aligned}
\int du\,\sup_{t'}\|B\|_{L^4(S)}^2+\int du\,\frac{dr}{r}\,\|B\|_{L^4(S)}^2&\les\int du\,\|B\|_{L^4(S_0)}^2\\
&\quad+\int du\int_{C_u\cap V}\left(r^{-2}|B|^2+|\nabs B|^2+|\nabs_4B|^2\right)d\mu_{C_u}.
\end{aligned}
\]
Indeed, each radius $r$ contributes for $R\in[r/2,2r]$, and the initial term contributes for $R\in[r_0/2,2r_0]$. Both intervals have $dR/R$ measure $\log 4$. The initial term is omitted on cones issuing from the axis.

On $\Si_0$, $a\simeq1$ and $r|\nabs\log a|\les\ep_0$, so the product rule, sphere Sobolev, and slice coarea bound the initial term by $C\ep_0^2$. Using \eqref{transport:bad-structure-equations}, \eqref{transport:connection-equations}, \eqref{transport:rescaled-betab}, and \eqref{curvature:outgoing-derivative}, together with \eqref{bootstrap:commute4}, $e_4\log a=-2\om$, and \eqref{bootstrap:angular-densities}, we obtain, for each choice of $B$,
\[
r^{-1}|B|+|\nabs B|+|\nabs_4 B|\les a^{1/2}\left(r^{-1}|\Ga^{(2)}|+|\Ga^{(1)}|^2+r|\Ga|^3\right).
\]
By coarea, \eqref{null:augmented-quadratic-bulk}, \eqref{null:full-quartic-bootstrap}, and \eqref{maximal:all-sixth}, the cone integral is $O(\ep_0^2+\ep^4+\ep^6)$. This gives \eqref{null:incoming-shear-strong-W14} and
\begin{equation}\label{transport:weighted-curvature-trace}
\int du\,\sup_{t'}r^2\|a^{1/2}(\Wb,\mu)\|_{L^4(S)}^2+\int du\,\frac{dr}{r}\,r^2\|a^{1/2}(\Wb,\mu)\|_{L^4(S)}^2\les\ep_0^2.
\end{equation}
For \eqref{null:lower-strong-fluxes}, the product rule and
$r|(\hch,\eps)|\les\ep$ give
\[
\int_S ar^{-1}|(r\nabs)^{\leq1}(\eta,\hch,\hchb)|^2d\slg\les\|a^{1/2}(\eta,\hch,\hchb)\|_{W^{1,4}(S)}^2+\int_S ar|\Gab|^4d\slg .
\]
Use \eqref{null:incoming-shear-strong-W14} and
\eqref{null:zeroth-strong-quartic-bootstrap} after taking the cone
supremum and integrating in $u$. Finally, the sphere bound for $\W$ gives
\[
\int_S ar|r\W|^4d\slg\leq r^5\|\W\|_{L^4(S)}^2\|a^{1/2}\W\|_{L^4(S)}^2\les\ep_0^2r^2\|a^{1/2}\W\|_{L^4(S)}^2.
\]
The two terms in \eqref{transport:weighted-curvature-trace} bound the cone supremum and the $du\,dr/r$ integral, respectively. By coarea, these give \eqref{transport:noncritical-quartic}.
\end{proof}
\begin{prop}
We have
\begin{align}\label{null:good-quartic-improved}
\sup_t\int_{\Si_t}r|\Gag^{(1)}|^4d\mu_g+\int_V|\Gag^{(1)}|^4d\mu_{\g}&\les\ep_0^4,\\
\int du\,\sup_{t'}\int_S ar|\Gab|^4d\slg+\int_V|\Gab|^4d\mu_{\g}&\les\ep_0^4,
\label{null:bad-zeroth-quartic}\\
\int du\,\sup_{t'}\int_S ar^3|\Gab|^6d\slg+\int_Vr^2|\Gab|^6d\mu_{\g}&\les\ep_0^6.\label{null:bad-sixth-closed}
\end{align}
\end{prop}
\begin{proof}
The zeroth-order quartic terms follow from $r|\Gag|\les\ep_0$ and
\eqref{null:augmented-quadratic-bulk}. For $r\nabs\Gag$, apply
\eqref{bootstrap:angular-interpolation} with $w=a,\phi a$, using
\eqref{null:augmented-quadratic-bulk} for the second angular derivatives
and \eqref{bootstrap:angular-densities} and \eqref{Gamma-quartic-flux-bulk}
for the weight terms. This gives $C(\ep_0^4+\ep_0^4\ep^4)$.
The curvature terms follow from \eqref{transport:noncritical-quartic}
and slice coarea, while $r\slu$ and $r\K-r^{-1}$ use
\eqref{bootstrap:mass-and-trace} and \eqref{bootstrap:gauss-coefficient}.
This proves \eqref{null:good-quartic-improved}.

Keeping $-\hch\c\eta$ and $2\om\hchb$ as small transport terms in
\eqref{transport:bad-structure-equations}, we write the remaining sources as
\[
F_\eta=-\b+\tfrac12\trch\,\etab+\hch\c\etab,\qquad F_{\hchb}=\nabs\hot\etab-\tfrac12\trchb\,\hch+\etab\hot\etab.
\]
By \eqref{null:good-quartic-improved},
\[
\int_V r^4\bigl(|F_\eta|^4+|F_{\hchb}|^4\bigr)d\mu_{\g}\les\ep_0^4.
\]
Using $\phi\simeq1$ and \eqref{normal:weighted-spatial-sobolev} with
$p=3,2$, interpolate the slice weights $r^9$ and $r^6$ with powers
$2/3$ and $1/3$ and integrate in $t$:
\begin{equation}\label{null:source-sixth-interpolation}
\int_Vr^8|Y|^6d\mu_{\g}\les\left[\sup_t\int_{\Si_t}\left(r|Y|^2+r^3|\nab Y|^2\right)d\mu_g\right]^2\int_V\left(|Y|^2+r^2|\nab Y|^2\right)d\mu_{\g}.
\end{equation}
Here $\nab$ is the projected metric connection on horizontal fields.

Using $2N=e_4-e_3$ in the Bianchi equations of
Proposition~\ref{bianchirenor}, with principal angular sources
$r\sld_2\a$ and $r\sld_1^*(\rhoc,-\sic)$, gives
\begin{equation}\label{null:beta-spatial-first-derivative}
 r|\nab\b|\les|\W^{(1)}|+|\Ga^{(1)}|^2+r|\Ga|^3.
\end{equation}
The expansion and outgoing-shear equations of Proposition~\ref{nulles}
likewise give
\[
|\nab\hch|\les r^{-1}|\Ga^{(1)}|+|\Ga|^2,\qquad|\nab\trch|+|\nab\trchb|\les r^{-2}+r^{-1}|\Ga^{(1)}|+|\Ga|^2.
\]
Thus the sources and their spatial derivatives satisfy
\[
|(F_\eta,F_{\hchb})|+r|\nab(F_\eta,F_{\hchb})|\les |\b|+r|\nab\b|+r^{-1}|\Ga^{(1)}|+|\nab\etab|+r|\nab\nabs\etab|+|\Ga^{(1)}|^2+r|\Ga|^3.
\]
The $\etab$ terms use $\etab=\nabs\log\phi-\eps$ and \eqref{normal:complete-mixed-estimate}. By \eqref{null:augmented-quadratic-bulk}, \eqref{null:full-quartic-bootstrap}, and \eqref{maximal:all-sixth}, the quadratic slice and bulk integrals in \eqref{null:source-sixth-interpolation} are $O(\ep_0^2)$ for $Y=F_\eta,F_{\hchb}$. Hence
\begin{equation}\label{null:remaining-sixth-sources}
\int_Vr^8\bigl(|F_\eta|^6+|F_{\hchb}|^6\bigr)d\mu_{\g}\les\ep_0^6.
\end{equation}
Apply \eqref{evolution:Lp} to $\eta$ and $\hchb$ for $p=4,6$. Absorb the small transport terms using $r|(\hch,\om)|\les\ep_0$. For the remaining sources, use
\[
r^{p-3}|X|^{p-1}|F|\leq\de_0 r^{p-4}|X|^p+C_{\de_0}r^{2p-4}|F|^p.
\]
The initial pointwise and quadratic bounds give initial fluxes $O(\ep_0^p)$. The quartic source bound and \eqref{null:remaining-sixth-sources} then give \eqref{null:bad-zeroth-quartic} and \eqref{null:bad-sixth-closed}.
\end{proof}
\begin{prop}\label{transport:betab-quartic-reduction-prop}
We have
\begin{equation}\label{transport:mass-incoming-quartic}
\int du\,\sup_{t'}\int_S ar|r(\mu,\bb)|^4d\slg+\int_V|r(\mu,\bb)|^4d\mu_{\g}\les\ep_0^4.
\end{equation}
\end{prop}
\begin{proof}
Denote the left-hand side of \eqref{transport:mass-incoming-quartic}
by $\mathcal Q$. By \eqref{null:augmented-quadratic-bulk} and
\eqref{null:full-quartic-bootstrap},
\begin{equation}\label{transport:betab-weighted-H1}
\begin{aligned}
 \int du\,dr\,r^{-2}\|a^{1/2}\Gab^{(1)}\|_{H^1(S)}^2\les \int_V\left(
 r^{-2}|\Gab^{(2)}|^2
 +|\nabs\log a|^2|\Gab^{(1)}|^2
 \right)d\mu_{\g}\les\ep_0^2+\ep^4\les\ep_0^2.
\end{aligned}
\end{equation}
Since $a\simeq1$ on $\Si_0$, the weighted trace estimate
\eqref{normal:weighted-L4-trace} with $p=3$ and sphere Sobolev give
\[
\int_{\Si_0}r^5|(\mu,\bb)|^4d\mu_g\les\left(\sup_{S_0}r^3\|(\mu,\bb)\|_{L^4(S_0)}^2\right)\int du\,r^2\|(\mu,\bb)\|_{L^4(S_0)}^2\les\II_{(0)}^4\les\ep_0^4.
\]
Applying \eqref{evolution:Lp} with $p=4$ to the equations for $r\mu$
and $r\bb$ in \eqref{transport:connection-equations} and
\eqref{transport:rescaled-betab}, we obtain
\[
\mathcal Q\les\ep_0^4+\int_V r|r(\mu,\bb)|^3
\left[
|-\sdivs\etab+\rhoc|
+r|\nabs(\rhoc,\sic)|
+|\Gab|\,|r(\W,\slu)|
+|\Gag|\,|\Ga^{(1)}|
+r|\Gag|^2|\Ga|
\right]d\mu_{\g}.
\]
The contribution of $-\sdivs\etab+\rhoc$ is bounded by
$\de_0\mathcal Q+C_{\de_0}\ep_0^4$ using
\eqref{null:good-quartic-improved} and Young's inequality.
The last two products satisfy the same bound by $r|\Gag|\les\ep_0$
and \eqref{null:full-quartic-bootstrap}.

\paragraph{Linear curvature source.}
Applying \eqref{bootstrap:sphere-W11} to $a|r(\mu,\bb)|^3$ gives
\[
 \|a|r(\mu,\bb)|^3\|_{L^2(S)}
 \les
 \left(\int_S a|r(\mu,\bb)|^4d\slg\right)^{1/2}
 \left(\int_S a\left[
 r^{-2}|\Gab^{(2)}|^2+|\nabs\log a|^2|\Gab^{(1)}|^2
 \right]d\slg\right)^{1/2}.
\]
Pair this with $r^2\nabs(\rhoc,\sic)$ on $S$ and take the cone
supremum of $r\int_S a|r(\mu,\bb)|^4d\slg$. Cauchy--Schwarz in $r$
pairs the outgoing flux of $r^{3/2}\nabs(\rhoc,\sic)$ with the bulk in
\eqref{transport:betab-weighted-H1}. Taking the curvature flux uniformly
in $u$ and applying Cauchy--Schwarz in $u$ yields
\begin{equation}\label{transport:betab-linear-source}
 \int_V r^2|r(\mu,\bb)|^3|\nabs(\rhoc,\sic)|d\mu_{\g}
 \les \ep_0^2\mathcal Q^{1/2}
 \leq \de_0\mathcal Q+C_{\de_0}\ep_0^4.
\end{equation}
\paragraph{Products with outgoing curvature.}
H\"older on $S$ gives
\[
\begin{aligned}
 \int_S a|r(\mu,\bb)|^3|\Gab|\,|\W|d\slg\les \|a^{1/4}r(\mu,\bb)\|_{L^4(S)}^{4/3}
       \|a^{1/4}\Gab^{(1)}\|_{L^4(S)}^{4/3}\|a^{1/2}\Gab^{(1)}\|_{L^6(S)}^{1/3}
       \|a^{1/6}\Gab\|_{L^6(S)}\|\W\|_{L^9(S)}.
\end{aligned}
\]
After multiplication by $r^2$, \eqref{bootstrap:sphere-sobolev}
bounds this by
\[
\begin{aligned}
 &\left(\sup_{t'}\int_S ar|r(\mu,\bb)|^4d\slg\right)^{1/3}
 \left(\sup_{t'}\int_S ar^3|\Gab|^6d\slg\right)^{1/6}\\
 &\quad\times
 \left(\int_S a|\Gab^{(1)}|^4d\slg\right)^{1/3}
 \left(r^{-2}\|a^{1/2}\Gab^{(1)}\|_{H^1(S)}^2\right)^{1/6}
 \left(r\|\W\|_{H^1(S)}^2\right)^{1/2}.
\end{aligned}
\]
Apply H\"older in $r$ with exponents $3,6,2$ to the last three
factors. Take the resulting curvature flux uniformly in $u$ and
apply H\"older in $u$ with exponents $3,6,3,6$ to the other factors.
By \eqref{null:bad-sixth-closed},
\eqref{null:full-quartic-bootstrap}, and \eqref{transport:betab-weighted-H1},
\begin{equation}\label{transport:sixth-curvature-product}
 \int_V r^2|r(\mu,\bb)|^3|\Gab|\,|\W|d\mu_{\g}
 \les \ep_0^{7/3}\ep^{4/3}\mathcal Q^{1/3}
 \leq \de_0\mathcal Q+C_{\de_0}\ep_0^4,
\end{equation}
where the last inequality uses Young's inequality and $\ep^3=\ep_0^2$.
\paragraph{Products with $\slu$.}
Use sphere H\"older with exponents $4/3,6,12$ on
$|a^{1/4}\Gab^{(1)}|^3$, $a^{1/6}|\Gab|$, and $a^{1/12}|\slu|$.
By \eqref{bootstrap:sphere-sobolev} and the sphere bound in
\eqref{transport:noncritical-quadratic},
\[
 \|a^{1/12}\slu\|_{L^{12}(S)}
 \leq \|\slu\|_{L^{12}(S)}^{5/6}
       \|a^{1/2}\slu\|_{L^{12}(S)}^{1/6}
 \les \ep_0^{5/6}r^{-5/3}\|a^{1/2}\slu\|_{H^1(S)}^{1/6}.
\]
The quadratic estimate for $\slu$ and the quartic bounds for
$\nabs\log a$ and $r\slu$ give
\begin{equation}\label{transport:U-weighted-H1}
 \int du\,dr\,\|a^{1/2}\slu\|_{H^1(S)}^2
 \les \int_V\left(
 |\slu|^2+|r\nabs\slu|^2+r^2|\nabs\log a|^2|\slu|^2
 \right)d\mu_{\g}\les\ep_0^2.
\end{equation}
The remaining factor $r^{1/3}$ combines with
$\|a^{1/6}\Gab\|_{L^6(S)}$ to give
$\bigl(\int_S ar^2|\Gab|^6d\slg\bigr)^{1/6}$.
H\"older in $(u,r)$ with exponents $4/3,6,12$, followed by
\eqref{null:bad-sixth-closed} and \eqref{transport:U-weighted-H1}, yields
\begin{equation}\label{transport:sixth-U-product}
\int_V r^2|\Gab^{(1)}|^3|\Gab|\,|\slu|d\mu_{\g}\les \ep_0^2\|\Gab^{(1)}\|_{L^4(V)}^3\les \ep_0^2\ep^3=\ep_0^4.
\end{equation}
Choosing $\de_0$ sufficiently small and absorbing the terms
proportional to $\mathcal Q$ proves \eqref{transport:mass-incoming-quartic}.
\end{proof}
Combining \eqref{transport:noncritical-quartic} and
\eqref{transport:mass-incoming-quartic}, and using slice coarea, gives
\begin{equation}\label{null:curvature-quartic-combined}
\int du\,\sup_{t'}\int_S ar|r(\Wb,\mu)|^4d\slg+\sup_{t'}\int_{\Si_{t'}}r|r(\Wb,\mu)|^4d\mu_g+\int_V|r(\Wb,\mu)|^4d\mu_{\g}\les\ep_0^4.
\end{equation}
The quadratic terms in \eqref{sssec8.1.3} follow from
\eqref{transport:noncritical-quadratic}, \eqref{null:good-sphere-estimates}, \eqref{null:incoming-shear-strong-W14}, and \eqref{null:augmented-quadratic-bulk}. The quartic terms follow from \eqref{null:curvature-quartic-combined} and \eqref{null:bad-zeroth-quartic}, and the sixth-power terms from \eqref{null:bad-sixth-closed} and slice
coarea. Thus $\mo\les\ep_0$.
\section{Conclusions}\label{sec:conclusions}
Throughout this section, we assume that the main estimates in Theorem \ref{mainth} hold.
\subsection{Future geodesic completeness}
\begin{prop}\label{conclusion:geodesic-completeness}
The future development $(\M,\g)$ is future
causally geodesically complete. In particular, every future-inextendible timelike or null geodesic has infinite affine length.
\end{prop}
\begin{proof}
Let $\gamma(\lambda)$ be a future-directed causal
geodesic, and write
\[
\dot\gamma=EU=E(T+vN+W),\qquad
E=-\g(\dot\gamma,T)>0,
\]
where $W$ is tangent to $S(t,u)$. Normalizing the
timelike case by proper time, we have
\[
v^2+|W|^2=1-\frac{m^2}{E^2},\qquad
m^2\in\{0,1\},
\]
with $m=0$ in the null case. Moreover,
\[
\frac{dt}{d\lambda}=\frac E\phi,\qquad
\frac{du}{dt}=\frac{\phi}{a}(1-v)\geq 0.
\]
Suppose that $\gamma$ has finite future affine length.
Since $t\to\infty$ along $\gamma$, we have
\begin{equation}\label{conclusion:finite-affine}
\int^\infty\frac{\phi}{E}\,dt<\infty.
\end{equation}
We first show that $r(t)\gtrsim t$ for large $t$.
Since $\D_{\dot\gamma}\dot\gamma=0$,
differentiating $Ev=\g(\dot\gamma,N)$ and
$E=-\g(\dot\gamma,T)$ gives
\[
\frac1\phi\frac{dv}{dt}
=\g(U,\D_Ue_4)-(1-v)\g(U,\D_UT).
\]
Using the null frame identities and $\xi=0$, we have
\begin{align*}
\g(U,\D_Ue_4)
={}&\chi(W,W)+(1-v)\bigl(
(1+v)\om-(1-v)\omb
+(\eta+\ze)\cdot W\bigr),\\
\g(U,\D_UT)
={}&v\nab_N\log\phi-v^2\de
+\nabs\log\phi\cdot W
-2v\eps\cdot W-\ka(W,W).
\end{align*}
Consequently, \eqref{6.6} yields
\begin{align*}
\frac1\phi\frac{dv}{dt}
=\chi(W,W)+(1-v)\Big[
(1+v)\big(v\de-\nab_N\log\phi
+2\eps\cdot W\big)
-\xib\cdot W+\ka(W,W)
\Big].
\end{align*}
Since $W(r)=0$ and $T(r)=e_4r-N(r)$,
we obtain
\begin{align*}
\frac1\phi\frac{d}{dt}(rv)
={}&v e_4r-v(1-v)N(r)
+\frac r\phi\frac{dv}{dt}\\
={}&v^2+|W|^2+v(e_4r-1)
-v(1-v)(N(r)-1)
+r\Gag+r(1-v)\Gab.
\end{align*}
Using $r|\Gag|\les\ep_0$, we conclude that
\begin{equation}\label{conclusion:radial-geodesic}
\frac1\phi\frac{d}{dt}(rv)
\geq 1-\frac{m^2}{E^2}-C\ep_0
-C(1-v)\left(
|N(r)-1|+r|\Gab|
\right).
\end{equation}
The factor $1-v$ allows us to use
$du=\phi a^{-1}(1-v)\,dt$.
By \eqref{wonderfulrp:coefficient-size} and
\eqref{Gab-quadratic-flux-bulk},
\[
\int du\,\sup_t\left(
r\|a^{1/2}\Gab\|_{L^\infty(S)}^2
+\|a^{1/2}\eta\|_{L^4(S)}^2
\right)\les\ep_0^2.
\]
Since
\[
\frac{dr}{dt}
=\phi\bigl(e_4r-(1-v)N(r)\bigr)\leq C,
\]
we have $r(t)\les1+t$. In particular,
\[
\int_{u(t_0)}^{u(t)}a\,du\les1+t,
\qquad
\int_{u(t_0)}^{u(t)}ar\,du\les(1+t)^2.
\]
Using \eqref{normal:radius-defect} and
Cauchy--Schwarz, with the $u$-integrals taken
over $[u(t_0),u(t)]$, we obtain
\begin{align*}
\int a|N(r)-1|\,du
&\les\ep_0\int a\,du
+\left(\int ar\,du\right)^{1/2}
 \left(\int du\,\sup_t
 \|a^{1/2}\eta\|_{L^4(S)}^2\right)^{1/2}\\
&\quad+\sup_{[t_0,t]}r
 \int du\,\sup_t
 \|a^{1/2}\eta\|_{L^4(S)}^2
\les\ep_0(1+t),\\
\int ar|\Gab|\,du
&\les
\left(\int du\,\sup_t
r\|a^{1/2}\Gab\|_{L^\infty(S)}^2\right)^{1/2}
\left(\int ar\,du\right)^{1/2}
\les\ep_0(1+t).
\end{align*}
Therefore,
\[
\int_{t_0}^t\phi(1-v)
\left(|N(r)-1|+r|\Gab|\right)\,ds
\les\ep_0(1+t).
\]
Integrating \eqref{conclusion:radial-geodesic}
and using $\phi\simeq1$, we obtain
\[
rv(t)\geq(1-C\ep_0)t-C_\gamma
-m^2\int_{t_0}^t\frac{\phi}{E^2}\,ds.
\]
The last term vanishes for $m=0$ and is bounded
by \eqref{conclusion:finite-affine} for $m=1$,
since $E\geq1$. As $|v|\leq1$, this gives
\begin{equation}\label{conclusion:radial-escape}
r(t)\gtrsim t
\qquad\text{for all sufficiently large }t.
\end{equation}
It remains to estimate $E$. The geodesic
equation and $E=-\g(\dot\gamma,T)$ give
\[
\frac{dE}{d\lambda}
=-E^2\g(U,\D_UT).
\]
Using $dt/d\lambda=E/\phi$,
$\g(T,\D_XT)=0$,
$\D_TT=\nab\log\phi$, and
$\g(vN+W,\D_{vN+W}T)
=-k(vN+W,vN+W)$, we obtain
\[
\frac{d}{dt}\log E
=\phi\left(
k(vN+W,vN+W)
-\nab\log\phi\cdot(vN+W)
\right).
\]
By \eqref{6.6} and the good coefficient bounds,
\[
\left|\frac{d}{dt}\log E\right|\les\phi\left(\frac{\ep_0}{r}+|\hchb|\,|W|^2
\right).
\]
Now take $t_0$ sufficiently large so that
\eqref{conclusion:radial-escape} holds for $t\geq t_0$.
Since $|W|^2\leq2(1-v)$, Cauchy--Schwarz
and the weighted bad-coefficient estimate yield
\begin{align*}
\int_{t_0}^t\phi|\hchb||W|^2\,ds\les\int_{u(t_0)}^{u(t)}a|\hchb|\,du\les
\left(\int ar|\hchb|^2\,du\right)^{1/2}
\left(\int\frac a r\,du\right)^{1/2}\les\ep_0\sqrt{\log(2+t)},
\end{align*}
where the last step follows from
\[
\int_{u(t_0)}^{u(t)}\frac a r\,du
=\int_{t_0}^t\frac{\phi(1-v)}{r(s)}\,ds
\les\int_{t_0}^t\frac{ds}{s}
\les\log(2+t).
\]
Consequently,
\[
\log\frac{E(t)}{E(t_0)}
\les\ep_0\log(2+t),
\qquad
E(t)\les C_\gamma(1+t)^{C\ep_0}.
\]
Therefore,
\[
\lambda(t)-\lambda(t_0)
=\int_{t_0}^t\frac{\phi}{E}\,ds
\gtrsim_\gamma
\int_{t_0}^t(1+s)^{-C\ep_0}\,ds
\longrightarrow\infty,
\]
provided $\ep_0$ is sufficiently small.
This contradicts \eqref{conclusion:finite-affine}
and completes the proof.
\end{proof}
\subsection{Asymptotically round spheres}
\label{sec:round-spheres}
\begin{prop}\label{conclusion:round-spheres}
Under the conclusions of Theorem~\ref{mainth}, for every fixed $u$, the Gauss curvature of $S(t,u)$ satisfies
\begin{equation}\label{conclusion:round-limit}
 \lim_{t\to+\infty}
 r^{3/2}\left\|\K-\frac1{r^2}\right\|_{L^4(S(t,u))}=0.
\end{equation}
\end{prop}
\begin{proof}
We use the established curvature and connection estimates,
in particular $r|\Gag|\les\ep_0$.
Applying the spatial exterior trace inequalities
\eqref{normal:exterior-sobolev-trace} and
\eqref{normal:weighted-L4-trace}, with $p=0,2$,
and then integrating in $t$, we obtain
\begin{equation}\label{conclusion:maximal-cone-traces}
\sup_u\int_{C_u\cap\{r\geq1\}}\left(r^{-1}|(\ze,\etab)|^2+r|\nabs(\de,\etab)|^2\right)d\mu_{C_u}\les\ep_0^2.
\end{equation}
Here we use sphere H\"older, $\phi\simeq1$,
\eqref{6.6}, and the maximal spatial bulk estimates.

We next establish the outgoing flux of $\slu$.
By \eqref{transport:codazzi} and the sphere Hodge
estimate, absorbing $\hch\c\ze$, we have
\[
 \|\nabs\hch\|_{L^2(S)}
 +r^{-1}\|\hch\|_{L^2(S)}
 \les\|(\slu,\b)\|_{L^2(S)}.
\]
Multiplying \eqref{U-transport} by $r^2\slu$
and integrating along $C_u$, the weighted energy
identity and \eqref{conclusion:maximal-cone-traces}
give
\[
 \int_{C_u\cap\{r\geq1\}}r|\slu|^2d\mu_{C_u}
 \les\ep_0^2+
 \ep_0^2\int_{C_u\cap\{r\geq1\}}
 r|\nabs\hch|^2d\mu_{C_u}.
\]
Indeed, the principal source is
$-\trch\,\b$, controlled by the curvature flux.
The terms containing $\om\slu$ and
$\hch\c\slu$ are absorbed using
$r|\Gag|\les\ep_0$, while
$\nabs|\hch|^2=2\hch\c\nabs\hch$ gives the
last integral. The remaining terms and the
initial cut are controlled by the established
estimates. The preceding Hodge inequality
and absorption therefore yield
\begin{equation}\label{conclusion:U-cone-flux}
 \sup_u\int_{C_u\cap\{r\geq1\}}
 r|\slu|^2d\mu_{C_u}
 \les\ep_0^2.
\end{equation}
By the Gauss equation \eqref{gauss}, using $\trchb=2\de-\trch$ and
$\nabs\trch=\slu-\trch\,\ze$, we have
\[
 |r\nabs\K|
 \les |r\nabs\rhoc|+|\slu|
 +|\nabs\de|+r^{-1}|\ze|.
\]
Thus the curvature flux and
\eqref{conclusion:maximal-cone-traces}--%
\eqref{conclusion:U-cone-flux} imply
\begin{equation}\label{conclusion:K-cone-flux}
 \sup_u\int_{C_u\cap\{r\geq1\}}
 r|r\nabs\K|^2d\mu_{C_u}
 \les\ep_0^2.
\end{equation}
Next, the Gauss equation, the outgoing null
structure equations, and \eqref{bianchirenor}
give schematically
\begin{align*}
 \nabs_4\left(\K-\frac1{r^2}\right)
 +\trch\left(\K-\frac1{r^2}\right)
 =-\sdivs\b-r^{-1}\sdivs\etab+\Gag\c\W+\Gag\c\nabs\Gag
 +r^{-1}\Gag\c\Gag+r^{-2}\Gag.
\end{align*}
The principal terms and nonlinear products are
controlled by the curvature flux,
\eqref{conclusion:maximal-cone-traces},
and $r|\Gag|\les\ep_0$.
For the last term, the radius identity \eqref{e3e4r}
and sphere Poincar\'e control
$e_4r-\frac r2\trch$ by
$r^2(\slu+r^{-1}\etab)$, so its cone flux
follows from \eqref{conclusion:maximal-cone-traces}
and \eqref{conclusion:U-cone-flux}.
Moreover, Gauss--Bonnet and sphere Poincar\'e give
\[
 \int_S r\left|\K-\frac1{r^2}\right|^2d\slg
 \les\int_S r|r\nabs\K|^2d\slg.
\]
Consequently, the preceding cone fluxes yield
\begin{equation}\label{conclusion:K-outgoing}
 \sup_u\int_{C_u\cap\{r\geq1\}}
 r\left|r\nabs_4\left(\K-\frac1{r^2}\right)\right|^2
 d\mu_{C_u}
 \les\ep_0^2.
\end{equation}
Finally, fix $u$. For sufficiently large $\Lambda$,
the initial cut term in \eqref{bootstrap:cone-trace}
is absent. Applying the cone trace estimate to
$r^{3/2}(\K-r^{-2})$ and using sphere Poincar\'e,
we obtain
\[
 \sup_{\Lambda\leq r(t,u)\leq2\Lambda}
 r^3\left\|\K-\frac1{r^2}\right\|_{L^4(S(t,u))}^2
 \les
 \int_{C_u\cap\{r\geq\Lambda/2\}}
 r\left(
 |r\nabs\K|^2
 +\left|r\nabs_4\left(\K-\frac1{r^2}\right)\right|^2
 \right)d\mu_{C_u}
 \longrightarrow0,
\]
by \eqref{conclusion:K-cone-flux} and \eqref{conclusion:K-outgoing}. Since $r(t,u)\to\infty$ as $t\to\infty$, this proves \eqref{conclusion:round-limit}.
\end{proof}

\end{document}